# First Efficient Convergence for Streaming k-PCA: a Global, Gap-Free, and Near-Optimal Rate*

Zeyuan Allen-Zhu
zeyuan@csail.mit.edu
Institute for Advanced Study

Yuanzhi Li
yuanzhil@cs.princeton.edu
Princeton University

July 26, 2016†

**Abstract**

We study streaming principal component analysis (PCA), that is to find, in $O(dk)$ space, the top $k$ eigenvectors of a $d \times d$ hidden matrix $\mathbf{\Sigma}$ with online vectors drawn from covariance matrix $\mathbf{\Sigma}$.

We provide *global* convergence for Oja's algorithm which is popularly used in practice but lacks theoretical understanding for $k > 1$. We also provide a modified variant $\mathsf{Oja}^{++}$ that runs *even faster* than Oja's. Our results match the information theoretic lower bound in terms of dependency on error, on eigengap, on rank $k$, and on dimension $d$, up to poly-log factors. In addition, our convergence rate can be made gap-free, that is proportional to the approximation error and independent of the eigengap.

In contrast, for general rank $k$, before our work (1) it was open to design any algorithm with efficient global convergence rate; and (2) it was open to design any algorithm with (even local) gap-free convergence rate in $O(dk)$ space.

## 1 Introduction

Principle component analysis (PCA) is the problem of finding the subspace of largest variance in a dataset consisting of vectors, and is a fundamental tool used to analyze and visualize data in machine learning, computer vision, statistics, and operations research. In the big-data scenario, since it can be unrealistic to store the entire dataset, it is interesting and more challenging to study the streaming model (a.k.a. the stochastic online model) of PCA.

Suppose the data vectors $x \in \mathbb{R}^d$ are drawn i.i.d. from an unknown distribution with covariance matrix $\mathbf{\Sigma} = \mathbb{E}[xx^\top] \in \mathbb{R}^{d \times d}$, and the vectors are presented to the algorithm in an online fashion. Following [10, 12], we assume the Euclidean norm $\|x\|_2 \leq 1$ with probability 1 (therefore $\mathbf{Tr}(\mathbf{\Sigma}) \leq 1$) and we are interested in approximately computing the top $k$ eigenvectors of $\mathbf{\Sigma}$. We are interested in algorithms with memory storage $O(dk)$, the same as the memory needed to store any $k$ vectors in $d$ dimensions. We call this the *streaming k-PCA problem*.

For streaming $k$-PCA, the popular and natural extension of Oja's algorithm originally designed for the $k = 1$ case works as follows. Beginning with a random Gaussian matrix $\mathbf{Q}_0 \in \mathbb{R}^{d \times k}$ (each

---

*We thank Jieming Mao for discussing our lower bound Theorem 6, and thank Dan Garber and Elad Hazan for useful conversations. Z. Allen-Zhu is partially supported by a Microsoft research award, no. 0518584, and an NSF grant, no. CCF-1412958.

†An earlier version of this paper appeared at https://arxiv.org/abs/1607.07837. This newer version contains a stronger Theorem 2, a new lower bound Theorem 6, as well as the new $\mathsf{Oja}^{++}$ results Theorem 4 and Theorem 5.

entry i.i.d $\sim \mathcal{N}(0,1)$), it repeatedly applies

$$\text{rank-}k \text{ Oja's algorithm:} \quad \mathbf{Q}_t \leftarrow (\mathbf{I} + \eta_t x_t x_t^\top)\mathbf{Q}_{t-1}, \quad \mathbf{Q}_t = \mathsf{QR}(\mathbf{Q}_t) \quad (1.1)$$

where $\eta_t > 0$ is some learning rate that may depend on $t$, vector $x_t$ is the random vector in iteration $t$, and $\mathsf{QR}(\mathbf{Q}_t)$ is the Gram-Schmidt decomposition that orthonormalizes the columns of $\mathbf{Q}_t$.

Although Oja's algorithm works reasonably well in practice, very limited theoretical results are known for its convergence in the $k > 1$ case. Even worse, little is known for *any* algorithm that solves streaming PCA in the $k > 1$. Specifically, there are three major challenges for this problem:

1. Provide an *efficient* convergence rate that only logarithmically depends on the dimension $d$.
2. Provide a *gap-free* convergence rate that is independent of the eigengap.
3. Provide a *global* convergence rate so the algorithm can start from a random initial point.

In the case of $k > 1$, to the best of our knowledge, only Shamir [18] successfully analyzed the original Oja's algorithm. His convergence result is only local and not gap-free.[1]

Other groups of researchers [3, 9, 12] studied a *block variant* of Oja's, that is to sample multiple vectors $x$ in each round $t$, and then use their empirical covariance to replace the use of $x_t x_t^\top$. This algorithm is more stable and easier to analyze, but only leads to suboptimal convergence.

We discuss them more formally below (and see Table 1):

- Shamir [18] implicitly provided a *local* but efficient convergence result for Oja's algorithm,[2] which requires a very accurate starting matrix $\mathbf{Q}_0$: his theorem relies on $\mathbf{Q}_0$ being correlated with the top $k$ eigenvectors by a correlation value at least $k-1/2$. If using random initialization, this event happens with probability at most $2^{-\Omega(d)}$.

- Hardt and Price [9] analyzed the block variant of Oja's,[3] and obtained a global convergence that linearly scales with the dimension $d$. Their result also has a cubic dependency on the gap between the $k$-th and $(k+1)$-th eigenvalue which is not optimal. They raised an open question regarding how to provide any convergence result that is gap-free.

- Balcan et al. [3] analyzed the block variant of Oja's. Their results are also not efficient and cubically scale with the eigengap. In the gap-free setting, their algorithm runs in space more than $O(kd)$, and also outputs more than $k$ vectors.[4] For such reason, we do not include their gap-free result in Table 1, and shall discuss it more in Section 4.

- Li et al. [12] also analyzed the block variant of Oja's. Their result also cubically scales with the eigengap, and their global convergence is not efficient.

- In practice, researchers observed that it is advantageous to choose the learning rate $\eta_t$ to be high at the beginning, and then gradually decreasing (c.f. [22]). To the best of our knowledge, there is no theoretical support behind this learning rate scheme for general $k$.

In sum, it remains open before our work to obtain (1) any gap-free convergence rate in space $O(kd)$, (2) any global convergence rate that is efficient, or (3) any global convergence rate that has the optimal quadratic dependence on eigengap.

---

[1] A local convergence rate means that the algorithm needs a warm start that is sufficiently close to the solution. However, the complexity to reach such a warm start is not clear.

[2] The original method of Shamir [18] is an offline one. One can translate his result into a streaming setting and this requires a lot of extra work including the martingale techniques we introduce in this paper.

[3] They are in fact only able to output $2k$ vectors, guaranteed to approximately include the top $k$ eigenvectors.

[4] They require space $O((k+m)d)$ where $k+m$ is the number of eigenvalues in the interval $[\lambda_k - \rho, 1]$ for some "virtual gap" parameter $\rho$. See our Theorem 2 for a definition. This may be as large as $O(d^2)$. Also, they output $k+m$ vectors which are only guaranteed to approximately "contain" the top $k$ eigenvectors.

| | Paper | Global Convergence | Is It "Efficient"? | Local Convergence |
|---|---|---|---|---|
| $k=1$ gap-dependent | Shamir [17] | $\widetilde{O}\big(\frac{d}{\mathsf{gap}^2}\cdot\frac{1}{\varepsilon}\big)$ ♭ | no | $\widetilde{O}\big(\frac{1}{\mathsf{gap}^2}\cdot\frac{1}{\varepsilon}\big)$ ♭ |
| | Sa et al. [16] | $\widetilde{O}\big(\frac{d}{\mathsf{gap}^2}\cdot\frac{1}{\varepsilon}\big)$ ♭ | no | $\widetilde{O}\big(\frac{d}{\mathsf{gap}^2}\cdot\frac{1}{\varepsilon}\big)$ ♭ |
| | Li et al. [11] [a] | $\widetilde{O}\big(\frac{d\lambda_1}{\mathsf{gap}^2}\cdot\frac{1}{\varepsilon}\big)$ ♭ | no | $\widetilde{O}\big(\frac{d\lambda_1}{\mathsf{gap}^2}\cdot\frac{1}{\varepsilon}\big)$ ♭ |
| | Jain et al. [10] | $\widetilde{O}\big(\frac{\lambda_1}{\mathsf{gap}^2}\cdot\frac{1}{\varepsilon}\big)$ | yes | $\widetilde{O}\big(\frac{\lambda_1}{\mathsf{gap}^2}\cdot\frac{1}{\varepsilon}\big)$ |
| | **Theorem 1 (Oja)** | $\widetilde{O}\big(\frac{\lambda_1}{\mathsf{gap}^2}\cdot\frac{1}{\varepsilon}\big)$ | yes | $\widetilde{O}\big(\frac{\lambda_1}{\mathsf{gap}^2}\cdot\frac{1}{\varepsilon}\big)$ |
| $k=1$ gap-free | Shamir [17] (Remark 1.3) | $\widetilde{O}\big(\frac{d}{\rho^2}\cdot\frac{1}{\varepsilon^2}\big)$ ♭ | no | $\widetilde{O}\big(\frac{1}{\rho^2}\cdot\frac{1}{\varepsilon^2}\big)$ ♭ |
| | **Theorem 2 (Oja)** | $\widetilde{O}\big(\frac{\lambda_{1\sim(1+m)}}{\rho^2}\cdot\frac{1}{\varepsilon}\big)$ | yes | $\widetilde{O}\big(\frac{\lambda_{1\sim(1+m)}}{\rho^2}\cdot\frac{1}{\varepsilon}\big)$ |
| $k\geq 1$ gap-dependent | Hardt-Price [9] [b] | $\widetilde{O}\big(\frac{d\lambda_k}{\mathsf{gap}^3}\cdot\frac{1}{\varepsilon}\big)$ ♭ | no | $\widetilde{O}\big(\frac{d\lambda_k}{\mathsf{gap}^3}\cdot\frac{1}{\varepsilon}\big)$ ♭ |
| | Li et al. [12] [b] | $\widetilde{O}\big(\frac{k\lambda_k}{\mathsf{gap}^3}\cdot\big(kd+\frac{1}{\varepsilon}\big)\big)$ ♭ | no | $\widetilde{O}\big(\frac{k\lambda_k}{\mathsf{gap}^3}\cdot\frac{1}{\varepsilon}\big)$ ♭ |
| | Shamir [18] | unknown ♭ | no | $O\big(\frac{1}{\mathsf{gap}^2}\cdot\frac{1}{\varepsilon}\big)$ ♭ |
| | Balcan et al. [3] [b] | $\widetilde{O}\big(\frac{d(\lambda_{1\sim k})^2\lambda_k}{\mathsf{gap}^3}\cdot\frac{1}{\varepsilon}\big)$ ♭ (when $\lambda_{1\sim k}\geq k/d$) [c] | no | $\widetilde{O}\big(\frac{d(\lambda_{1\sim k})^2\lambda_k}{\mathsf{gap}^3}\cdot\frac{1}{\varepsilon}\big)$ ♭ (when $\lambda_{1\sim k}\geq k/d$) |
| | **Theorem 1 (Oja)** | $\widetilde{O}\big(\frac{\lambda_{1\sim k}}{\mathsf{gap}^2}\cdot\big(\frac{1}{\varepsilon}+k\big)\big)$ | yes | $\widetilde{O}\big(\frac{\lambda_{1\sim k}}{\mathsf{gap}^2}\cdot\frac{1}{\varepsilon}\big)$ |
| | **Theorem 4 (Oja$^{++}$)** | $\widetilde{O}\big(\frac{\lambda_{1\sim k}}{\mathsf{gap}^2}\cdot\frac{1}{\varepsilon}\big)$ | yes | $\widetilde{O}\big(\frac{\lambda_{1\sim k}}{\mathsf{gap}^2}\cdot\frac{1}{\varepsilon}\big)$ |
| | **Theorem 6 (LB)** | $\Omega\big(\frac{k\lambda_k}{\mathsf{gap}^2}\cdot\frac{1}{\varepsilon}\big)$ (lower bound) | | |
| $k\geq 1$ gap-free | **Theorem 2 (Oja)** | $\widetilde{O}\big(\frac{\min\{1,(\lambda_{1\sim k}+k\cdot\lambda_{(k+1)\sim(k+m)})\}}{\rho^2}\cdot k\big)$ $+\widetilde{O}\big(\frac{\lambda_{1\sim k+m}}{\rho^2}\cdot\frac{1}{\varepsilon}\big)$ | yes | $\widetilde{O}\big(\frac{\lambda_{1\sim k+m}}{\rho^2}\cdot\frac{1}{\varepsilon}\big)$ |
| | **Theorem 5 (Oja$^{++}$)** | $\widetilde{O}\big(\frac{\lambda_{1\sim k+m}}{\rho^2}\cdot\frac{1}{\varepsilon}\big)$ | yes | $\widetilde{O}\big(\frac{\lambda_{1\sim k+m}}{\rho^2}\cdot\frac{1}{\varepsilon}\big)$ |
| | **Theorem 6 (LB)** | $\Omega\big(\frac{k\lambda_k}{\rho^2}\cdot\frac{1}{\varepsilon}\big)$ (lower bound) | | |

Table 1: Comparison of known results. For $\mathsf{gap}\stackrel{\text{def}}{=}\lambda_k-\lambda_{k+1}$, every $\varepsilon\in(0,1)$ and $\rho\in(0,1)$:

- "gap-dependent convergence" means $\|\mathbf{Q}_T^\top\mathbf{Z}\|_F^2\leq\varepsilon$ where $\mathbf{Z}$ consists of the last $d-k$ eigenvectors.
- "gap-free convergence" means $\|\mathbf{Q}_T^\top\mathbf{W}\|_F^2\leq\varepsilon$ where $\mathbf{W}$ consists of all eigenvectors with eigenvalues $\leq\lambda_k-\rho$.
- a global convergence is "efficient" if it only (poly-)logarithmically depend on the dimension $d$.
- $k$ is the target rank; in gap-free settings $m$ be the largest index so that $\lambda_{k+m}>\lambda_k-\rho$.
- we denote by $\lambda_{a\sim b}\stackrel{\text{def}}{=}\sum_{i=a}^b\lambda_i$ in this table. Since $\|x\|_2\leq 1$ for each sample vector, we have
  $$\mathsf{gap}\in[0,1/k],\quad \lambda_i\in[0,1/i],\quad k\mathsf{gap}\leq k\lambda_k\leq\lambda_{1\sim k}\leq\lambda_{1\sim k+m}\leq 1\ .$$
- we use ♭ to indicate the result is outperformed.
- some results in this table (both ours and prior work) depend on $\lambda_{1\sim k}$. In principle, this requires the algorithm to know a constant approximation of $\lambda_{1\sim k}$ upfront. In practice, however, since one always tunes the learning rate $\eta$ (for any algorithm in the table), we do not need additional knowledge on $\lambda_{1\sim k}$.

---

[a]The result of [11] is in fact $\widetilde{O}\big(\frac{d\lambda_1^2}{\mathsf{gap}^2}\cdot\frac{1}{\varepsilon}\big)$ by under a stronger 4-th moment assumption. It slows down at least by a factor $1/\lambda_1$ if the 4-th moment assumption is removed.

[b]These results give guarantees on spectral norm $\|\mathbf{Q}_T^\top\mathbf{W}\|_2^2$, so we increased them by a factor $k$ for a fair comparison.

[c]If $\|x_t\|_2$ is always 1 then $\lambda_{1\sim k}\geq k/d$ always holds. Otherwise, even in the rare case of $\lambda_{1\sim k}<k/d$, their complexity becomes $\widetilde{O}\big(\frac{k^2\lambda_k}{d\cdot\mathsf{gap}^3}\big)$ and is still worse than ours.

**Over Sampling.** Let us emphasize that it is often desirable to directly output a $d \times k$ matrix $\mathbf{Q}_T$. Some of the previous results, such as Hardt and Price [9], or the gap-free case of Balcan et al. [3], are only capable of finding an over-sampled matrix $d \times k'$ for some $k' > k$, with the guarantee that these $k'$ columns approximately contain the top $k$ eigenvectors of $\mathbf{\Sigma}$. However, it is not clear how to find "the best $k$ vectors" out of this $k'$-dimensional subspace.

**Special Case of $k = 1$.** Jain [10] obtained the first convergence result that is both efficient and global for streaming 1-PCA. Shamir [17] obtained the first gap-free result for streaming 1-PCA, but his result is not efficient. Both these results are based on Oja's algorithm, and it remains open before our work to obtain a gap-free result that is also efficient even when $k = 1$.

**Other Related Results.** Mitliagkas et al. [13] obtained a streaming PCA result but in the restricted spiked covariance model. Balsubramani et al. [4] analyzed a modified variant of Oja's algorithm and needed an extra $O(d^5)$ factor in the complexity.

The offline problem of PCA (and SVD) can be solved via iterative algorithms that are based on variance-reduction techniques on top of stochastic gradient methods [2, 18] (see also [6, 7] for the $k = 1$ case); these methods do multiple passes on the input data so are not relevant to the streaming model. Offline PCA can also be solved via power method or block Krylov method [14], but since each iteration of these methods relies on one full pass on the dataset, they are not suitable for streaming setting either. Other offline problems and efficient algorithms relevant to PCA include canonical correlation analysis and generalized eigenvector decomposition [1, 8, 21].

Offline PCA is *significantly easier* to solve because one can (although non-trivially) reduce a general $k$-PCA problem to $k$ times of 1-PCA using the techniques of [2]. However, this is *not the case* in streaming PCA because one can lose a large polynomial factor in the sampling complexity.

## 1.1 Results on Oja's Algorithm

We denote by $\lambda_1 \geq \cdots \geq \lambda_d \geq 0$ the eigenvalues of $\mathbf{\Sigma}$, and it satisfies $\lambda_1 + \cdots + \lambda_d = \mathbf{Tr}(\mathbf{\Sigma}) \leq 1$. We present convergence results on Oja's algorithm that are *global, efficient and gap-free*.

Our first theorem works when there is a eigengap between $\lambda_k$ and $\lambda_{k+1}$:

**Theorem 1** (Oja, gap-dependent). *Letting* $\mathsf{gap} \stackrel{\text{def}}{=} \lambda_k - \lambda_{k+1} \in \left(0, \frac{1}{k}\right]$ *and* $\Lambda \stackrel{\text{def}}{=} \sum_{i=1}^{k} \lambda_i \in (0, 1]$, *for every* $\varepsilon, p \in (0, 1)$ *define learning rates*

$$T_0 = \widetilde{\Theta}\left(\frac{k\Lambda}{\mathsf{gap}^2 p^2}\right), \quad T_1 = \widetilde{\Theta}\left(\frac{\Lambda}{\mathsf{gap}^2}\right), \quad \eta_t = \begin{cases} \widetilde{\Theta}\left(\frac{1}{\mathsf{gap} \cdot T_0}\right) & 1 \leq t \leq T_0; \\ \widetilde{\Theta}\left(\frac{1}{\mathsf{gap} \cdot T_1}\right) & T_0 < t \leq T_0 + T_1; \\ \widetilde{\Theta}\left(\frac{1}{\mathsf{gap} \cdot (t - T_0)}\right) & t > T_0 + T_1. \end{cases}$$ [5]

*Let $\mathbf{Z}$ be the column orthonormal matrix consisting of all eigenvectors of $\mathbf{\Sigma}$ with values no more than $\lambda_{k+1}$. Then, the output $\mathbf{Q}_T \in \mathbb{R}^{d \times k}$ of Oja's algorithm satisfies with probability at least $1-p$:*

*for every*[6] $\quad T = T_0 + T_1 + \widetilde{\Theta}\left(\frac{T_1}{\varepsilon}\right) \quad$ *it satisfies* $\quad \|\mathbf{Z}^\top \mathbf{Q}_T\|_F^2 \leq \varepsilon$ .

*Above, $\widetilde{\Theta}$ hides poly-log factors in $\frac{1}{p}, \frac{1}{\mathsf{gap}}$ and $d$.*

In other words, after a warm up phase of length $T_0$, we obtain a $\frac{\lambda_1 + \cdots + \lambda_k}{\mathsf{gap}^2} \cdot \frac{1}{T}$ convergence rate for the quantity $\|\mathbf{Z}^\top \mathbf{Q}_T\|_F^2$. We make several observations (see also Table 1):

- In the $k = 1$ case, Theorem 1 matches the best known result of Jain et al. [10].

---
[5] The intermediate stage $[T_0, T_0 + T_1]$ is in fact unnecessary, but we add this phase to simplify proofs.
[6] Theorem also applies to every $T \geq T_0 + T_1 + \widetilde{\Omega}(T_1/\varepsilon)$ by making $\eta_t$ poly-logarithmically dependent on $T$.

- In the $k > 1$ case, Theorem 1 gives the first efficient global convergence rate.
- In the $k > 1$ case, even in terms of local convergence rate, Theorem 1 is faster than the best known result of Shamir [18] by a factor $\lambda_1 + \cdots + \lambda_k \in (0, 1)$.

*Remark* 1.1. The quantity $\|\mathbf{Z}^\top \mathbf{Q}_T\|_F^2$ captures the correlation between the resulting matrix $\mathbf{Q}_T \in \mathbb{R}^{d \times k}$ and the smallest $d - k$ eigenvectors of $\mathbf{\Sigma}$. It is a natural generalization of the sin-square quantity widely used in the $k = 1$ case, because if $k = 1$ then $\|\mathbf{Z}^\top \mathbf{Q}_T\|_F^2 = \sin^2(q, \nu_1)$ where $q$ is the only column of $\mathbf{Q}$ and $\nu_1$ is the leading eigenvector of $\mathbf{\Sigma}$.

Some literatures instead adopt the spectral-norm guarantee (i.e., bounds on $\|\mathbf{Z}^\top \mathbf{Q}_T\|_2^2$) as opposed to the Frobenius-norm one. The two guarantees are only up to a factor $k$ away. We choose to prove Frobenius-norm results because: (1) it makes the analysis significantly simpler, and (2) $k$ is usually small comparing to $d$, so if one can design an efficient (i.e., dimension free) convergence rate for the Frobenius norm that also implies an efficient convergence rate for the spectral norm.

*Remark* 1.2. Our lower bound later (i.e. Theorem 6) implies, at least when $\lambda_1$ and $\lambda_k$ are within a constant factor of each other, the local convergence rate in Theorem 1 is optimal up to log factors.

**Gap-Free Streaming $k$-PCA.** When the eigengap is small which is usually true in practice, it is desirable to obtain *gap-free* convergence [14, 17]. We have the following theorem which answers the open question of Hardt and Price [9] regarding gap-free convergence rate for streaming $k$-PCA.

> **Theorem 2** (Oja, gap-free)**.** *For every $\rho, \varepsilon, p \in (0, 1)$, let $\lambda_1, \ldots, \lambda_m$ be all eigenvalues of $\mathbf{\Sigma}$ that are $> \lambda_k - \rho$, let $\Lambda_1 \stackrel{\text{def}}{=} \sum_{i=1}^k \lambda_i \in (0, 1]$, $\Lambda_2 \stackrel{\text{def}}{=} \sum_{j=k+1}^{k+m} \lambda_j \in (0, 1]$, define learning rates*
> $$T_0 = \widetilde{\Theta}\left(\frac{k \cdot \min\{1, \Lambda_1 + \frac{k\Lambda_2}{p^2}\}}{\rho^2 \cdot p^2}\right), \quad T_1 = \widetilde{\Theta}\left(\frac{\Lambda_1 + \Lambda_2}{\rho^2}\right), \quad \eta_t = \begin{cases} \widetilde{\Theta}\left(\frac{1}{\rho \cdot T_0}\right) & t \leq T_0; \\ \widetilde{\Theta}\left(\frac{1}{\rho \cdot T_1}\right) & t \in (T_0, T_0 + T_1]; \\ \widetilde{\Theta}\left(\frac{1}{\rho \cdot (t - T_0)}\right) & t > T_0 + T_1. \end{cases}$$
> *Let $\mathbf{W}$ be the column orthonormal matrix consisting of all eigenvectors of $\mathbf{\Sigma}$ with values no more than $\lambda_k - \rho$. Then, the output $\mathbf{Q}_T \in \mathbb{R}^{d \times k}$ of Oja's algorithm satisfies with prob. at least $1 - p$:*
> $$\text{for every}^7 \quad T = T_0 + T_1 + \widetilde{\Theta}\left(\frac{T_1}{\varepsilon}\right) \quad \text{it satisfies} \quad \|\mathbf{W}^\top \mathbf{Q}_T\|_F^2 \leq \varepsilon \ .$$
> *Above, $\widetilde{\Theta}$ hides poly-log factors in $\frac{1}{p}, \frac{1}{\rho}$ and $d$.*

Note that the above theorem is a *double approximation*. The number of iterations depend both on $\rho$ and $\varepsilon$, where $\varepsilon$ is an upper bound on the correlation between $\mathbf{Q}_T$ and all eigenvectors in $\mathbf{W}$ (which depends on $\rho$). This is the first known gap-free result for the $k > 1$ case using $O(kd)$ space.

One may also be interested in single-approximation guarantees, such as the rayleigh-quotient guarantee. Note that a single-approximation guarantee by definition loses information about the $\varepsilon$-$\rho$ tradeoff; furthermore, (good) single-approximation guarantees are not easy to obtain.[8]

We show in this paper the following theorem regarding the rayleigh-quotient guarantee:

> **Theorem 3** (Oja, rayleigh quotient, informal)**.** *There exist learning rate choices so that, for every $T = \widetilde{\Theta}\left(\frac{k}{\rho^2 \cdot p^2}\right)$, letting $q_i$ be the $i$-th column of the output matrix $\mathbf{Q}_T$, then*
> $$\mathbf{Pr}\left[\forall i \in [k], \quad q_i^\top \mathbf{\Sigma} q_i \geq \lambda_i - \widetilde{\Theta}(\rho)\right] \geq 1 - p \ .$$
> *Again, $\widetilde{\Theta}$ hides poly-log factors in $\frac{1}{p}, \frac{1}{\rho}$ and $d$.*

---

[7]Theorem also applies to every $T \geq T_0 + \widetilde{\Omega}(T_0/\varepsilon)$ by making $\eta_t$ poly-logarithmically dependent on $T$.

[8]Pointed out by [10], a direct translation from double approximation to a rayleigh-quotient type of convergence loses a factor on the approximation error. They raised it as an open question regarding how to design a direct proof without sacrificing this loss. Our next theorem answers this open question (at least in the gap-free case).

*Remark* 1.3. Before our work, the only gap-free result with space $O(kd)$ is Shamir [17] — but it is not efficient and only for $k = 1$. His result is in Rayleigh quotient but not double-approximation. If the initialization phase is ignored, Shamir's local convergence rate matches our *global* one in Theorem 3. However, if one translates his result into double approximation, the running time loses a factor $\varepsilon$. This is why in Table 1 Shamir's result is in terms of $1/\varepsilon^2$ as opposed to $1/\varepsilon$.

## 1.2 Results on Our New Oja$^{++}$ Algorithm

Oja's algorithm has a slow initialization phase (which is also observed in practice [22]). For example, in the gap-dependent case, Oja's running time $\widetilde{O}\big(\frac{\lambda_1+\cdots+\lambda_k}{\rho^2}\cdot\big(k+\frac{1}{\varepsilon}\big)\big)$ is dominated by its initialization when $\varepsilon > 1/k$. We propose in this paper a modified variant of Oja's that initializes *gradually*.

**Our Oja$^{++}$ Algorithm.** At iteration 0, instead putting all the $dk$ random Gaussians into $\mathbf{Q}_0$ like Oja's, our Oja$^{++}$ only fills the first $k/2$ columns of $\mathbf{Q}_0$ with random Gaussians, and sets the remaining columns be zeros. It applies the same iterative rule as Oja's to go from $\mathbf{Q}_t$ to $\mathbf{Q}_{t+1}$, but after every $T_0$ iterations for some $T_0 \in \mathbb{N}^*$, it replaces the zeros in the next $k/4, k/8, \ldots$ columns with random Gaussians and continues.[9] This gradual initialization ends when all the $k$ columns become nonzero, and the remaining algorithm of Oja$^{++}$ works exactly the same as Oja's.

We provide pseudocode of Oja$^{++}$ in Algorithm 2 on page 58, and state below its main theorems:

**Theorem 4** (Oja$^{++}$, gap-dependent, informal)**.** *Letting* $\mathsf{gap} \stackrel{\text{def}}{=} \lambda_k - \lambda_{k+1} \in \big(0, \frac{1}{k}\big]$, *our* Oja$^{++}$ *outputs a column-orthonormal* $\mathbf{Q}_T \in \mathbb{R}^{d\times k}$ *with* $\|\mathbf{Z}^\top \mathbf{Q}_T\|_F^2 \leq \varepsilon$ *in* $T = \widetilde{\Theta}\left(\frac{\lambda_1+\cdots+\lambda_k}{\mathsf{gap}^2\varepsilon}\right)$ *iterations.*

**Theorem 5** (Oja$^{++}$, gap-free, informal)**.** *Given* $\rho \in (0,1)$, *our* Oja$^{++}$ *outputs a column-orthonormal* $\mathbf{Q}_T \in \mathbb{R}^{d\times k}$ *with* $\|\mathbf{W}^\top \mathbf{Q}_T\|_F^2 \leq \varepsilon$ *in* $T = \widetilde{\Theta}\left(\frac{\lambda_1+\cdots+\lambda_{k+m}}{\rho^2\varepsilon}\right)$ *iterations.*

## 1.3 Result on Lower Bound

We have the following information-theoretical lower bound for any (possibly offline) algorithm:

**Theorem 6** (lower bound, informal)**.** *For every integer $k \geq 1$, integer $m \geq 0$, every $0 < \rho < \lambda < 1/k$, every (possibly randomized) algorithm $\mathcal{A}$, we can construct a distribution $\mu$ over unit vectors with $\lambda_{k+m+1}(\mathbb{E}_\mu[xx^\top]) \leq \lambda - \rho$ and $\lambda_k(\mathbb{E}_\mu[xx^\top]) \geq \lambda$. The output $\mathbf{Q}_T$ of $\mathcal{A}$ with samples $x_1,\ldots,x_T$ i.i.d. drawn from $\mu$ satisfies*

$$\mathbb{E}_{x_1,\ldots,x_T,\mathcal{A}}\left[\|\mathbf{W}^\top \mathbf{Q}_T\|_F^2\right] = \Omega\left(\frac{k\lambda}{\rho^2 \cdot T}\right) \ .$$

(**W** *consists of the last $d - (k+m)$ eigenvectors of $\mathbb{E}_\mu[xx^\top]$.*)

Our Theorem 6 (with $m = 0$ and $\rho = \mathsf{gap}$) implies that, in the gap-dependent setting, the global convergence rate of Oja$^{++}$ is optimal up to log factors, at least when $\lambda_1 = O(\lambda_k)$. Our gap-free result does not match this lower bound. We explain in Section 4 that if one increases the space from $O(kd)$ to $O((k+m)d)$ in the gap-free case, our Oja$^{++}$ can also match this lower bound.

---

[9]Zeros columns will remain zero according to the usage of Gram-Schmidt in Oja's algorithm.

## 2 Preliminaries

We denote by $1 \geq \lambda_1 \geq \cdots \geq \lambda_d \geq 0$ the eigenvalues of the positive semidefinite (PSD) matrix $\boldsymbol{\Sigma}$, and $\nu_1, \nu_2, \ldots, \nu_d$ the corresponding normalized eigenvectors. Since we assumed $\|x\|_2 \leq 1$ for each data vector it satisfies $\lambda_1 + \cdots + \lambda_d = \mathbf{Tr}(\boldsymbol{\Sigma}) \leq 1$. We define $\mathsf{gap} \stackrel{\text{def}}{=} \lambda_k - \lambda_{k+1} \in \left[0, \frac{1}{k}\right]$. Slightly abusing notations, we also use $\lambda_k(\mathbf{M})$ to denote the $k$-th largest eigenvalue of an arbitrary $\mathbf{M}$.

Unless otherwise noted, we denote by $\mathbf{V} \stackrel{\text{def}}{=} [\nu_1, \ldots, \nu_k] \in \mathbb{R}^{d \times k}$ and $\mathbf{Z} \stackrel{\text{def}}{=} [\nu_{k+1}, \ldots, \nu_d] \in \mathbb{R}^{d \times (d-k)}$. For a given parameter $\rho > 0$ in our gap-free results, we also define $\mathbf{W} = [\nu_{k+m+1}, \ldots, \nu_d] \in \mathbb{R}^{d \times (d-k-m)}$ where $m$ is the largest index so that $\lambda_{k+m} > \lambda_k - \rho$. We write $\boldsymbol{\Sigma}_{\leq k} = \mathbf{V}\mathsf{Diag}\{\lambda_1, \ldots, \lambda_k\}\mathbf{V}^\top$ and $\boldsymbol{\Sigma}_{>k} \stackrel{\text{def}}{=} \mathbf{Z}\mathsf{Diag}\{\lambda_{k+1}, \ldots, \lambda_d\}\mathbf{Z}^\top$ so $\boldsymbol{\Sigma} = \boldsymbol{\Sigma}_{\leq k} + \boldsymbol{\Sigma}_{>k}$.

For a vector $y$, we sometimes denote by $y[i]$ or $y^{(i)}$ the $i$-th coordinate of $y$. We may use different notations in different lemmas in order to obtain the cleanest representations; when we do so, we shall clearly point out in the statement of the lemmas.

We denote by $\mathbf{P}_t \stackrel{\text{def}}{=} \prod_{s=1}^t (\mathbf{I} + \eta_s x_s x_s^\top)$ where $x_s$ is the $s$-th data sample and $\eta_s$ is the learning rate of iteration $s$. We denote by $\mathbf{Q} \in \mathbb{R}^{d \times k}$ (or $\mathbf{Q}_0$) the random initial matrix, and by $\mathbf{Q}_t \stackrel{\text{def}}{=} \mathsf{QR}((\mathbf{I} + \eta_t x_t x_t^\top)\mathbf{Q}_{t-1}) = \mathsf{QR}(\mathbf{P}_t \mathbf{Q}_0)$ the output of Oja's algorithm for every $t \geq 1$.[10] We use the notation $\mathcal{F}_t$ to denote the sigma-algebra generated by $x_t$. We denote $\mathcal{F}_{\leq t}$ to be the sigma-algebra generated by $x_1, \ldots, x_t$, i.e. $\mathcal{F}_{\leq t} = \vee_{s=1}^t \mathcal{F}_s$. In other words, whenever we condition on $\mathcal{F}_{\leq t}$ it means we have fixed $x_1, \ldots, x_t$.

For a vector $x$ we denote by $\|x\|$ or $\|x\|_2$ the Euclidean norm of $x$. We write $\mathbf{A} \succeq \mathbf{B}$ if $\mathbf{A}, \mathbf{B}$ are symmetric matrices and $\mathbf{A} - \mathbf{B}$ is PSD. We denote by $\|\mathbf{A}\|_{S_1}$ the Schatten-1 norm which is the summation of the (nonnegative) singular values of $\mathbf{A}$. It satisfies the following simple properties:

**Proposition 2.1.** *For not necessarily symmetric matrices $\mathbf{A}, \mathbf{B} \in \mathbb{R}^{d \times d}$ we have*

$$(1): \big|\mathbf{Tr}(\mathbf{A})\big| \leq \|\mathbf{A}\|_{S_1} \qquad (2): \big|\mathbf{Tr}(\mathbf{AB})\big| \leq \|\mathbf{AB}\|_{S_1} \leq \|\mathbf{A}\|_{S_1}\|\mathbf{B}\|_2 \ .$$

$$(3): \mathbf{Tr}(\mathbf{AB}) \leq \|A\|_F \|B\|_F = \big(\mathbf{Tr}(\mathbf{A}^\top \mathbf{A})\mathbf{Tr}(\mathbf{B}^\top \mathbf{B})\big)^{1/2} \ .$$

*Proof.* (1) is because $\mathbf{Tr}(\mathbf{A}) = \frac{1}{2}\mathbf{Tr}(\mathbf{A} + \mathbf{A}^\top) \leq \frac{1}{2}\|\mathbf{A} + \mathbf{A}^\top\|_{S_1} \leq \frac{1}{2}(\|\mathbf{A}\|_{S_1} + \|\mathbf{A}^\top\|_{S_1}) = \|\mathbf{A}\|_{S_1}$. (2) is because of (1) and the matrix Holder's inequality. (3) is owing to von Neumann's trace inequality (together with Cauchy's) which says $\mathbf{Tr}(\mathbf{AB}) \leq \sum_i \sigma_{A,i} \cdot \sigma_{B,i} \leq \|A\|_F \|B\|_F$. (Here, we have noted by $\sigma_{A,i}$ the $i$-th largest eigenvalue of $A$ and similarly for $B$. □

### 2.1 A Matrix View of Oja's Algorithm

The following lemma tells us that, for analysis purpose only, we can push the $\mathsf{QR}$ orthogonalization step in Oja's algorithm to the end:

**Lemma 2.2** (Oja's algorithm). *For every $s \in [d]$, every $\mathbf{X} \in \mathbf{R}^{d \times s}$, every $t \geq 1$, every initialization matrix $\mathbf{Q} \in \mathbb{R}^{d \times k}$, it satisfies $\|\mathbf{X}^\top \mathbf{Q}_t\|_F \leq \|\mathbf{X}^\top \mathbf{P}_t \mathbf{Q} (\mathbf{V}^\top \mathbf{P}_t \mathbf{Q})^{-1}\|_F$ .*

*Proof of Lemma 2.2.* Denoting by $\widetilde{\mathbf{Q}}_t = \mathbf{P}_t \mathbf{Q}$, we first observe that for every $t \geq 0$, $\mathbf{Q}_t = \widetilde{\mathbf{Q}}_t \mathbf{R}_t$ for some (upper triangular) invertible matrix $\mathbf{R}_t \in \mathbb{R}^{k \times k}$ (if $\mathbf{R}_t$ is not invertible, then the right hand side of the statement is $+\infty$ so we already done). The claim is true for $t = 0$. Suppose it holds for $t$ by induction, then

$$\mathbf{Q}_{t+1} = \mathsf{QR}[(\mathbf{I} + \eta_{t+1} x_{t+1} x_{t+1}^\top)\mathbf{Q}_t] = (\mathbf{I} + \eta_{t+1} x_{t+1} x_{t+1}^\top)\mathbf{Q}_t \mathbf{S}_t$$

for some $\mathbf{S}_t \in \mathbb{R}^{k \times k}$ by the definition of Gram-Schmidt. This implies that

$$\mathbf{Q}_{t+1} = (\mathbf{I} + \eta_{t+1} x_{t+1} x_{t+1}^\top)\widetilde{\mathbf{Q}}_t \mathbf{R}_t \mathbf{S}_t = \mathbf{P}_{t+1}\mathbf{Q}\mathbf{R}_t \mathbf{S}_t = \widetilde{\mathbf{Q}}_{t+1}\mathbf{R}_t \mathbf{S}_t = \widetilde{\mathbf{Q}}_{t+1}\mathbf{R}_{t+1}$$

---
[10] The second equality is a simple fact but anyways proved in Lemma 2.2 later.

if we define $\mathbf{R}_{t+1} = \mathbf{R}_t \mathbf{S}_t$. This completes the proof that $\mathbf{Q}_t = \widetilde{\mathbf{Q}}_t \mathbf{R}_t$. As a result, since each $\mathbf{Q}_t$ is column orthogonal for $t \geq 1$, we have $\|\mathbf{Q}_t^\top \mathbf{V}\|_2 \leq 1$ and therefore $\|\mathbf{X}^\top \mathbf{Q}_t\|_F \leq \|\mathbf{X}^\top \mathbf{Q}_t (\mathbf{V}^\top \mathbf{Q}_t)^{-1}\|_F \|\mathbf{V}^\top \mathbf{Q}_t\|_2 \leq \|\mathbf{X}^\top \mathbf{Q}_t (\mathbf{V}^\top \mathbf{Q}_t)^{-1}\|_F$. Finally,
$$\|\mathbf{X}^\top \mathbf{Q}_t\|_F \leq \|\mathbf{X}^\top \mathbf{Q}_t (\mathbf{V}^\top \mathbf{Q}_t)^{-1}\|_F = \|\mathbf{X}^\top \widetilde{\mathbf{Q}}_t \mathbf{R}_t (\mathbf{V}^\top \widetilde{\mathbf{Q}}_t \mathbf{R}_t)^{-1}\|_F \leq \|\mathbf{X}^\top \widetilde{\mathbf{Q}}_t (\mathbf{V}^\top \widetilde{\mathbf{Q}}_t)^{-1}\|_F \ . \quad \square$$

**Observation.** Due to Lemma 2.2, in order to prove our upper bound theorems, it suffices to upper bound the quantity $\|\mathbf{X}^\top \mathbf{P}_t \mathbf{Q} (\mathbf{V}^\top \mathbf{P}_t \mathbf{Q})^{-1}\|_F$ for $\mathbf{X} = \mathbf{W}$ (gap-free) or $\mathbf{X} = \mathbf{Z}$ (gap-dependent).

## 3 Overview of Our Proofs and Techniques

**Oja's Algorithm.** To illustrate the idea, let us simply focus on the gap-dependent case. Denoting in this section by $s_t \stackrel{\text{def}}{=} \|\mathbf{Z}^\top \mathbf{P}_t \mathbf{Q} (\mathbf{V}^\top \mathbf{P}_t \mathbf{Q})^{-1}\|_F$, owing to Lemma 2.2, we want to bound $s_t$ in terms of $x_t$ and $s_{t-1}$. A simple calculation using the Sherman-Morrison formula gives

$$\mathbb{E}[s_t^2] \leq (1 - \eta_t \mathsf{gap}) \mathbb{E}[s_{t-1}^2] + \mathbb{E}\left[\left(\frac{\eta_t a_t}{1 - \eta_t a_t}\right)^2\right] \quad \text{where} \quad a_t = \|x_t^\top \mathbf{P}_{t-1} \mathbf{Q} (\mathbf{V}^\top \mathbf{P}_{t-1} \mathbf{Q})^{-1}\|_2 \quad (3.1)$$

At a first look, $\mathbb{E}[s_t^2]$ is decaying by a multiplicative factor $(1 - \eta_t \mathsf{gap})$ at every iteration; however, this bound could be *problematic* when $\eta_t a_t$ is close to 1 and thus we need to ensure $\eta_t \leq \frac{1}{a_t}$ with high probability for every step.

One can naively bound $a_t \leq \|\mathbf{P}_{t-1}\mathbf{Q}(\mathbf{V}^\top \mathbf{P}_{t-1}\mathbf{Q})^{-1}\|_2 \leq s_{t-1} + 1$. However, since $s_{t-1}$ can be $\Omega(\sqrt{d})$ even at $t = 1$, we must choose $\eta_t \leq O(1/\sqrt{d})$ and the resulting convergence rate will be *not* efficient (i.e., proportional to $d$). This is why most known global convergence results are not efficient (see Table 1). On the other hand, if one ignores initialization and starts from a point $t_0$ when $s_{t_0} \leq 1$ is already satisfied, then we can prove a *local* convergence rate that is efficient (c.f. [18]). Note that this local rate is still slower than ours by a factor $\lambda_1 + \cdots + \lambda_k$.

Our *first contribution* is the following crucial observation: for a random initial matrix $\mathbf{Q}$, $a_1 = \|x_1^\top \mathbf{Q}(\mathbf{V}^\top \mathbf{Q})^{-1}\|_2$ is actually quite small. A simple fact on the singular value distribution of inverse-Wishart distribution implies $a_1 = O(\sqrt{k})$ with high probability. Thus, at least in the first iteration, we can set $\eta_1 = \Omega(1/\sqrt{k})$ independent of the dimension $d$. Unfortunately, in subsequent iterations, it is not clear whether $a_t$ remains small or increases.

Our *second contribution* is to control $a_t$ using the fact that $a_t$ itself "forms another random process." More precisely, denoting by $a_{t,s} = \|x_t^\top \mathbf{P}_s \mathbf{Q} (\mathbf{V}^\top \mathbf{P}_s \mathbf{Q})^{-1}\|_2$ for $0 \leq s \leq t-1$, we wish to bound $a_{t,s}$ in terms of $a_{t,s-1}$ and show that it does not increase by much. (If we could achieve so, combining it with the initialization $a_{t,0} \leq O(\sqrt{k})$, we would know that all $a_{t,s}$ are small for $s \leq t - 1$.) Unfortunately, since $x_t$ is not an eigenvector of $\mathbf{\Sigma}$, the recursion looks like

$$\mathbb{E}[a_{t,s}^2] \leq (1 - \eta_s \lambda_k) \mathbb{E}[a_{t,s-1}^2] + \eta_s \lambda_k \mathbb{E}[b_{t,s-1}^2] + \mathbb{E}\left[\left(\frac{\eta_s a_{s,s-1}}{1 - \eta_s a_{s,s-1}}\right)^2\right] \quad (3.2)$$

where $b_{t,s} = \|x_t^\top \mathbf{\Sigma} \mathbf{P}_s \mathbf{Q} (\mathbf{V}^\top \mathbf{P}_s \mathbf{Q})^{-1}\|_2$. Now three difficulties arise from formula (3.2):

- $b_{t,s}$ can be very different from $a_{t,s}$ — in worse case, the ratio between them can be unbounded.
- the problematic term now becomes $a_{s,s-1}$ (rather than the original $a_t = a_{t,t-1}$ in (3.1)) which is not present in the chain $\{a_{t,s}\}_{s=1}^{t-1}$.
- since $b_{t,s}$ differs from $a_{t,s}$ by an additional factor $\mathbf{\Sigma}$ in the middle, to analyze $b_{t,s}$, we need to further study $\|x_t^\top \mathbf{\Sigma}^2 \mathbf{P}_s \mathbf{Q} (\mathbf{V}^\top \mathbf{P}_s \mathbf{Q})^{-1}\|_2$ and so on.

We solve these issues by carefully considering a multi-dimensional random process $c_{t,s}$ with $c_{t,s}^{(i)} \stackrel{\text{def}}{=} \|x_t^\top \mathbf{\Sigma}^i \mathbf{P}_s \mathbf{Q} (\mathbf{V}^\top \mathbf{P}_s \mathbf{Q})^{-1}\|_2$. Ignoring the last term, we can derive that

$$\forall t, \forall s \leq t-1, \quad \mathbb{E}\left[\left(c_{t,s}^{(i)}\right)^2\right] \lesssim (1 - \eta_s \lambda_k) \mathbb{E}\left[\left(c_{t,s-1}^{(i)}\right)^2\right] + \eta_s \lambda_k \mathbb{E}\left[\left(c_{t,s-1}^{(i+1)}\right)^2\right] \ . \quad (3.3)$$

Our *third contribution* is a new random process concentration bound to control the change in this multi-dimensional chain (3.3). To achieve this, we adapt the prove of standard Chernoff bound to multi dimensions (which is not the same as matrix concentration bound). After establishing this non-trivial concentration result, all terms of $a_t = c_{t,t-1}^{(0)}$ can be simultaneously bounded by a constant. This ensures that the problematic term in (3.1) is well-controlled.

The overall plan looks promising; however, there are holes in the above thought experiment.

- In order to apply a random-process concentration bound (e.g., Azuma concentration), we need the process to *not depend on the future*. However, the random vector $c_{t,s}$ is not $\mathcal{F}_{\leq s}$ measurable but $\mathcal{F}_{\leq s} \vee \mathcal{F}_t$ measurable (i.e., it depends on $x_t$ for a future $t > s$).

- Furthermore, the expectation bounds such as (3.1), (3.2), (3.3) only hold if $\mathbb{E}[x_t x_t] = \mathbf{\Sigma}$; however, if we take away a failure event $\mathcal{C}$ —$\mathcal{C}$ may correspond to the event when $a_t$ is large— the conditional expectation $\mathbb{E}[x_t x_t \mid \overline{\mathcal{C}}]$ becomes $\mathbf{\Sigma} + \mathbf{\Delta}$ where $\mathbf{\Delta}$ is some error matrix. This can amplify the failure probability in next iteration.

Our *fourth contribution* is a "decoupling" framework to deal with the above issues (Appendix i.D). At a high level, to deal with the first issue we fix $x_t$ and study $\{c_{t,s}\}_{s=0,1,\ldots,t-1}$ conditioning on $x_t$; in this way the process decouples and each $c_{t,s}$ becomes $\mathcal{F}_{\leq s}$ measurable. We can do so because we can carefully ensure that the failure events only depend on $x_s$ for $s \leq t-1$ but not on $x_t$. To deal with the second issue, we convert the random process into an unconditional random process (see (i.D.2)); this is a generalization of using stopping time on martingales. Using these tools, we manage to show that the failure probability only grows linearly with respect to $T$ and henceforth bound the value of $c_{t,s}^{(i)}$ for all $t, s$ and $i$.

Putting them together, we are able to show that Oja's algorithm achieves convergence rate $\widetilde{O}\left(\frac{\lambda_1 + \ldots + \lambda_k}{\mathsf{gap}^2}(\frac{1}{\varepsilon} + k)\right)$. The rate matches our lower bound asymptotically when $\lambda_1$ and $\lambda_k$ are within a constant factor of each other, however, if we only care about crude approximation of the eigenvectors (e.g. for constant $\varepsilon$), then the Oja's algorithm is off by a factor $k$.

*Remark* 3.1. The ideas above are insufficient for our gap-free results. To prove Theorem 2 and 3, we also need to bound quantities $s'_t \stackrel{\text{def}}{=} \|\mathbf{W}^\top \mathbf{P}_t \mathbf{Q} (\mathbf{V}^\top \mathbf{P}_t \mathbf{Q})^{-1}\|_F$ where $\mathbf{W}$ consists of all eigenvectors of $\mathbf{\Sigma}$ with values no more than $\lambda_k - \rho$. This is so because the interesting quantity in a gap-free case changes from $s_t$ to $s'_t$ according to Lemma 2.2. Now, to bound $s'_t$ one has to bound $c_{t,s}$; however, the $c_{t,s}$ process still depends on the original $s_t$ as opposed to $s'_t$. In sum, we unavoidably have to bound $s_t, s'_t$, and $c_{t,s}$ all together, making the proofs even more sophisticated.

**Our Oja$^{++}$ Algorithm.** The factor $k$ in Oja's algorithm comes from the fact that the earlier quantity $a_1 = \|x_1^\top \mathbf{Q}(\mathbf{V}^\top \mathbf{Q})^{-1}\|_2$ is at least $\Omega(\sqrt{k})$ at $t=1$, so we must set $\eta_1 \leq O(1/\sqrt{k})$ and this incurs a factor $k$ in the running time. After warm start, $a_t$ drops to $O(1)$ and we can choose $\eta_t \leq 1$.

A similar issue was also observed by Hardt and Price [9] and they solved it using "oversampling". Namely, to put it into our setting, we can use a $d \times 2k$ random starting matrix $\mathbf{Q}_0$ and run Oja's to produce $\mathbf{Q}_T \in \mathbb{R}^{d \times 2k}$. In this way, the quantity $a_1$ becomes $O(1)$ even at the beginning due to some property of the inverse-Wishart distribution. However, the output $\mathbf{Q}_T$ is now a $2k$ dimensional space that is only guaranteed to "approximately contain" the top $k$ eigenvectors. It is not clear how to find this $k$-subspace (recall the algorithm does not see $\mathbf{\Sigma}$).

> Our *key observation* behind Oja$^{++}$ is that a similar effect also occurs via "under-sampling".

If we initialize $\mathbf{Q}_0$ randomly with dimension $d \times k/2$, we can also obtain a speed-up factor of $k$. Unlike Hardt and Price, this time the output $\mathbf{Q}_{T_0}$ is a $k/2$-dimensional subspace that approximately lies *entirely* in the column span of $\mathbf{V} \in \mathbb{R}^{d \times k}$.

After getting $\mathbf{Q}_{T_0}$, one could hope to get the rest by running the same algorithm again, but restricted to the orthogonal complement of $\mathbf{Q}_{T_0}$. This approach would work if $\mathbf{Q}_{T_0}$ were exactly the eigenvectors of $\mathbf{\Sigma}$; however, due to approximation error, this approach would eventually lose a factor $1/\texttt{gap}$ in the sample complexity which is even bigger than the factor $k$ that we could gain.

Instead, our $\mathsf{Oja}^{++}$ algorithm is divided into $\log k$ epochs. At each epoch $i = 1, 2, \ldots, \log k$, we attach $k/2^i$ new random columns to the working matrix $\mathbf{Q}_t$ in Oja's algorithm, and then run Oja's for a fixed number of iterations. Note that every time (except the last time) we attach new random columns, we are in an "under-sampling" mode because if we add $k/2^i$ columns there must be $k/2^i$ remaining dimensions. This ensures that the quantity $a_t$ only increases by a constant so we have $a_t = O(1)$ throughout execution of $\mathsf{Oja}^{++}$. Finally, there are only $\log k$ epochs so the total running time is still $\widetilde{O}\left(\frac{\lambda_1 + \ldots + \lambda_k}{\rho^2} \frac{1}{\varepsilon}\right)$ and this $\widetilde{O}$ notion hides a $\log k$ factor.

**Roadmap.** Our proofs are highly technical so we carefully choose what to present in this main body. In Section 5 we state properties of the initial matrix $\mathbf{Q}$ which corresponds to our first contribution. In Section 6 we provide expected guarantees on $s_t$, $s'_t$ and $a_{t,s}$ and they correspond to our second contribution. The third (martingale lemmas) and fourth contributions (decoupling lemma) are deferred to the appendix.

Most importantly, in Section 7 we present (although in weaker forms) two Main Lemmas to deal with the convergence one for $t \leq T_0$ (before warm start) and one for $t > T_0$ (after warm start). These sections, when put together, directly imply two weaker variants of Theorem 1 and 2. We state these weaker variants in Appendix i and include all the mathematical details there.

Appendix ii includes our Rayleigh quotient Theorem 3 and lower bound Theorem 6. Appendix iii strengthens the main lemmas into their stronger forms, and prove Theorem 1 and 2 formally. Our $\mathsf{Oja}^{++}$ results, namely Theorem 4 and 5, are also proved in Appendix iii.

In Figure 1 on page 14, we present a dependency graph of all of our main theorems and lemmas. We hope that the readers could appreciate our organization of this paper.

## 4 Discussions, Extensions and Future Directions

In this paper we give global convergence analysis of the Oja's algorithm, and a twisted version $\mathsf{Oja}^{++}$ which has better complexity. We also give an information-theoretic lower bound showing that *any* algorithm, *offline or online*, must have final accuracy $\mathbb{E}_{x_1, \ldots, x_T, \mathcal{A}}\left[\|\mathbf{W}^\top \mathbf{Q}_T\|_F^2\right] = \Omega\big(\frac{k\lambda_k}{\texttt{gap}^2 \cdot T}\big)$. This matches our gap-dependent result on $\mathsf{Oja}^{++}$ when $\lambda_1 + \cdots + \lambda_k = O(k\lambda_k)$; that is, when there is an eigengap and when *"the spectrum is flat."*

When the spectrum is not flat, our algorithm can be improved to have better accuracy. However, this requires good prior knowledge of $\lambda_1, \cdots \lambda_k$, and may not be realistic.

In the gap-free case, $\mathsf{Oja}^{++}$ only achieves accuracy $O\big(\frac{\lambda_1 + \cdots \lambda_{k+m}}{\rho^2 \cdot T}\big)$, which appears worse than the lower bound $O\big(\frac{k\lambda_k}{\rho^2 \cdot T}\big)$. In fact, we can also achieve $O\big(\frac{\lambda_1 + \cdots \lambda_k}{\rho^2 \cdot T}\big)$ if we allow more space, namely, space up to $O((k+m)d)$ as opposed to $O(kd)$. More generally, we have a space-accuracy tradeoff.

If we run $\mathsf{Oja}^{++}$ on $k'$ initial random vectors, and thus using space $O(dk')$ for $k' \in [k, k+m]$, we can randomly pick $k$ columns from the output and have the following accuracy:

**Theorem 4.1** (tradeoff)**.** *For every $k' \in [k, k+m]$ with $\lambda_{k'} - \lambda_{k+m} \geq \frac{\rho}{\log d}$, let $\mathbf{Q} \in \mathbb{R}^{d \times k'}$ be a random gaussian matrix and $\mathbf{Q}_T \in \mathbb{R}^{d \times k'}$ be the output of $\mathsf{Oja}^{++}$ with random input $\mathbf{Q}$. Then, letting $\mathbf{Q}'_T \in \mathbb{R}^{d \times k}$ be $k$ random columns of $\mathbf{Q}_T$ (chosen uniformly at random), we have*

$$\mathbb{E}\left[\|\mathbf{W}^\top \mathbf{Q}'_T\|_F^2\right] = \widetilde{O}\left(\frac{k}{k'} \frac{\lambda_1 + \cdots \lambda_{k+m'}}{\rho^2 T}\right)$$

where $m' \leq m$ is any index satisfying $\lambda_{k'} - \lambda_{k+m'} \geq \frac{\rho}{\log d}$.

*Proof of Theorem 4.1.* Observe that $\mathsf{Oja}^{++}$ guarantees (using $\rho/\log d$ instead of $\rho$ as the gap-free parameter) $\mathbb{E}[\|\mathbf{W}^\top \mathbf{Q}_T\|_F^2] = \widetilde{O}\left(\frac{\lambda_1 + \cdots \lambda_{k+m'}}{\rho^2 T}\right)$. Then, $k$ random columns of $\mathbf{Q}'_T$ decreases the squared Frobenius norm by a factor of $k'/k$. □

We also have the following crucial observation:

**Corollary 4.2.** *There exists $k' \in [k, k+m]$ such that $\lambda_{k'} - \lambda_{k+m} \geq \frac{\rho}{\log d}$ and*
$$\frac{k}{k'}(\lambda_1 + \cdots \lambda_{k+m'}) = O(\lambda_1 + \cdots \lambda_k) \ .$$

*Proof of Corollary 4.2.* The proof is by counting. Divide $[\lambda_k - \rho, \lambda_k]$ into $\log d$ intervals of equal span in descending order $[\lambda_k - \frac{\rho}{\log d}, \lambda_k), [\lambda_k - 2\frac{\rho}{\log d}, \lambda_k - \frac{\rho}{\log d}), \cdots, [\lambda_k - \rho, \lambda_k - (1 - \frac{1}{\log d})\rho)$, and let $S_i \subseteq \{k+1, \ldots, k+m\}$ be the indices of $\lambda_j$ is in the $i$-th interval above, for $i = 1, 2, \ldots, \log d$.

Define $\Lambda_i = \sum_{j \in S_i} \lambda_j$ and $\Lambda_0 = \Lambda = \lambda_1 + \cdots \lambda_k$. Also define $k_i = k + \sum_{1 \leq j < i} |S_j|$. We then have, for every $k' = k_i$, we can choose $m'$ such that $\lambda_1 + \cdots \lambda_{k+m'} = \Lambda + \sum_{j=1}^{i} \Lambda_j$. By $\Lambda_0 \leq d\Lambda_{\log d}$, we know that there must exist some $i$ such that $\Lambda_i \leq 100\Lambda_{i-1}$. We compute
$$\frac{k}{k_i}\left(\Lambda + \sum_{j=1}^{i} \Lambda_j\right) \leq 100\frac{k}{k_i}\left(\Lambda + \sum_{j=1}^{i-1} \Lambda_i\right) = 100\frac{k}{k_i}(\lambda_1 + \cdots + \lambda_{k_i}) \leq 100\Lambda$$
so we have found such $k'$ satisfying the statement. □

**In sum.** In the gap-free case, if we increase the space of $\mathsf{Oja}^{++}$ to at most $O((k+m)d)$, we can achieve accuracy $O\left(\frac{\lambda_1 + \cdots \lambda_k}{\rho^2 \cdot T}\right)$, and thus match the lower bound when the spectrum is flat.[11]

It was also studied by Balcan et al. [3] that increasing the space could enhance the performance. However, their algorithm always uses space $\Omega((k+m)d)$, and furthermore, in the Frobenius-norm accuracy, their performance is always worse than $\mathsf{Oja}^{++}$, not to say $\mathsf{Oja}^{++}$ only uses space $O(kd)$.[12] It is an important future direction to directly get $O\left(\frac{\lambda_1 + \cdots \lambda_k}{\rho^2 \cdot T}\right)$ without increasing the space.

**Matrix Bernstein.** Using matrix Bernstein and a gap-free variant of the Wedin theorem [2], one can show that, if we simply let $\mathbf{Q}_T$ be the top-$k$ eigenvectors of the empirical covariance matrix $\sum_{t=1}^{T} x_t x_t^\top$, then we have $\mathbb{E}[\|\mathbf{W}^\top \mathbf{Q}_T\|_2^2] = \widetilde{O}\left(\frac{\lambda_1}{\rho^2 T}\right)$. If one directly translates this to a Frobenius norm bound, it gives $\mathbb{E}[\|\mathbf{W}^\top \mathbf{Q}_T\|_F^2] = \widetilde{O}\left(\frac{\lambda_1 k}{\rho^2 T}\right)$ and is worse than ours. However, our result, if naively translated to spectral norm, also loses a factor $k$. It is a future direction to directly get a spectral-norm guarantee for streaming PCA.

## 5 Random Initialization

We state our main lemma for initialization. Let $\mathbf{Q} \in \mathbb{R}^{d \times k}$ be a matrix with each entry i.i.d drawn from $\mathcal{N}(0, 1)$, the standard gaussian.

**Lemma 5.1** (initialization). *For every $p, q \in (0, 1)$, every $T \in \mathbb{N}^*$, every distribution on vector set $\{x_t\}_{t=1}^T$ with $\|x_t\|_2 \leq 1$, with probability at least $1 - p - 2q$ over the random choice of $\mathbf{Q}$:*
$$\begin{cases} \left\|(\mathbf{Z}^\top \mathbf{Q})(\mathbf{V}^\top \mathbf{Q})^{-1}\right\|_F^2 \leq \frac{576dk}{p^2} \ln \frac{d}{p} & \text{and} \\ \mathbf{Pr}_{x_1, \ldots, x_T}\left[\exists i \in [T], \exists t \in [T], \left\|x_t^\top \mathbf{Z}\mathbf{Z}^\top (\mathbf{\Sigma}/\lambda_{k+1})^{i-1} \mathbf{Q}(\mathbf{V}^\top \mathbf{Q})^{-1}\right\|_2 \geq \frac{18}{p}\left(2k \ln \frac{T}{q}\right)^{1/2}\right] \leq q \end{cases}$$

---
[11] Of course, this requires the algorithm to know $k'$ which can be done by trying $k' = k+2, k+4, k+8$, etc.
[12] For instance, when $\lambda_{1 \sim k+m} \geq \frac{k+m}{d}$, their global convergence is $\widetilde{O}\left(\frac{dk(\lambda_{1 \sim k+m})^2 \lambda_k}{(k+m)\rho^3 T}\right)$, but ours is only $\widetilde{O}\left(\frac{\lambda_{1 \sim k+m}}{\rho^2 T}\right)$.

The two statements of the above lemma correspond to $s_0$ and $c_{t,0}^{(i)}$ that we defined in Section 3. The second statement is of the form "$\mathbf{Pr}[event] \leq q$" instead of "for every fixed $x_t$, $event$ holds with probability $\leq q$" because we cannot afford taking union bound on $x_t$.

## 6 Expected Results

We now formalize inequalities (3.1), (3.2), (3.3) which characterize to the behavior of our target random processes. Let $\mathbf{X} \in \mathbb{R}^{d \times r}$ be a generic matrix that shall later be chosen as either $\mathbf{X} = \mathbf{W}$ (corresponding to $s_t'$), $\mathbf{X} = \mathbf{Z}$ (corresponding to $s_t$), or $\mathbf{X} = [w]$ for some vector $w$ (corresponding to $c_{t,s}^{(i)}$). We introduce the following notions that shall be used extensively:

$$\mathbf{L}_t = \mathbf{P}_t \mathbf{Q} (\mathbf{V}^\top \mathbf{P}_t \mathbf{Q})^{-1} \in \mathbb{R}^{d \times k} \qquad \mathbf{R}_t' = \mathbf{X}^\top x_t x_t^\top \mathbf{L}_{t-1} \in \mathbb{R}^{r \times k}$$
$$\mathbf{S}_t = \mathbf{X}^\top \mathbf{L}_t \in \mathbb{R}^{r \times k} \qquad \mathbf{H}_t' = \mathbf{V}^\top x_t x_t^\top \mathbf{L}_{t-1} \in \mathbb{R}^{k \times k}$$

**Lemma 6.1** (Appendix i.B). *For every $t \in [T]$, suppose $\mathcal{C}_{\leq t}$ is an event on random $x_1, \ldots, x_t$ and*

$$\mathcal{C}_{\leq t} \text{ implies } \|x_t^\top \mathbf{L}_{t-1}\|_2 = \|x_t^\top \mathbf{P}_{t-1} \mathbf{Q} (\mathbf{V}^\top \mathbf{P}_{t-1} \mathbf{Q})^{-1}\|_2 \leq \phi_t \quad \text{where} \quad \eta_t \phi_t \leq \frac{1}{2} \;,$$

*and suppose $\mathbb{E}_{x_t}\left[x_t x_t^\top \mid \mathcal{F}_{\leq t-1}, \mathcal{C}_{\leq t}\right] = \mathbf{\Sigma} + \mathbf{\Delta}$. Then, we have:*
*(a) When $\mathbf{X} = \mathbf{Z}$,*

$$\mathbb{E}\left[\mathbf{Tr}(\mathbf{S}_t^\top \mathbf{S}_t) \mid \mathcal{F}_{\leq t-1}, \mathcal{C}_{\leq t}\right] \leq (1 - 2\eta_t \mathsf{gap} + 14\eta_t^2 \phi_t^2) \mathbf{Tr}(\mathbf{S}_{t-1} \mathbf{S}_{t-1}^\top) + 10\eta_t^2 \phi_t^2$$
$$+ 2\eta_t \|\mathbf{\Delta}\|_2 \left( \left[\mathbf{Tr}(\mathbf{S}_{t-1}^\top \mathbf{S}_{t-1})\right]^{3/2} + 2\mathbf{Tr}(\mathbf{S}_{t-1}^\top \mathbf{S}_{t-1}) + \left[\mathbf{Tr}(\mathbf{S}_{t-1}^\top \mathbf{S}_{t-1})\right]^{1/2} \right)$$

*(b) When $\mathbf{X} = \mathbf{W}$,*

$$\mathbb{E}\left[\mathbf{Tr}(\mathbf{S}_t^\top \mathbf{S}_t) \mid \mathcal{F}_{\leq t-1}, \mathcal{C}_{\leq t}\right] \leq (1 - 2\eta_t \rho + 14\eta_t^2 \phi_t^2) \mathbf{Tr}(\mathbf{S}_{t-1} \mathbf{S}_{t-1}^\top) + 10\eta_t^2 \phi_t^2$$
$$+ 2\eta_t \|\mathbf{\Delta}\|_2 \left( \left[\mathbf{Tr}(\mathbf{S}_{t-1}^\top \mathbf{S}_{t-1})\right]^{1/2} + \mathbf{Tr}(\mathbf{S}_{t-1}^\top \mathbf{S}_{t-1}) \right) \left( 1 + \left[\mathbf{Tr}(\mathbf{Z}^\top \mathbf{L}_{t-1} \mathbf{L}_{t-1}^\top \mathbf{Z})\right]^{1/2} \right)$$

*(c) When $\mathbf{X} = [w] \in \mathbb{R}^{d \times 1}$ where $w$ is a vector with Euclidean norm at most 1,*

$$\mathbb{E}\left[\mathbf{Tr}(\mathbf{S}_t^\top \mathbf{S}_t) \mid \mathcal{F}_{\leq t-1}, \mathcal{C}_{\leq t}\right] \leq \left(1 - \eta_t \lambda_k + 14\eta_t^2 \phi_t^2\right) \mathbf{Tr}(\mathbf{S}_{t-1} \mathbf{S}_{t-1}^\top) + 10\eta_t^2 \phi_t^2 + \frac{\eta_t}{\lambda_k} \|w^\top \mathbf{\Sigma} \mathbf{L}_{t-1}\|_2^2$$
$$+ 2\eta_t \|\mathbf{\Delta}\|_2 \left( \left[\mathbf{Tr}(\mathbf{S}_{t-1}^\top \mathbf{S}_{t-1})\right]^{1/2} + \mathbf{Tr}(\mathbf{S}_{t-1}^\top \mathbf{S}_{t-1}) \right) \left( 1 + \left[\mathbf{Tr}(\mathbf{Z}^\top \mathbf{L}_{t-1} \mathbf{L}_{t-1}^\top \mathbf{Z})\right]^{1/2} \right)$$

The above three expectation results will be used to provide upper bounds on the quantities we care about (i.e., $s_t$, $s_t'$, $c_{t,s}^{(i)}$). In the appendix, to enable proper use of martingale concentration, we also bound their absolute changes $|\mathbf{Tr}(\mathbf{S}_t^\top \mathbf{S}_t) - \mathbf{Tr}(\mathbf{S}_{t-1} \mathbf{S}_{t-1}^\top)|$ and variance $\mathbb{E}\left[|\mathbf{Tr}(\mathbf{S}_t^\top \mathbf{S}_t) - \mathbf{Tr}(\mathbf{S}_{t-1} \mathbf{S}_{t-1}^\top)|^2\right]$ in changes.[13]

## 7 Main Lemmas

Our main lemmas in this section can be proved by combining (1) the expectation results in Section 6, (2) the martingale concentrations in Appendix i.C, and (3) our decoupling lemma in Appendix i.D.

---

[13]Recall that even in the simplest martingale concentration, one needs upper bounds on the absolute difference between consecutive variables; furthermore, the concentration can be tightened if one also has an (expected) variance upper bound between variables.

Our first lemma describes the behavior of quantities $s_t = \|\mathbf{Z}^\top \mathbf{P}_t \mathbf{Q}(\mathbf{V}^\top \mathbf{P}_t \mathbf{Q})^{-1}\|_F$ and $s'_t = \|\mathbf{W}^\top \mathbf{P}_t \mathbf{Q}(\mathbf{V}^\top \mathbf{P}_t \mathbf{Q})^{-1}\|_F$ before warm start. At a high level, it shows if the $s_t$ sequence starts from $s_0^2 \leq \Xi_\mathbf{Z}$, under mild conditions, $s_t^2$ never increases to more than $2\Xi_\mathbf{Z}$. Note that $\Xi_\mathbf{Z} = O(d)$ according to Lemma 5.1. The other sequence $(s'_t)^2$ also never increases to more than $2\Xi_\mathbf{Z}$ because $s'_t \leq s_t$, but most importantly, $(s'_t)^2$ drops below 2 after $t \geq T_0$. Therefore, at point $t = T_0$ we need to adjust the learning rate so the algorithm achieves best convergence rate, and this is the goal of our Lemma Main 2. (We emphasize that although we are only interested in $s_t$ and $s'_t$, our proof of the lemma also needs to bound the multi-dimensional $c_{t,s}$ sequence discussed in Section 3.)

**Lemma Main 1** (before warm start). *For every $\rho \in (0,1)$, $q \in \left(0, \frac{1}{2}\right]$, $\Xi_\mathbf{Z} \geq 2$, $\Xi_x \geq 2$, and fixed matrix $\mathbf{Q} \in \mathbb{R}^{d \times k}$, suppose it satisfies*

- $\|\mathbf{Z}^\top \mathbf{Q}(\mathbf{V}^\top \mathbf{Q})^{-1}\|_F^2 \leq \Xi_\mathbf{Z}$, *and*
- $\mathbf{Pr}_{x_t}\left[\forall j \in [T], \left\|x_t^\top \mathbf{Z}\mathbf{Z}^\top \left(\mathbf{\Sigma}/\lambda_{k+1}\right)^{j-1} \mathbf{Q}(\mathbf{V}^\top \mathbf{Q})^{-1}\right\|_2 \leq \Xi_x\right] \geq 1 - q^2/2$ *for every $t \in [T]$.*

*Suppose also the learning rates $\{\eta_s\}_{s \in [T]}$ satisfy*

$$(1): \quad \forall s \in [T], q\Xi_\mathbf{Z}^{3/2} \leq \eta_s \leq \frac{\rho}{4000\Xi_x^2 \ln \frac{24T}{q^2}} \qquad (2): \quad \sum_{t=1}^T \eta_t^2 \Xi_x^2 \leq \frac{1}{100 \ln \frac{32T}{q^2}}$$

$$(3): \quad \exists T_0 \in [T] \text{ such that } \sum_{t=1}^{T_0} \eta_t \geq \frac{\ln(3\Xi_\mathbf{Z})}{\rho} \ .$$

*Then, for every $t \in [T-1]$, with probability at least $1 - 2qT$ (over the randomness of $x_1, \ldots, x_t$):*

- $\|\mathbf{Z}^\top \mathbf{P}_t \mathbf{Q}(\mathbf{V}^\top \mathbf{P}_t \mathbf{Q})^{-1}\|_F^2 \leq 2\Xi_\mathbf{Z}$, *and*
- *if $t \geq T_0$ then $\left\|\mathbf{W}^\top \mathbf{P}_t \mathbf{Q}(\mathbf{V}^\top \mathbf{P}_t \mathbf{Q})^{-1}\right\|_F^2 \leq 2$.*

Our second lemma asks for a stronger assumption on the learning rates and shows that after warm start (i.e., for $t \geq T_0$), the quantity $(s'_t)^2$ scales as $1/t$.

**Lemma Main 2** (after warm start). *In the same setting as Lemma Main 1, if there exists $\delta \leq 1/\sqrt{8}$ s.t.*

$$\frac{T_0}{\ln^2 T_0} \geq \frac{9\ln(8/q^2)}{\delta^2} \ , \quad \forall s \in \{T_0+1, \ldots, T\}: \quad 2\eta_s\rho - 56\eta_s^2 \Xi_x^2 \geq \frac{1}{s-1} \quad \text{and} \quad \eta_s \leq \frac{1}{20(s-1)\delta\Xi_x} \ ,$$

*then, with probability at least $1 - 2qT$ (over the randomness of $x_1, \ldots, x_T$):*

- $\|\mathbf{Z}^\top \mathbf{P}_t \mathbf{Q}(\mathbf{V}^\top \mathbf{P}_t \mathbf{Q})^{-1}\|_F^2 \leq 2\Xi_\mathbf{Z}$ *for every $t \in \{T_0, \ldots, T\}$, and*
- $\|\mathbf{W}^\top \mathbf{P}_t \mathbf{Q}(\mathbf{V}^\top \mathbf{P}_t \mathbf{Q})^{-1}\|_F^2 \leq \frac{5T_0/\ln^2(T_0)}{t/\ln^2 t}$ *for every $t \in \{T_0, \ldots, T\}$.*

**Parameter 7.1.** There exist constants $C_1, C_2, C_3 > 0$ such that for every $q > 0$ that is sufficiently small (meaning $q < 1/\mathsf{poly}(T, \Xi_\mathbf{Z}, \Xi_x, 1/\rho)$), the following parameters satisfy both Lemma Main 1 and Lemma Main 2:

$$\frac{T_0}{\ln^2(T_0)} = C_1 \cdot \frac{\Xi_x^2 \ln \frac{T}{q} \ln^2 \Xi_\mathbf{Z}}{\rho^2} \ , \quad \eta_t = C_2 \cdot \begin{cases} \frac{\ln \Xi_\mathbf{Z}}{T_0 \cdot \rho} & t \leq T_0; \\ \frac{1}{t \cdot \rho} & t > T_0. \end{cases} \ , \quad \text{and} \quad \delta = C_3 \cdot \frac{\rho}{\Xi_x} \ .$$

Using such learning rates for our main lemmas, one can prove in one page (see Appendix i.F)

- a weaker version of Theorem 2 where $(\Lambda_1, \Lambda_2)$ are replaced by $(1, 0)$, and
- a weaker version of Theorem 1 where $\Lambda = \lambda_1 + \cdots + \lambda_k$ is replaced by 1.

# Appendix Overview

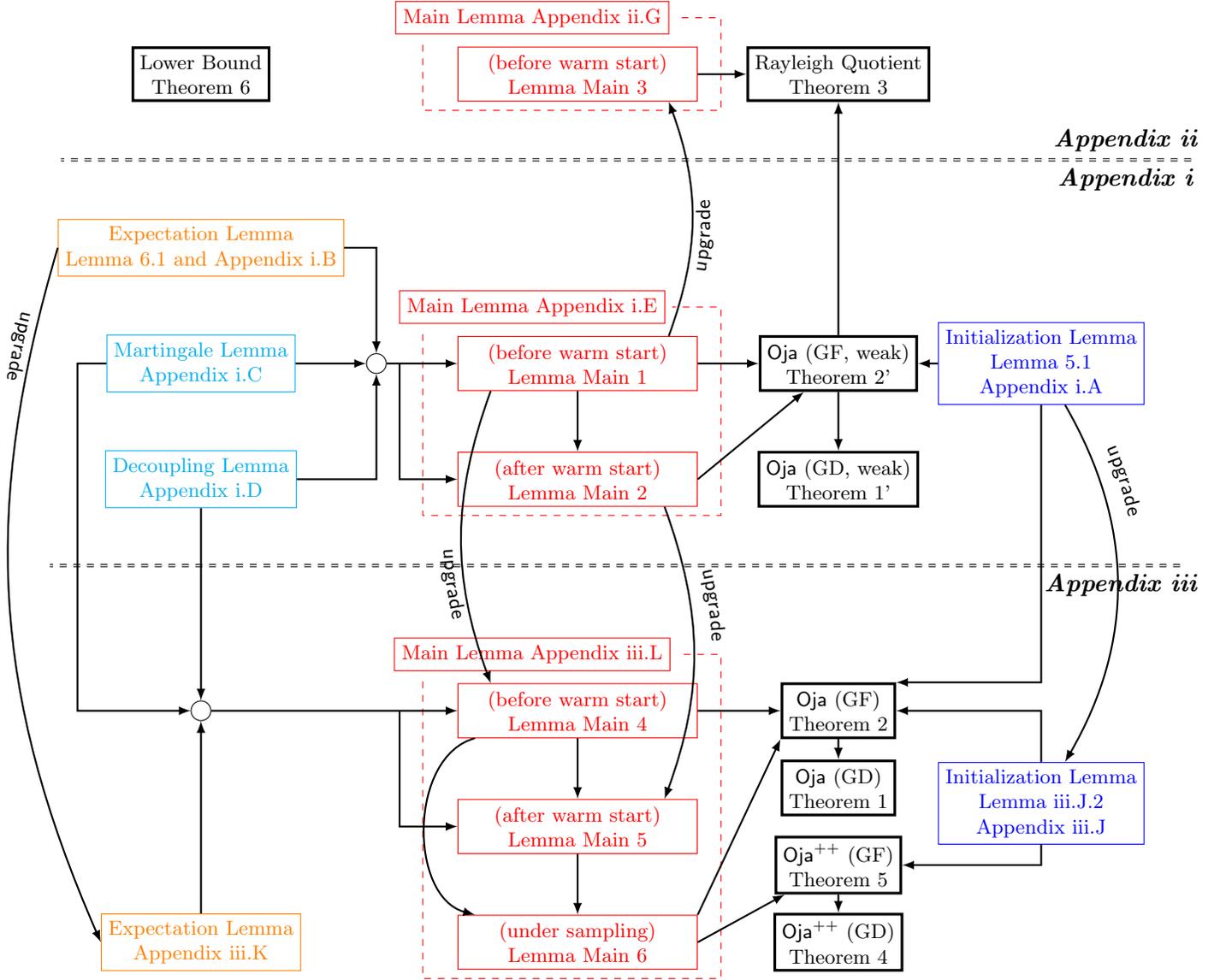

Figure 1: Overall structure of this paper. GF and GD stand for gap-free and gap-dependent respectively.

We divide our appendix sections into three parts, Appendix i, ii, and iii.

- Appendix i (page 16) provides complete proof but for two weaker versions of our Theorem 1 and 2.
  - Appendix i.A and i.B give missing proofs for Section 5 and 6;
  - Appendix i.C and i.D provide details for our martingale and decoupling lemmas;
  - Appendix i.E proves main lemmas in Section 7 and Appendix i.F puts everything together.
- Appendix ii (page 35) includes proofs for Theorem 6 and Theorem 3.
  - Appendix ii.G extends our main lemmas to better serve for the rayleigh quotient setting;

- Appendix ii.H provides the final proof for our Rayleigh Quotient Theorem 3;
- Appendix ii.I includes a three-paged proof of our lower bound Theorem 6.

- Appendix iii (page 42) provide full proofs not only to the stronger Theorem 1 and Theorem 2 for Oja's algorithm, but also to Theorem 4 and Theorem 5 for Oja$^{++}$.

  - Appendix iii.J extends our initialization lemma in Appendix i.A to stronger settings;
  - Appendix iii.K extends our expectation lemmas in Appendix i.B to stronger settings;
  - Appendix iii.L extends our main lemmas in Appendix i.E to stronger settings;
  - Appendix iii.M provides the final proofs for our Theorem 1 and Theorem 2;
  - Appendix iii.N provides the final proofs for our Theorem 4 and Theorem 5.

We include the dependency graphs of all of our main sections, lemmas and theorems in Figure 1 for a quick reference.

# Appendix (Part I)

In this Part I of the appendix, we provide complete proof but two weaker versions of our Theorem 1 and 2. We state these weaker versions Theorem 1' and 2' here, and meanwhile:

- Appendix i.A and i.B give missing proofs for Section 5 and 6;
- Appendix i.C and i.D provide details for our martingale and decoupling lemmas;
- Appendix i.E proves main lemmas in Section 7; and
- Appendix i.F puts everything together and proves Theorem 1' and 2'.

---

**Theorem 1'** (gap-dependent streaming $k$-PCA). *Letting* $\mathsf{gap} \overset{\text{def}}{=} \lambda_k - \lambda_{k+1} \in \big(0, \frac{1}{k}\big]$, *for every* $\varepsilon, p \in (0,1)$ *define learning rates*

$$T_0 = \widetilde{\Theta}\left(\frac{k}{\mathsf{gap}^2 \cdot p^2}\right), \quad \eta_t = \begin{cases} \widetilde{\Theta}\left(\frac{1}{\mathsf{gap} \cdot T_0}\right) & 1 \leq t \leq T_0; \\ \widetilde{\Theta}\left(\frac{1}{\mathsf{gap} \cdot t}\right) & t > T_0. \end{cases}$$

*Let* $\mathbf{Z}$ *be the column orthonormal matrix consisting of all eigenvectors of* $\mathbf{\Sigma}$ *with values no more than* $\lambda_{k+1}$. *Then, the output* $\mathbf{Q}_T \in \mathbb{R}^{d \times k}$ *of Oja's algorithm satisfies with prob. at least* $1 - p$:

$$\text{for every} \quad T = T_0 + \widetilde{\Theta}\left(\frac{T_0}{\varepsilon}\right) \quad \text{it satisfies} \quad \|\mathbf{Z}^\top \mathbf{Q}_T\|_F^2 \leq \varepsilon \ .$$

*Above,* $\widetilde{\Theta}$ *hides poly-log factors in* $\frac{1}{p}, \frac{1}{\mathsf{gap}}$ *and* $d$.

---

**Theorem 2'** (gap-free streaming $k$-PCA). *For every* $\rho, \varepsilon, p \in (0,1)$, *define learning rates*

$$T_0 = \widetilde{\Theta}\left(\frac{k}{\rho^2 \cdot p^2}\right), \quad \eta_t = \begin{cases} \widetilde{\Theta}\left(\frac{1}{\rho \cdot T_0}\right) & t \leq T_0; \\ \widetilde{\Theta}\left(\frac{1}{\rho \cdot t}\right) & t > T_0. \end{cases}$$

*Let* $\mathbf{W}$ *be the column orthonormal matrix consisting of all eigenvectors of* $\mathbf{\Sigma}$ *with values no more than* $\lambda_k - \rho$. *Then, the output* $\mathbf{Q}_T \in \mathbb{R}^{d \times k}$ *of Oja's algorithm satisfies with prob. at least* $1 - p$:

$$\text{for every} \quad T = T_0 + \widetilde{\Theta}\left(\frac{T_0}{\varepsilon}\right) \quad \text{it satisfies} \quad \|\mathbf{W}^\top \mathbf{Q}_T\|_F^2 \leq \varepsilon \ .$$

*Above,* $\widetilde{\Theta}$ *hides poly-log factors in* $\frac{1}{p}, \frac{1}{\rho}$ *and* $d$.

---

## i.A  Random Initialization (for Section 5)

Recall that $\mathbf{Q} \in \mathbb{R}^{d \times k}$ is a matrix with each entry i.i.d drawn from $\mathcal{N}(0,1)$, the standard gaussian.

### i.A.1  Preparation Lemmas

**Lemma i.A.1.** *For every $x \in \mathbb{R}^d$ that has Euclidean norm $\|x\|_2 \leq 1$, every PSD matrix $\mathbf{A} \in \mathbb{R}^{k \times k}$, and every $\lambda \geq 1$, we have*

$$\Pr_{\mathbf{Q}}\left[x^\top \mathbf{ZZ}^\top \mathbf{QAQ}^\top \mathbf{ZZ}^\top x \geq \mathbf{Tr}(\mathbf{A}) + \lambda\right] \leq e^{-\frac{\lambda}{8\mathbf{Tr}(\mathbf{A})}} \ .$$

*Proof of Lemma i.A.1.* Let $\mathbf{A} = \mathbf{U}\boldsymbol{\Sigma}_{\mathbf{A}}\mathbf{U}^\top$ be the eigendecomposition of $\mathbf{A}$, and we denote by $\mathbf{Q}_z = \mathbf{Z}^\top \mathbf{QU} \in \mathbb{R}^{(d-k) \times d}$. Since a random Gaussian matrix is rotation invariant, and since $\mathbf{U}$ is unitary and $\mathbf{Z}$ is column orthonormal, we know that each entry of $\mathbf{Q}_z$ draw i.i.d. from $\mathcal{N}(0,1)$.

Next, since we have $\|\mathbf{Z}^\top x\|_2 \leq 1$, it satisfies that $y = x^\top \mathbf{ZZ}^\top \mathbf{QU}$ is a vector with each coordinate $i$ independently drawn from distribution $\mathcal{N}(0, \sigma_i)$ for $\sigma_i \leq 1$. This implies

$$x^\top \mathbf{ZZ}^\top \mathbf{QAQ}^\top \mathbf{ZZ}^\top x = y^\top \boldsymbol{\Sigma}_{\mathbf{A}} y = \sum_{i=1}^k [\boldsymbol{\Sigma}_{\mathbf{A}}]_{i,i}(y_i)^2 \ .$$

Now, $\sum_{i \in [k]}[\boldsymbol{\Sigma}_{\mathbf{A}}]_{i,i}(y_i)^2$ is a subexponential distribution[14] with parameter $(\sigma^2, b)$ where $\sigma^2, b \leq 4\sum_{i=1}^k [\boldsymbol{\Sigma}_{\mathbf{A}}]_{i,i}$. Using the subexponential concentration bound [20], we have for every $\lambda \geq 1$,

$$\Pr\left[\sum_{i=1}^k [\boldsymbol{\Sigma}_{\mathbf{A}}]_{i,i}(y_i)^2 \geq \sum_{i=1}^k [\boldsymbol{\Sigma}_{\mathbf{A}}]_{i,i} + \lambda\right] \leq \exp\left\{-\frac{\lambda}{8\sum_{i=1}^k [\boldsymbol{\Sigma}_{\mathbf{A}}]_{i,i}}\right\} \ .$$

After rearranging, we have

$$\Pr[x^\top \mathbf{ZZ}^\top \mathbf{QAQ}^\top \mathbf{ZZ}^\top x \geq \mathbf{Tr}(\mathbf{A}) + \lambda] \leq e^{-\frac{\lambda}{8\mathbf{Tr}(\mathbf{A})}} \ . \qquad \square$$

The following lemma is on the singular value distribution of a random Gaussian matrix:

**Lemma i.A.2** (Theorem 1.2 of [19]). *Let $\mathbf{Q} \in \mathbb{R}^{k \times k}$ be a random matrix with each entry i.i.d. drawn from $\mathcal{N}(0,1)$, and $\sigma_1 \leq \sigma_2 \leq \cdots \leq \sigma_k$ be its singular values. We have for every $j \in [k]$ and $\alpha \geq 0$:*

$$\Pr\left[\sigma_j \leq \frac{\alpha j}{\sqrt{k}}\right] \leq \left((2e)^{1/2}\alpha\right)^{j^2} \ .$$

**Lemma i.A.3.** *Let $\mathbf{Q}$ be our initial matrix, then for every $p \in (0,1)$:*

$$\Pr_{\mathbf{Q}}\left[\mathbf{Tr}\left[((\mathbf{V}^\top \mathbf{Q})^\top(\mathbf{V}^\top \mathbf{Q}))^{-1}\right] \geq \frac{\pi^2 ek}{3p}\right] \leq \frac{\sqrt{p}}{1-p} \ .$$

*Proof of Lemma i.A.3.* Using Lemma i.A.2, we know that (using the famous equation $\sum_{j=1}^\infty \frac{1}{j^2} = \frac{\pi^2}{6}$)

$$\Pr\left[\mathbf{Tr}\left[\left((\mathbf{V}^\top \mathbf{Q})^\top(\mathbf{V}^\top \mathbf{Q})\right)^{-1}\right] \geq \frac{\pi^2 ek}{3p}\right] \leq \Pr\left[\exists j \in [k], \sigma_j^{-2}(\mathbf{V}^\top \mathbf{Q}) \geq \frac{2ek}{j^2 p}\right]$$

$$= \Pr\left[\exists j \in [k], \sigma_j(\mathbf{V}^\top \mathbf{Q}) \leq \frac{j\sqrt{p}}{\sqrt{2ek}}\right] \leq \sum_{j=1}^k p^{j^2/2} \leq \frac{\sqrt{p}}{1-p} \ . \qquad \square$$

---

[14]Recall that a random variable $X$ is $(\sigma^2, b)$-subexponential if $\log \mathbb{E}\exp(\lambda(X - \mathbb{E}X)) \leq \lambda^2 \sigma^2/2$ for all $\lambda \in [0, 1/b]$. The squared standard Gaussian variable is $(4, 4)$-subexponential.

### i.A.2 Proof of Lemma 5.1

*Proof of Lemma 5.1.* Applying Lemma i.A.3 with the choice of probability $= \frac{p^2}{4}$, we know that

$$\Pr_{\mathbf{Q}}\left[\mathbf{Tr}(\mathbf{A}) \geq \frac{36k}{p^2}\right] \leq p \quad \text{where} \quad \mathbf{A} \stackrel{\text{def}}{=} \left((\mathbf{V}^\top\mathbf{Q})^\top(\mathbf{V}^\top\mathbf{Q})\right)^{-1} \ .$$

Conditioning on event $\mathcal{C} = \left\{\mathbf{Tr}(\mathbf{A}) \leq \frac{36k}{p^2}\right\}$, and setting $r = \frac{36k}{p^2}$, we have for every fixed $x_1, ..., x_T$ and fixed $i \in [T]$, it satisfies

$$\Pr_{\mathbf{Q}}\left[\left\|x_t^\top \mathbf{Z}\mathbf{Z}^\top (\mathbf{\Sigma}/\lambda_{k+1})^{i-1} \mathbf{Q}(\mathbf{V}^\top\mathbf{Q})^{-1}\right\|_2 \geq \left(18r \ln \frac{T}{q}\right)^{1/2} \bigg| \mathcal{C}, x_t\right]$$

$$\stackrel{\text{①}}{\leq} \Pr\left[\left\|y_t \mathbf{Z}\mathbf{Z}^\top \mathbf{Q}(\mathbf{V}^\top\mathbf{Q})^{-1}\right\|_2 \geq \left(18r \ln \frac{T}{q}\right)^{1/2} \bigg| \mathcal{C}, x_1, ..., x_t\right]$$

$$\stackrel{\text{②}}{\leq} \Pr\left[y_t^\top \mathbf{Z}\mathbf{Z}^\top \mathbf{Q}\mathbf{A}\mathbf{Q}^\top \mathbf{Z}\mathbf{Z}^\top y_t \geq 9r \ln \frac{T^2}{q^2} \bigg| \mathcal{C}, x_1, ..., x_t\right] \stackrel{\text{③}}{\leq} \frac{q^2}{T^2} \ .$$

Above, ① uses the definition $y_t \stackrel{\text{def}}{=} x_t^\top \mathbf{Z}\mathbf{Z}^\top (\mathbf{\Sigma}/\lambda_{k+1})^{i-1}$; ② is from the definition of $\mathbf{A}$; and ③ is owing to Lemma i.A.1 together with the fact that $\|y_t\|_2 \leq \|x_t\|_2 \cdot \|(\frac{\mathbf{Z}\mathbf{Z}^\top \mathbf{\Sigma}}{\lambda_{k+1}})^{i-1}\|_2 \leq 1$ and the fact that $\mathbf{Z}^\top\mathbf{Q}$ is independent of $\mathbf{V}^\top\mathbf{Q}$.[15] Next, define event

$$\mathcal{C}_2 = \left\{\exists i \in [T], \exists t \in [T], \left\|x_t^\top \mathbf{Z}\mathbf{Z}^\top (\mathbf{\Sigma}/\lambda_{k+1})^{i-1} \mathbf{Q}(\mathbf{V}^\top\mathbf{Q})^{-1}\right\|_2 \geq \left(18r \ln \frac{T}{q}\right)^{1/2}\right\} \ .$$

The above derivation, after taking union bound, implies that for every fixed $x_1, ..., x_T$, it satisfies $\mathbf{Pr}_{\mathbf{Q}}[\mathcal{C}_2 \mid \mathcal{C}, x_1, ..., x_T] \leq q^2$. Therefore, denoting by $\mathbb{1}_{\mathcal{C}_2}$ the indicator function of event $\mathcal{C}_2$,

$$\Pr_{\mathbf{Q}}\left[\Pr_{x_1,...,x_T}[\mathcal{C}_2 \mid \mathbf{Q}] \geq q \bigg| \mathcal{C}\right] \leq \frac{1}{q}\mathbb{E}_{\mathbf{Q}}\left[\Pr_{x_1,...,x_T}[\mathcal{C}_2 \mid \mathbf{Q}] \bigg| \mathcal{C}\right]$$

$$= \frac{1}{q}\mathbb{E}_{\mathbf{Q}}\left[\mathbb{E}_{x_1,...,x_T}[\mathbb{1}_{\mathcal{C}_2} \mid \mathbf{Q}] \bigg| \mathcal{C}\right]$$

$$= \frac{1}{q}\mathbb{E}_{x_1,...,x_T}\left[\mathbb{E}_{\mathbf{Q}}[\mathbb{1}_{\mathcal{C}_2} \mid \mathcal{C}, x_1, ..., x_T]\right]$$

$$= \frac{1}{q}\mathbb{E}_{x_1,...,x_T}\left[\Pr_{\mathbf{Q}}[\mathcal{C}_2 \mid \mathcal{C}, x_1, ..., x_T]\right] \leq q \ .$$

Above, the first inequality uses Markov's bound. In an analogous manner, we define event

$$\mathcal{C}_3 = \left\{\exists j \in [d], j \geq k+1, \|\nu_j^\top \mathbf{Q}(\mathbf{V}^\top\mathbf{Q})^{-1}\|_2 \geq \left(18r \ln \frac{d}{p}\right)^{1/2}\right\}$$

where $\nu_j$ is the $j$-th eigenvector of $\mathbf{\Sigma}$ corresponding to eigenvalue $\lambda_j$. A completely analogous proof as the lines above also shows $\mathbf{Pr}_{\mathbf{Q}}[\mathcal{C}_3 \mid \mathcal{C}] \leq q$. Finally, using union bound

$$\Pr_{\mathbf{Q}}\left[\mathcal{C}_3 \bigwedge \Pr_{x_1,...,x_T}[\mathcal{C}_2 \mid \mathbf{Q}] \geq q\right] \leq \Pr_{\mathbf{Q}}[\mathcal{C}_3 \mid \mathcal{C}] + \Pr_{\mathbf{Q}}\left[\Pr_{x_1,...,x_T}[\mathcal{C}_2 \mid \mathbf{Q}] \geq q \bigg| \mathcal{C}\right] + \Pr_{\mathbf{Q}}[\overline{\mathcal{C}}] \leq q + q + p \ ,$$

we conclude that with probability at least $1 - p - 2q$ over the random choice of $\mathbf{Q}$, it satisfies

- $\mathbf{Pr}_{x_1,...,x_T}[\mathcal{C}_2 \mid \mathbf{Q}] < q$, and
- $\overline{\mathcal{C}_3}$ holds (which implies $\|\mathbf{Z}^\top\mathbf{Q}(\mathbf{V}^\top\mathbf{Q})^{-1}\|_F^2 < 18rd \ln \frac{d}{p}$ as desired).

□

---

[15] In principle, we only proved Lemma i.A.1 when $\mathbf{Q}$ is a random matrix, independent of $\mathbf{A}$. Here, $\mathbf{A}$ also depends on $\mathbf{Q}$ but only on $\mathbf{V}^\top\mathbf{Q}$. Therefore, $\mathbf{A}$ is independent from $\mathbf{Z}^\top\mathbf{Q}$, so we can still safely apply Lemma i.A.1.

## i.B  Expectation Lemmas (for Section 6)

Let $\mathbf{X} \in \mathbb{R}^{d \times r}$ be a generic matrix that shall later be chosen as either $\mathbf{X} = \mathbf{W}$, $\mathbf{X} = \mathbf{Z}$, or $\mathbf{X} = [w]$ for some vector $w$. We recall the following notions from Section 6

$$
\begin{array}{ll}
\mathbf{L}_t = \mathbf{P}_t \mathbf{Q}(\mathbf{V}^\top \mathbf{P}_t \mathbf{Q})^{-1} \in \mathbb{R}^{d \times k} & \mathbf{R}'_t = \mathbf{X}^\top x_t x_t^\top \mathbf{L}_{t-1} \in \mathbb{R}^{r \times k} \\
\mathbf{S}_t = \mathbf{X}^\top \mathbf{L}_t \in \mathbb{R}^{r \times k} & \mathbf{H}'_t = \mathbf{V}^\top x_t x_t^\top \mathbf{L}_{t-1} \in \mathbb{R}^{k \times k}
\end{array}
$$

**Lemma i.B.1.** *For every $\mathbf{Q} \in \mathbb{R}^{d \times k}$ and every $t \in [T]$, suppose for $\phi_t \geq 0$, $x_t$ satisfies:*

$$\|x_t^\top \mathbf{L}_{t-1}\|_2 = \|x_t^\top \mathbf{P}_{t-1}\mathbf{Q}(\mathbf{V}^\top \mathbf{P}_{t-1}\mathbf{Q})^{-1}\|_2 \leq \phi_t \quad \text{and} \quad \eta_t \phi_t \leq \frac{1}{2} \ .$$

*Then the following holds:*

(a) $\mathbf{Tr}(\mathbf{S}_t^\top \mathbf{S}_t) \leq \mathbf{Tr}(\mathbf{S}_{t-1}^\top \mathbf{S}_{t-1}) - 2\eta_t \mathbf{Tr}(\mathbf{S}_{t-1}^\top \mathbf{S}_{t-1} \mathbf{H}'_t) + 2\eta_t \mathbf{Tr}(\mathbf{S}_{t-1}^\top \mathbf{R}'_t)$
$\qquad + (12\eta_t^2 \|\mathbf{H}'_t\|_2^2 + 2\eta_t^2 \|\mathbf{R}'_t\|_2 \|\mathbf{H}'_t\|_2) \mathbf{Tr}(\mathbf{S}_{t-1}^\top \mathbf{S}_{t-1}) + 8\eta_t^2 \|\mathbf{R}'_t\|_2^2 + 2\eta_t^2 \|\mathbf{R}'_t\|_2 \|\mathbf{H}'_t\|_2$

(b) $|\mathbf{Tr}(\mathbf{S}_t^\top \mathbf{S}_t) - \mathbf{Tr}(\mathbf{S}_{t-1}^\top \mathbf{S}_{t-1})|^2 \leq 243\eta_t^2 \|\mathbf{H}'_t\|_2^2 \mathbf{Tr}(\mathbf{S}_{t-1}^\top \mathbf{S}_{t-1})^2 + 12\eta_t^2 \|\mathbf{R}'_t\|_2^2 \mathbf{Tr}(\mathbf{S}_{t-1}^\top \mathbf{S}_{t-1}) + 300\eta_t^4 \phi_t^2 \|\mathbf{R}'_t\|_2^2$

(c) $|\mathbf{Tr}(\mathbf{S}_t^\top \mathbf{S}_t) - \mathbf{Tr}(\mathbf{S}_{t-1}^\top \mathbf{S}_{t-1})| \leq 9\eta_t \phi_t \mathbf{Tr}(\mathbf{S}_{t-1}^\top \mathbf{S}_{t-1}) + 2\eta_t \phi_t \sqrt{\mathbf{Tr}(\mathbf{S}_{t-1}^\top \mathbf{S}_{t-1})} + 10\eta_t^2 \phi_t^2$

*Proof of Lemma i.B.1.* We first notice that

$$\mathbf{X}^\top \mathbf{P}_t \mathbf{Q} = \mathbf{X}^\top \mathbf{P}_{t-1} \mathbf{Q} + \eta_t \mathbf{X}^\top x_t x_t^\top \mathbf{P}_{t-1} \mathbf{Q} \quad \text{and}$$
$$\mathbf{V}^\top \mathbf{P}_t \mathbf{Q} = \mathbf{V}^\top \mathbf{P}_{t-1} \mathbf{Q} + \eta_t \mathbf{V}^\top x_t x_t^\top \mathbf{P}_{t-1} \mathbf{Q} \ ,$$

where the second equality further implies (using the Sherman-Morrison formula) that

$$
\begin{aligned}
(\mathbf{V}^\top \mathbf{P}_t \mathbf{Q})^{-1} &= (\mathbf{V}^\top \mathbf{P}_{t-1} \mathbf{Q})^{-1} - \frac{\eta_t (\mathbf{V}^\top \mathbf{P}_{t-1}\mathbf{Q})^{-1} \mathbf{V}^\top x_t x_t^\top \mathbf{P}_{t-1}\mathbf{Q}(\mathbf{V}^\top \mathbf{P}_{t-1}\mathbf{Q})^{-1}}{1 + \eta_t x_t^\top \mathbf{P}_{t-1}\mathbf{Q}(\mathbf{V}^\top \mathbf{P}_{t-1}\mathbf{Q})^{-1} \mathbf{V}^\top x_t} \\
&= (\mathbf{V}^\top \mathbf{P}_{t-1}\mathbf{Q})^{-1} - (\eta_t - \alpha_t \eta_t^2)(\mathbf{V}^\top \mathbf{P}_{t-1}\mathbf{Q})^{-1} \mathbf{H}'_t \ ,
\end{aligned}
$$

and above we denote by $\alpha_t \stackrel{\text{def}}{=} \frac{\psi_t}{1 + \eta_t \psi_t}$ where $\psi_t \stackrel{\text{def}}{=} x_t^\top \mathbf{L}_{t-1} \mathbf{V}^\top x_t$. Therefore, we can write

$$
\begin{aligned}
\mathbf{S}_t &\stackrel{\text{①}}{=} \mathbf{X}^\top \mathbf{P}_t \mathbf{Q}(\mathbf{V}^\top \mathbf{P}_t \mathbf{Q})^{-1} \\
&\stackrel{\text{②}}{=} \mathbf{S}_{t-1} - (\eta_t - \alpha_t \eta_t^2)\mathbf{S}_{t-1}\mathbf{H}'_t + \eta_t \mathbf{R}'_t - (\eta_t^2 - \alpha_t \eta_t^3)\mathbf{R}'_t \mathbf{H}'_t \\
&\stackrel{\text{③}}{=} \mathbf{S}_{t-1} - (\eta_t - \alpha_t \eta_t^2)\mathbf{S}_{t-1}\mathbf{H}'_t + (\eta_t - \psi_t \eta_t^2 + \alpha_t \psi_t \eta_t^3)\mathbf{R}'_t \stackrel{\text{④}}{=} \mathbf{S}_{t-1} - \eta_t \mathbf{S}_{t-1}\mathbf{H}_t + \eta_t \mathbf{R}_t \ .
\end{aligned}
$$

Above, equality ① uses the definition of $\mathbf{S}_t$ and $\mathbf{L}_t$; equality ② uses our derived equations for $\mathbf{X}^\top \mathbf{P}_t \mathbf{Q}$ and $(\mathbf{V}^\top \mathbf{P}_t \mathbf{Q})^{-1}$; equality ③ uses $\mathbf{R}'_t \mathbf{H}'_t = \psi_t \mathbf{R}'_t$; and in quality ④ we have denoted by

$$\mathbf{H}_t = (1 - \alpha_t \eta_t)\mathbf{H}'_t \quad \text{and} \quad \mathbf{R}_t = (1 - \psi_t \eta_t + \alpha_t \psi_t \eta_t^2)\mathbf{R}'_t$$

to simplify the notations. Note that $\mathbf{H}'_t$, $\mathbf{R}'_t$ are rank one matrices so $\|\mathbf{H}'_t\|_F = \|\mathbf{H}'_t\|_2$ and $\|\mathbf{R}'_t\|_F =$

$\|\mathbf{R}'_t\|_2$. We now proceed and compute

$$
\begin{aligned}
\mathbf{Tr}(\mathbf{S}_t^\top \mathbf{S}_t) &= \mathbf{Tr}(\mathbf{S}_{t-1}\mathbf{S}_{t-1}^\top) - 2\eta_t \mathbf{Tr}(\mathbf{S}_{t-1}^\top \mathbf{S}_{t-1}\mathbf{H}_t) + 2\eta_t \mathbf{Tr}(\mathbf{S}_{t-1}^\top \mathbf{R}_t) \\
&\quad + \eta_t^2 \mathbf{Tr}(\mathbf{H}_t^\top \mathbf{S}_{t-1}^\top \mathbf{S}_{t-1}\mathbf{H}_t) + \eta_t^2 \mathbf{Tr}(\mathbf{R}_t^\top \mathbf{R}_t) - 2\eta_t^2 \mathbf{Tr}(\mathbf{R}_t^\top \mathbf{S}_{t-1}\mathbf{H}_t) \\
&\overset{①}{\leq} \mathbf{Tr}(\mathbf{S}_{t-1}\mathbf{S}_{t-1}^\top) - 2\eta_t \mathbf{Tr}(\mathbf{S}_{t-1}^\top \mathbf{S}_{t-1}\mathbf{H}_t) + 2\eta_t \mathbf{Tr}(\mathbf{S}_{t-1}^\top \mathbf{R}_t) \\
&\quad + 2\eta_t^2 \mathbf{Tr}(\mathbf{H}_t^\top \mathbf{S}_{t-1}^\top \mathbf{S}_{t-1}\mathbf{H}_t) + 2\eta_t^2 \mathbf{Tr}(\mathbf{R}_t^\top \mathbf{R}_t) \\
&\overset{②}{\leq} \mathbf{Tr}(\mathbf{S}_{t-1}\mathbf{S}_{t-1}^\top) - 2\eta_t \mathbf{Tr}(\mathbf{S}_{t-1}^\top \mathbf{S}_{t-1}\mathbf{H}_t) + 2\eta_t \mathbf{Tr}(\mathbf{S}_{t-1}^\top \mathbf{R}_t) \\
&\quad + 2\eta_t^2 (1 - \alpha_t \eta_t)^2 \|\mathbf{H}'_t\|_2^2 \mathbf{Tr}(\mathbf{S}_{t-1}\mathbf{S}_{t-1}^\top) + 2\eta_t^2 (1 - \psi_t \eta_t + \alpha_t \psi_t \eta_t^2)^2 \|\mathbf{R}'_t\|_2^2 \\
&\overset{③}{\leq} \mathbf{Tr}(\mathbf{S}_{t-1}\mathbf{S}_{t-1}^\top) - 2\eta_t \mathbf{Tr}(\mathbf{S}_{t-1}^\top \mathbf{S}_{t-1}\mathbf{H}'_t) + 2\eta_t \mathbf{Tr}(\mathbf{S}_{t-1}^\top \mathbf{R}'_t) \\
&\quad + 2\eta_t^2 |\alpha_t| \left|\mathbf{Tr}(\mathbf{S}_{t-1}^\top \mathbf{S}_{t-1}\mathbf{H}'_t)\right| + 2\eta_t (\eta_t |\psi_t| + \eta_t^2 |\alpha_t||\psi_t|) \left|\mathbf{Tr}(\mathbf{S}_{t-1}^\top \mathbf{R}'_t)\right| \\
&\quad + 2\eta_t^2 (1 + 2\phi_t \eta_t)^2 \|\mathbf{H}'_t\|_2^2 \mathbf{Tr}(\mathbf{S}_{t-1}\mathbf{S}_{t-1}^\top) + 2\eta_t^2 (1 + \phi_t \eta_t + 2\phi_t^2 \eta_t^2)^2 \|\mathbf{R}'_t\|_2^2 \\
&\overset{④}{\leq} \mathbf{Tr}(\mathbf{S}_{t-1}\mathbf{S}_{t-1}^\top) - 2\eta_t \mathbf{Tr}(\mathbf{S}_{t-1}^\top \mathbf{S}_{t-1}\mathbf{H}'_t) + 2\eta_t \mathbf{Tr}(\mathbf{S}_{t-1}^\top \mathbf{R}'_t) \\
&\quad + 4\eta_t^2 \|\mathbf{H}'_t\|_2 \left|\mathbf{Tr}(\mathbf{S}_{t-1}^\top \mathbf{S}_{t-1}\mathbf{H}'_t)\right| + 4\eta_t^2 \|\mathbf{H}'_t\|_2 \left|\mathbf{Tr}(\mathbf{S}_{t-1}^\top \mathbf{R}'_t)\right| \\
&\quad + 8\eta_t^2 \|\mathbf{H}'_t\|_2^2 \mathbf{Tr}(\mathbf{S}_{t-1}\mathbf{S}_{t-1}^\top) + 8\eta_t^2 \|\mathbf{R}'_t\|_2^2 \\
&\overset{⑤}{\leq} \mathbf{Tr}(\mathbf{S}_{t-1}\mathbf{S}_{t-1}^\top) - 2\eta_t \mathbf{Tr}(\mathbf{S}_{t-1}^\top \mathbf{S}_{t-1}\mathbf{H}'_t) + 2\eta_t \mathbf{Tr}(\mathbf{S}_{t-1}^\top \mathbf{R}'_t) \\
&\quad + 4\eta_t^2 \|\mathbf{H}'_t\|_2 \left|\mathbf{Tr}(\mathbf{S}_{t-1}^\top \mathbf{R}'_t)\right| + 12\eta_t^2 \|\mathbf{H}'_t\|_2^2 \mathbf{Tr}(\mathbf{S}_{t-1}\mathbf{S}_{t-1}^\top) + 8\eta_t^2 \|\mathbf{R}'_t\|_2^2 \ . \quad \text{(i.B.1)}
\end{aligned}
$$

Above, ① is because $2\mathbf{Tr}(\mathbf{A}^\top \mathbf{B}) \leq \mathbf{Tr}(\mathbf{A}^\top \mathbf{A}) + \mathbf{Tr}(\mathbf{B}^\top \mathbf{B})$ which is Young's inequality in the matrix case; ② and ③ are both because $\mathbf{H}_t = (1 - \alpha_t \eta_t)\mathbf{H}'_t$ and $\mathbf{R}_t = (1 - \psi_t \eta_t + \alpha_t \psi_t \eta_t^2)\mathbf{R}'_t$; ④ follow from the parameter properties $|\psi_t| \leq \|\mathbf{H}'_t\|_2 \leq \phi_t$, $|\alpha_t| \leq 2\|\mathbf{H}'_t\|_2 \leq 2\phi_t$, and $0 \leq \eta_t \phi_t \leq \frac{1}{2}$; ⑤ follows from $|\mathbf{Tr}(\mathbf{S}_{t-1}^\top \mathbf{S}_{t-1}\mathbf{H}'_t)| \leq \mathbf{Tr}(\mathbf{S}_{t-1}^\top \mathbf{S}_{t-1})\|\mathbf{H}'_t\|_2$ which uses Proposition 2.1.

Next, Proposition 2.1 tells us

$$
|\mathbf{Tr}(\mathbf{S}_{t-1}^\top \mathbf{R}'_t)| \leq \|\mathbf{R}'_t\|_{S_1} \|\mathbf{S}_{t-1}\|_2 \leq \|\mathbf{R}'_t\|_2 \sqrt{\mathbf{Tr}(\mathbf{S}_{t-1}^\top \mathbf{S}_{t-1})} \leq \frac{\|\mathbf{R}'_t\|_2}{2} \left(\mathbf{Tr}(\mathbf{S}_{t-1}^\top \mathbf{S}_{t-1}) + 1\right) \ , \quad \text{(i.B.2)}
$$

(the second inequality is because $\mathbf{R}'_t$ is rank 1, and the spectral norm of a matrix is no greater than its Frobenius norm.) we can further simplify the upper bound in (i.B.1) as

$$
\begin{aligned}
\mathbf{Tr}(\mathbf{S}_t^\top \mathbf{S}_t) &\leq \mathbf{Tr}(\mathbf{S}_{t-1}\mathbf{S}_{t-1}^\top) - 2\eta_t \mathbf{Tr}(\mathbf{S}_{t-1}^\top \mathbf{S}_{t-1}\mathbf{H}'_t) + 2\eta_t \mathbf{Tr}(\mathbf{S}_{t-1}^\top \mathbf{R}'_t) \\
&\quad + 2\eta_t^2 \|\mathbf{R}'_t\|_2 \|\mathbf{H}'_t\|_2 \left(\mathbf{Tr}(\mathbf{S}_{t-1}^\top \mathbf{S}_{t-1}) + 1\right) + 12\eta_t^2 \|\mathbf{H}'_t\|_2^2 \mathbf{Tr}(\mathbf{S}_{t-1}\mathbf{S}_{t-1}^\top) + 8\eta_t^2 \|\mathbf{R}'_t\|_2^2 \\
&= \mathbf{Tr}(\mathbf{S}_{t-1}\mathbf{S}_{t-1}^\top) - 2\eta_t \mathbf{Tr}(\mathbf{S}_{t-1}^\top \mathbf{S}_{t-1}\mathbf{H}'_t) + 2\eta_t \mathbf{Tr}(\mathbf{S}_{t-1}^\top \mathbf{R}'_t) \\
&\quad + (12\eta_t^2 \|\mathbf{H}'_t\|_2^2 + 2\eta_t^2 \|\mathbf{R}'_t\|_2 \|\mathbf{H}'_t\|_2) \mathbf{Tr}(\mathbf{S}_{t-1}^\top \mathbf{S}_{t-1}) + 8\eta_t^2 \|\mathbf{R}'_t\|_2^2 + 2\eta_t^2 \|\mathbf{R}'_t\|_2 \|\mathbf{H}'_t\|_2 \ .
\end{aligned}
$$

This finishes the proof of Lemma i.B.1-(a).

A completely symmetric analysis of the above derivation also gives

$$
\begin{aligned}
\mathbf{Tr}(\mathbf{S}_t^\top \mathbf{S}_t) &\geq \mathbf{Tr}(\mathbf{S}_{t-1}\mathbf{S}_{t-1}^\top) - 2\eta_t \mathbf{Tr}(\mathbf{S}_{t-1}^\top \mathbf{S}_{t-1}\mathbf{H}'_t) + 2\eta_t \mathbf{Tr}(\mathbf{S}_{t-1}^\top \mathbf{R}'_t) \\
&\quad - (12\eta_t^2 \|\mathbf{H}'_t\|_2^2 + 2\eta_t^2 \|\mathbf{R}'_t\|_2 \|\mathbf{H}'_t\|_2) \mathbf{Tr}(\mathbf{S}_{t-1}^\top \mathbf{S}_{t-1}) - 8\eta_t^2 \|\mathbf{R}'_t\|_2^2 - 2\eta_t^2 \|\mathbf{R}'_t\|_2 \|\mathbf{H}'_t\|_2 \ ,
\end{aligned}
$$

and thus combining the upper and lower bounds we have

$$|\mathbf{Tr}(\mathbf{S}_t^\top \mathbf{S}_t) - \mathbf{Tr}(\mathbf{S}_{t-1}\mathbf{S}_{t-1}^\top)| \leq 2\eta_t |\mathbf{Tr}(\mathbf{S}_{t-1}^\top \mathbf{S}_{t-1}\mathbf{H}_t')| + 2\eta_t |\mathbf{Tr}(\mathbf{S}_{t-1}^\top \mathbf{R}_t')| \quad \text{(i.B.3)}$$
$$+ (12\eta_t^2 \|\mathbf{H}_t'\|_2^2 + 2\eta_t^2 \|\mathbf{R}_t'\|_2 \|\mathbf{H}_t'\|_2)\mathbf{Tr}(\mathbf{S}_{t-1}^\top \mathbf{S}_{t-1}) + 8\eta_t^2 \|\mathbf{R}_t'\|_2^2 + 2\eta_t^2 \|\mathbf{R}_t'\|_2 \|\mathbf{H}_t'\|_2$$
$$\stackrel{①}{\leq} (2\eta_t \|\mathbf{H}_t'\|_2 + 12\eta_t^2 \|\mathbf{H}_t'\|_2^2 + 2\eta_t^2 \|\mathbf{R}_t'\|_2 \|\mathbf{H}_t'\|_2)\mathbf{Tr}(\mathbf{S}_{t-1}^\top \mathbf{S}_{t-1}) + 2\eta_t \|\mathbf{R}_t'\|_2 \sqrt{\mathbf{Tr}(\mathbf{S}_{t-1}^\top \mathbf{S}_{t-1})} \quad \text{(i.B.4)}$$
$$+ 8\eta_t^2 \|\mathbf{R}_t'\|_2^2 + 2\eta_t^2 \|\mathbf{R}_t'\|_2 \|\mathbf{H}_t'\|_2$$
$$\stackrel{②}{\leq} 9\eta_t \|\mathbf{H}_t'\|_2 \mathbf{Tr}(\mathbf{S}_{t-1}^\top \mathbf{S}_{t-1}) + 2\eta_t \|\mathbf{R}_t'\|_2 \sqrt{\mathbf{Tr}(\mathbf{S}_{t-1}^\top \mathbf{S}_{t-1})} + 10\eta_t^2 \phi_t \|\mathbf{R}_t'\|_2 \ . \quad \text{(i.B.5)}$$

Above, ① again uses Proposition 2.1 and (i.B.2); ② uses $\eta_t \phi_t \leq 1/2$ and $\|\mathbf{H}_t'\|_2, \|\mathbf{R}_t'\|_2 \leq \phi_t$.

Finally, if we take square on both sides of (i.B.5), we have (using again $\eta_t \|\mathbf{R}_t'\|_2 \leq \frac{1}{2}$):

$$|\mathbf{Tr}(\mathbf{S}_t^\top \mathbf{S}_t) - \mathbf{Tr}(\mathbf{S}_{t-1}\mathbf{S}_{t-1}^\top)|^2 \leq 243\eta_t^2 \|\mathbf{H}_t'\|_2^2 \mathbf{Tr}(\mathbf{S}_{t-1}^\top \mathbf{S}_{t-1})^2 + 12\eta_t^2 \|\mathbf{R}_t'\|_2^2 \mathbf{Tr}(\mathbf{S}_{t-1}^\top \mathbf{S}_{t-1}) + 300\eta_t^4 \phi_t^2 \|\mathbf{R}_t'\|_2^2$$

and this finishes the proof of Lemma i.B.1-(b). If we continue to use $\|\mathbf{H}_t'\|_2, \|\mathbf{R}_t'\|_2 \leq \phi_t$ to upper bound the right hand side of (i.B.5), we finish the proof of Lemma i.B.1-(c).
$\square$

*Proof of Lemma 6.1 from Lemma i.B.1.* According to the expectation we have $\mathbb{E}[\mathbf{H}_t' \mid \mathcal{F}_{\leq t-1}, \mathcal{C}_{\leq t}] = \mathbf{V}^\top (\mathbf{\Sigma} + \mathbf{\Delta})\mathbf{L}_{t-1}$ and $\mathbb{E}[\mathbf{R}_t' \mid \mathcal{F}_{\leq t-1}, \mathcal{C}_{\leq t}] = \mathbf{X}^\top (\mathbf{\Sigma} + \mathbf{\Delta})\mathbf{L}_{t-1}$. Now we consider the subcases separately:

(a) By Lemma i.B.1-(a),

$$\mathbb{E}\left[\mathbf{Tr}(\mathbf{S}_t^\top \mathbf{S}_t) \mid \mathcal{F}_{\leq t-1}, \mathcal{C}_{\leq t}\right] \stackrel{①}{\leq} (1 + 14\eta_t^2 \phi_t^2)\mathbf{Tr}(\mathbf{S}_{t-1}\mathbf{S}_{t-1}^\top) + 10\eta_t^2 \phi_t^2$$
$$+ \mathbb{E}\left[-2\eta_t \mathbf{Tr}(\mathbf{S}_{t-1}^\top \mathbf{S}_{t-1}\mathbf{H}_t') + 2\eta_t \mathbf{Tr}(\mathbf{S}_{t-1}^\top \mathbf{R}_t') \mid \mathcal{F}_{\leq t-1}, \mathcal{C}_{\leq t}\right] \ .$$
$$\text{(i.B.6)}$$

Above, ① uses $\|\mathbf{R}_t'\|_2, \|\mathbf{H}_t'\|_2 \leq \phi_t$. Next, we compute the expectation

$$\mathbb{E}\left[-2\eta_t \mathbf{Tr}(\mathbf{S}_{t-1}^\top \mathbf{S}_{t-1}\mathbf{H}_t') + 2\eta_t \mathbf{Tr}(\mathbf{S}_{t-1}^\top \mathbf{R}_t') \mid \mathcal{F}_{\leq t-1}, \mathcal{C}_{\leq t}\right]$$
$$= -2\eta_t \mathbf{Tr}(\mathbf{S}_{t-1}^\top \mathbf{S}_{t-1}\mathbf{V}^\top (\mathbf{\Sigma} + \mathbf{\Delta})\mathbf{L}_{t-1}) + 2\eta_t \mathbf{Tr}(\mathbf{S}_{t-1}^\top \mathbf{Z}^\top (\mathbf{\Sigma} + \mathbf{\Delta})\mathbf{L}_{t-1})$$
$$\stackrel{②}{\leq} -2\eta_t \mathsf{gap} \cdot \mathbf{Tr}(\mathbf{S}_{t-1}\mathbf{S}_{t-1}^\top) - 2\eta_t \mathbf{Tr}(\mathbf{S}_{t-1}^\top \mathbf{S}_{t-1}\mathbf{V}^\top \mathbf{\Delta}\mathbf{L}_{t-1}) + 2\eta_t \mathbf{Tr}(\mathbf{S}_{t-1}^\top \mathbf{Z}^\top \mathbf{\Delta}\mathbf{L}_{t-1}) \ . \text{(i.B.7)}$$

Above, ② is because $\mathbf{Tr}(\mathbf{S}_{t-1}^\top \mathbf{Z}^\top \mathbf{\Sigma}\mathbf{L}_{t-1}) = \mathbf{Tr}(\mathbf{S}_{t-1}^\top \mathbf{\Sigma}_{>k}\mathbf{Z}^\top \mathbf{L}_{t-1}) = \mathbf{Tr}(\mathbf{S}_{t-1}^\top \mathbf{\Sigma}_{>k}\mathbf{S}_{t-1}) \leq \lambda_{k+1}\mathbf{Tr}(\mathbf{S}_{t-1}^\top \mathbf{S}_{t-1})$, as well as $\mathbf{Tr}(\mathbf{S}_{t-1}^\top \mathbf{S}_{t-1}\mathbf{V}^\top \mathbf{\Sigma}\mathbf{L}_{t-1}) = \mathbf{Tr}(\mathbf{S}_{t-1}^\top \mathbf{S}_{t-1}\mathbf{\Sigma}_{\leq k}\mathbf{V}^\top \mathbf{L}_{t-1}) = \mathbf{Tr}(\mathbf{S}_{t-1}^\top \mathbf{S}_{t-1}\mathbf{\Sigma}_{\leq k}) \geq \lambda_k \mathbf{Tr}(\mathbf{S}_{t-1}^\top \mathbf{S}_{t-1})$.

Using the decomposition $\mathbf{I} = \mathbf{V}\mathbf{V}^\top + \mathbf{Z}\mathbf{Z}^\top$, $\|\mathbf{V}\|_2 \leq 1, \|\mathbf{Z}\|_2 \leq 1$, and Proposition 2.1 multiple times, we have

$$\begin{aligned}
\mathbf{Tr}(\mathbf{S}_{t-1}^\top \mathbf{S}_{t-1}\mathbf{V}^\top \mathbf{\Delta}\mathbf{L}_{t-1}) &= \mathbf{Tr}(\mathbf{S}_{t-1}^\top \mathbf{S}_{t-1}\mathbf{V}^\top \mathbf{\Delta}(\mathbf{V}\mathbf{V}^\top + \mathbf{Z}\mathbf{Z}^\top)\mathbf{L}_{t-1}) \\
&\leq \mathbf{Tr}(\mathbf{S}_{t-1}^\top \mathbf{S}_{t-1}\mathbf{V}^\top \mathbf{\Delta}\mathbf{V}) + \mathbf{Tr}(\mathbf{S}_{t-1}^\top \mathbf{S}_{t-1}\mathbf{V}^\top \mathbf{\Delta}\mathbf{Z}\mathbf{S}_{t-1}) \\
&\stackrel{①}{\leq} \|\mathbf{\Delta}\|_2 \left(\mathbf{Tr}(\mathbf{S}_{t-1}^\top \mathbf{S}_{t-1}) + \left[\mathbf{Tr}(\mathbf{S}_{t-1}^\top \mathbf{S}_{t-1})\right]^{3/2}\right) \\
\mathbf{Tr}(\mathbf{S}_{t-1}^\top \mathbf{Z}^\top \mathbf{\Delta}\mathbf{L}_{t-1}) &= \mathbf{Tr}(\mathbf{S}_{t-1}^\top \mathbf{Z}^\top \mathbf{\Delta}(\mathbf{V}\mathbf{V}^\top + \mathbf{Z}\mathbf{Z}^\top)\mathbf{L}_{t-1}) \\
&\leq \mathbf{Tr}(\mathbf{S}_{t-1}^\top \mathbf{Z}^\top \mathbf{\Delta}\mathbf{V}) + \mathbf{Tr}(\mathbf{S}_{t-1}^\top \mathbf{Z}^\top \mathbf{\Delta}\mathbf{Z}\mathbf{S}_{t-1}) \\
&\leq \|\mathbf{\Delta}\|_2 \left(\mathbf{Tr}(\mathbf{S}_{t-1}^\top \mathbf{S}_{t-1}) + \mathbf{Tr}(\mathbf{S}_{t-1}^\top \mathbf{S}_{t-1})^{1/2}\right) \ .
\end{aligned}$$

Above, ① uses the fact that $\|\mathbf{S}_{t-1}\mathbf{S}_{t-1}^\top \mathbf{S}_{t-1}\|_{S_1} \leq \|\mathbf{S}_{t-1}\mathbf{S}_{t-1}^\top\|_{S_1}\|\mathbf{S}_{t-1}\|_2 \leq \left[\mathbf{Tr}(\mathbf{S}_{t-1}^\top \mathbf{S}_{t-1})\right]^{3/2}$

Plugging them into (i.B.7) gives

$$\mathbb{E}\left[-2\eta_t \mathbf{Tr}(\mathbf{S}_{t-1}^\top \mathbf{S}_{t-1}\mathbf{H}_t') + 2\eta_t \mathbf{Tr}(\mathbf{S}_{t-1}^\top \mathbf{R}_t') \mid \mathcal{F}_{\leq t-1}, \mathcal{C}_{\leq t}\right] \leq -2\eta_t \mathsf{gap} \cdot \mathbf{Tr}(\mathbf{S}_{t-1}\mathbf{S}_{t-1}^\top)$$
$$+2\eta_t \|\mathbf{\Delta}\|_2 \left( \left[\mathbf{Tr}(\mathbf{S}_{t-1}^\top \mathbf{S}_{t-1})\right]^{3/2} + 2\mathbf{Tr}(\mathbf{S}_{t-1}^\top \mathbf{S}_{t-1}) + \left[\mathbf{Tr}(\mathbf{S}_{t-1}^\top \mathbf{S}_{t-1})\right]^{1/2} \right) . \quad \text{(i.B.8)}$$

Putting this back to (i.B.6) finishes the proof of Corollary 6.1-(a).

(b) In this case (i.B.7) also holds but one needs to replace $\mathsf{gap}$ with $\rho$ because of the definitional difference between $\mathbf{W}$ and $\mathbf{Z}$. We compute the following upper bounds similar to case (a):

$$\begin{aligned}
\mathbf{Tr}(\mathbf{S}_{t-1}^\top \mathbf{S}_{t-1}\mathbf{V}^\top \mathbf{\Delta}\mathbf{L}_{t-1}) &= \mathbf{Tr}(\mathbf{S}_{t-1}^\top \mathbf{S}_{t-1}\mathbf{V}^\top \mathbf{\Delta}(\mathbf{V}\mathbf{V}^\top + \mathbf{Z}\mathbf{Z}^\top)\mathbf{L}_{t-1}) \\
&\leq \mathbf{Tr}(\mathbf{S}_{t-1}^\top \mathbf{S}_{t-1}\mathbf{V}^\top \mathbf{\Delta}\mathbf{V}) + \mathbf{Tr}(\mathbf{S}_{t-1}^\top \mathbf{S}_{t-1}\mathbf{V}^\top \mathbf{\Delta}\mathbf{Z}\mathbf{Z}^\top \mathbf{L}_{t-1}) \\
&\stackrel{①}{\leq} \|\mathbf{\Delta}\|_2 \mathbf{Tr}(\mathbf{S}_{t-1}^\top \mathbf{S}_{t-1}) \left(1 + \left[\mathbf{Tr}(\mathbf{Z}^\top \mathbf{L}_{t-1}\mathbf{L}_{t-1}^\top \mathbf{Z})\right]^{1/2}\right) \\
\mathbf{Tr}(\mathbf{S}_{t-1}^\top \mathbf{Z}^\top \mathbf{\Delta}\mathbf{L}_{t-1}) &= \mathbf{Tr}(\mathbf{S}_{t-1}^\top \mathbf{Z}^\top \mathbf{\Delta}(\mathbf{V}\mathbf{V}^\top + \mathbf{Z}\mathbf{Z}^\top)\mathbf{L}_{t-1}) \\
&\leq \mathbf{Tr}(\mathbf{S}_{t-1}^\top \mathbf{Z}^\top \mathbf{\Delta}\mathbf{V}) + \mathbf{Tr}(\mathbf{S}_{t-1}^\top \mathbf{Z}^\top \mathbf{\Delta}\mathbf{Z}\mathbf{Z}^\top \mathbf{L}_{t-1}) \\
&\stackrel{②}{\leq} \|\mathbf{\Delta}\|_2 \mathbf{Tr}(\mathbf{S}_{t-1}^\top \mathbf{S}_{t-1})^{1/2} \left(1 + \mathbf{Tr}(\mathbf{Z}^\top \mathbf{L}_{t-1}\mathbf{L}_{t-1}^\top \mathbf{Z})^{1/2}\right) \quad \text{(i.B.9)}
\end{aligned}$$

Above, ① is because (using Proposition 2.1)

$$\mathbf{Tr}(\mathbf{S}_{t-1}^\top \mathbf{S}_{t-1}\mathbf{V}^\top \mathbf{\Delta}\mathbf{Z}\mathbf{Z}^\top \mathbf{L}_{t-1}) \leq \mathbf{Tr}\left((\mathbf{S}_{t-1}^\top \mathbf{S}_{t-1})^2\right)^{1/2} \cdot \left[\mathbf{Tr}(\mathbf{V}^\top \mathbf{\Delta}\mathbf{Z}\mathbf{Z}^\top \mathbf{L}_{t-1}\mathbf{L}_{t-1}^\top \mathbf{Z}\mathbf{Z}^\top \mathbf{\Delta}^\top \mathbf{V})\right]^{1/2}$$
$$\leq \|\mathbf{\Delta}\|_2 \mathbf{Tr}(\mathbf{S}_{t-1}^\top \mathbf{S}_{t-1}) \cdot \left[\mathbf{Tr}(\mathbf{Z}^\top \mathbf{L}_{t-1}\mathbf{L}_{t-1}^\top \mathbf{Z})\right]^{1/2}$$

and ② holds for a similar reason.

Putting these upper bounds into (i.B.7) finishes the proof of Corollary 6.1-(b).

(c) When $\mathbf{X} = [w]$, a slightly different derivation of (i.B.7) gives

$$\mathbb{E}\left[\mathbf{Tr}(\mathbf{S}_t^\top \mathbf{S}_t) \mid \mathcal{F}_{t-1}, \mathcal{C}_{\leq t}\right] \leq (1 - 2\eta_t \lambda_k + 14\eta_t^2 \phi_t^2) \mathbf{Tr}(\mathbf{S}_{t-1}\mathbf{S}_{t-1}^\top) + 10\eta_t^2 \phi_t^2$$
$$- 2\eta_t \mathbf{Tr}(\mathbf{S}_{t-1}^\top \mathbf{S}_{t-1}\mathbf{V}^\top \mathbf{\Delta}\mathbf{L}_{t-1}) + 2\eta_t \mathbf{Tr}(\mathbf{S}_{t-1}^\top w^\top \mathbf{\Delta}\mathbf{L}_{t-1}) + 2\eta_t \mathbf{Tr}(\mathbf{S}_{t-1}^\top w^\top \mathbf{\Sigma}\mathbf{L}_{t-1}) . \quad \text{(i.B.10)}$$

Note that the third and fourth terms can be upper bounded similarly using (i.B.9). As for the fifth term, we have

$$\mathbf{Tr}(\mathbf{S}_{t-1}^\top w^\top \mathbf{\Sigma}\mathbf{L}_{t-1}) \leq \frac{\lambda_k}{2} \mathbf{Tr}(\mathbf{S}_{t-1}^\top \mathbf{S}_{t-1}) + \frac{1}{2\lambda_k} \mathbf{Tr}(w^\top \mathbf{\Sigma}\mathbf{L}_{t-1}\mathbf{L}_{t-1}^\top \mathbf{\Sigma}w)$$

Putting these together, we have:

$$\mathbb{E}\left[\mathbf{Tr}(\mathbf{S}_t^\top \mathbf{S}_t) \mid \mathcal{F}_{\leq t-1}, \mathcal{C}_{\leq t}\right] \leq \left(1 - \eta_t \lambda_k + 14\eta_t^2 \phi_t^2\right) \mathbf{Tr}(\mathbf{S}_{t-1}\mathbf{S}_{t-1}^\top) + 10\eta_t^2 \phi_t^2 + \frac{\eta_t}{\lambda_k} \|w^\top \mathbf{\Sigma}\mathbf{L}_{t-1}\|_2^2$$
$$+ 2\eta_t \|\mathbf{\Delta}\|_2 \left( \left[\mathbf{Tr}(\mathbf{S}_{t-1}^\top \mathbf{S}_{t-1})\right]^{1/2} + \mathbf{Tr}(\mathbf{S}_{t-1}^\top \mathbf{S}_{t-1}) \right) \left(1 + \left[\mathbf{Tr}(\mathbf{Z}^\top \mathbf{L}_{t-1}\mathbf{L}_{t-1}^\top \mathbf{Z})\right]^{1/2}\right) \quad \square$$

## i.C Martingale Concentrations

We prove in the appendix the following two martingale concentration lemmas. Both of them are stated in their most general form for the purpose of this paper. The first lemma is for 1-d martingales and the second is for multi-d martingales.

At a high level, Lemma i.C.1 will only be used to analyze the sequences $s_t$ or $s_t'$ (see Section 3) after warm start — that is, after $t \geq T_0$. Our Lemma i.C.2 can be used to analyze $c_{t,s}$ as well as $s_t$ and $s_t'$ before warm start.

**Lemma i.C.1** (1-d martingale). *Let $\{z_t\}_{t=t_0}^{\infty}$ be a non-negative random process with starting time $t_0 \in \mathbb{N}^*$. Suppose there exists $\delta > 0$, $\kappa \geq 2$, and $\tau_t = \frac{1}{\delta t}$ such that*

$$\forall t \geq t_0: \left\{ \begin{array}{rl} \mathbb{E}[z_{t+1} \mid \mathcal{F}_{\leq t}] & \leq (1 - \delta\tau_t)z_t + \tau_t^2 \\ \mathbb{E}[(z_{t+1} - z_t)^2 \mid \mathcal{F}_{\leq t}] & \leq \tau_t^2 z_t + \kappa^2 \tau_t^4 \\ |z_{t+1} - z_t| & \leq \kappa \tau_t \sqrt{z_t} + \kappa^2 \tau_t^2 \end{array} \right\} \quad \text{(i.C.1)}$$

*If there exists $\phi \geq 36$ satisfying $\frac{t_0}{\ln^2 t_0} \geq 7.5\kappa^2(\phi + 1)$ with $z_{t_0} \leq \frac{\phi \ln^2 t_0}{2\delta^2 t_0}$, we have:*

$$\mathbf{Pr}\left[\exists t \geq t_0, z_t > \frac{(\phi+1)\ln^2 t}{\delta^2 t}\right] \leq \frac{\exp\{-\left(\frac{\phi}{36} - 1\right)\ln t_0\}}{\frac{\phi}{36} - 1} \; .$$

**Lemma i.C.2** (multi-dimensional martingale). *Let $\{z_t\}_{t=0}^T$ be a random process where each $z_t \in \mathbb{R}_{\geq 0}^D$ is $\mathcal{F}_{\leq t}$-measurable. Suppose there exist nonnegative parameters $\{\beta_t, \delta_t, \tau_t\}_{t=0}^{T-1}$ satisfying $\kappa \geq 0$ and $\kappa \tau_t \leq 1/6$ such that, $\forall i \in [D], \forall t \in \{0, 1, \ldots, T-1\}$,*

*(denoting by $[z_t]_i$ is the i-th coordinate of $z_t$ and $[z_t]_{D+1} = 0$)*

$$\left\{ \begin{array}{rl} \mathbb{E}\left[[z_{t+1}]_i \mid \mathcal{F}_{\leq t}\right] & \leq (1 - \beta_t - \delta_t + \tau_t^2)[z_t]_i + \delta_t[z_t]_{i+1} + \tau_t^2 \;, \\ \mathbb{E}\left[|[z_{t+1}]_i - [z_t]_i|^2 \mid \mathcal{F}_{\leq t}\right] & \leq \tau_t^2 \left([z_t]_i^2 + [z_t]_i\right) + \kappa^2 \tau_t^4 \;, \text{and} \\ |[z_{t+1}]_i - [z_t]_i| & \leq \kappa \tau_t \left([z_t]_i + \sqrt{[z_t]_i}\right) + \kappa^2 \tau_t^2 \;. \end{array} \right\} \quad \text{(i.C.2)}$$

*Then, we have: for every $\lambda > 0$, every $p \in \left[1, \min_{s \in [t]}\{\frac{1}{6\kappa \tau_{s-1}}\}\right]$:*

$$\mathbf{Pr}\left[[z_t]_1 \geq \lambda\right] \leq \lambda^{-p}\Big(\max_{j \in [t+1]}\{[z_0]_j^p\} \exp\left\{\sum_{s=0}^{t-1} 5p^2\tau_s^2 - p\beta_s\right\}$$

$$+ 1.4 \sum_{s=0}^{t-1} \exp\left\{\sum_{u=s+1}^{t-1} 5p^2\tau_u^2 - p\beta_u\right\}\Big) \;.$$

The above two lemmas are stated in the most general way in order to be used towards all of our three theorems each requiring different parameter choices of $\beta_t, \delta_t, \tau_t, \kappa$. For instance, to prove Theorem 2 it suffices to use $\kappa = O(1)$.

### i.C.1 Martingale Corollaries

We provide below four instantiations of these lemmas, each of them can be verified by plugging in the specific parameters.

**Corollary i.C.3** (1-d martingale). *Consider the same setting as Lemma i.C.1. Suppose $p \in (0, \frac{1}{e^2})$, $\delta \leq \frac{1}{\sqrt{8}}$, $\tau_t = \frac{1}{\delta t}$, $\kappa \in \left[2, \frac{1}{\sqrt{2}\delta}\right]$, $\frac{t_0}{\ln^2 t_0} \geq \frac{9\ln(1/p)}{\delta^2}$, and $z_{t_0} \leq 2$ we have:*

$$\mathbf{Pr}\left[\exists t \geq t_0, z_t > \frac{5(t_0/\ln^2 t_0)}{t/\ln^2 t}\right] \leq p \;.$$

**Corollary i.C.4** (multi-d martingale). *Consider the same setting as Lemma i.C.2. Suppose $q \in (0, 1)$, $\min_{s \in [t]}\{\frac{1}{6\kappa \tau_{s-1}}\} \geq 4\ln \frac{4t}{q}$ and $\sum_{s=0}^{t-1} \tau_s^2 \leq \frac{1}{100}\ln^{-1} \frac{4t}{q}$, then*

$$\mathbf{Pr}\left[[z_t]_1 \geq 2\max\left\{1, \max_{j \in [t+1]}\{[z_0]_j\}\right\}\right] \leq q \;.$$

**Corollary i.C.5** (multi-d martingale). *Consider the same setting as Lemma i.C.2. Given $q \in (0, 1)$, suppose there exists parameter $\gamma \geq 1$ such that, denoting by $l \stackrel{\text{def}}{=} 10\gamma \ln \frac{3t}{q}$,*

$$\sum_{s=0}^{t-1} \beta_s - l\tau_s^2 \geq \ln\left(\max_{j \in [t+1]}\{[z_0]_j\}\right) \quad \text{and} \quad \forall s \in \{0, 1, \ldots, t-1\} : \beta_s \geq l\tau_s^2 \bigwedge \kappa\tau_s \leq \frac{1}{12\ln \frac{3t}{q}} \;.$$

*Then, we have* $\quad \mathbf{Pr}\left[[z_t]_1 \geq 2/\gamma\right] \leq q \;.$

### i.C.2 Proofs for One-Dimensional Martingale

*Proof of Lemma i.C.1.* Define $y_t = \frac{\delta^2 t z_t}{\ln t} - \ln t$, then we have:

$$\begin{aligned}
\mathbb{E}[y_{t+1} \mid \mathcal{F}_{\leq t}] &= \frac{\delta^2(t+1)\mathbb{E}[z_{t+1} \mid \mathcal{F}_{\leq t}]}{\ln(t+1)} - \ln(t+1) \\
&\leq \frac{\delta^2(t+1)(1-\delta\tau_t)z_t}{\ln(t+1)} + \frac{\delta^2(t+1)\tau_t^2}{\ln(t+1)} - \ln(t+1) \\
&\leq \frac{\delta^2(t+1)\left(1-\frac{1}{t}\right)}{\ln(t+1)}z_t + \frac{t+1}{t^2\ln(t+1)} - \ln(t+1) \stackrel{①}{\leq} \frac{\delta^2 t z_t}{\ln t} - \ln t = y_t ,
\end{aligned}$$

where ① is because for every $t \geq 4$ it satisfies $\frac{(t+1)(t-1)}{\ln(t+1)} \leq \frac{t^2}{\ln t}$ and $\frac{t+1}{t^2 \ln(t+1)} \leq \ln\left(1 + \frac{1}{t}\right)$.

At the same time, we have

$$|y_{t+1} - y_t| \stackrel{②}{\leq} \frac{\delta^2 t}{\ln t}|z_{t+1} - z_t| + \frac{\delta^2}{\ln t} z_{t+1} + \frac{1}{t} , \tag{i.C.3}$$

where ② is because for every $t \geq 3$ it satisfies $0 \leq \frac{t+1}{\ln(t+1)} - \frac{t}{\ln t} \leq \frac{1}{\ln t}$ and $\ln(t+1) - \ln(t) \leq 1/t$.

Taking square on both sides, we have

$$|y_{t+1} - y_t|^2 \leq 3\left(\frac{\delta^2 t}{\ln t}\right)^2 |z_{t+1} - z_t|^2 + 3\left(\frac{\delta^2}{\ln t}\right)^2 z_{t+1}^2 + \frac{3}{t^2} .$$

Taking expectation on both sides, we have

$$\begin{aligned}
\mathbb{E}[|y_{t+1} - y_t|^2 \mid \mathcal{F}_{\leq t}] &\leq 3\left(\frac{\delta^2 t}{\ln t}\right)^2 (\tau_t^2 z_t + \kappa^2 \tau_t^4) + 3\frac{(y_t + \ln t)^2}{t^2} + \frac{3}{t^2} \\
&< \frac{3(y_t + \ln t)}{t \ln t} + \frac{3(y_t + \ln t)^2}{t^2} + \frac{3(1+\kappa^2)}{t^2} \\
&\stackrel{③}{\leq} \frac{3(\phi+1)}{t} + \frac{3(\phi+1)^2 \ln^2 t}{t^2} + \frac{15\kappa^2}{4t^2} \stackrel{④}{\leq} \frac{4(\phi_t + 1)}{t} .
\end{aligned}$$

Above, ③ uses $y_t \leq \phi \ln t$ and $\kappa \geq 2$; ④ uses $\frac{t}{\ln^2 t} \geq \frac{t_0}{\ln^2 t_0} \geq \max\{7.5\kappa^2, 6(\phi+1)\}$ and $\ln t \geq 1$.

Therefore, if $y_t \leq \phi \ln t$ holds true for $t = t_0, ..., T$ and $t_0 \geq 8$ (which implies $\frac{t}{\ln^2 t} \geq \frac{t_0}{\ln^2 t_0}$), then

$$\sum_{t=t_0}^{T} \mathbb{E}[|y_{t+1} - y_t|^2 \mid \mathcal{F}_{\leq t}] \leq \sum_{t=t_0}^{T} \frac{4(\phi+1)}{t} \leq 4(\phi+1) \int_{t=t_0-1}^{T} \frac{dt}{t} \leq 4(\phi+1)\ln(T) .$$

Now we can check about the absolute difference. We continue from (i.C.3) and derive that, if $y_t \leq \phi \ln t$, then

$$\begin{aligned}
|y_{t+1} - y_t| &\leq \frac{\delta^2 t}{\ln t}|z_{t+1} - z_t| + \frac{\delta^2}{\ln t}z_{t+1} + \frac{1}{t} \leq \frac{\delta^2 t}{\ln t}\left(\kappa \tau_t \sqrt{z_t} + \kappa^2 \tau_t^2\right) + \frac{\delta^2}{\ln t}z_{t+1} + \frac{1}{t} \\
&\leq \kappa\left(\sqrt{\frac{y_t + \ln t}{t \ln t}} + \frac{\kappa}{t \ln t} + \frac{(y_t + \ln t)}{t} + \frac{1}{t}\right) \stackrel{⑤}{\leq} \kappa\left(\sqrt{\frac{y_t + \ln t}{t \ln t}} + \frac{y_t + \ln t + \kappa}{t}\right) \\
&\stackrel{⑥}{\leq} \kappa\left(\sqrt{\frac{(\phi+1)}{t}} + \frac{(\phi+1)\ln t + \kappa}{t}\right) \stackrel{⑦}{\leq} 2\kappa\sqrt{\frac{(\phi+1)}{t}}
\end{aligned}$$

where ⑤ uses $\ln t \geq 2$ and $\kappa \geq 2$, ⑥ uses $y_t \leq \phi \ln t$, and ⑦ uses $\frac{t}{\ln^2 t} \geq \frac{t_0}{\ln^2 t_0} \geq 4\max\{\phi+1, \kappa\}$.

From the above inequality, we have that if $t_0 \geq 4\kappa^2(\phi+1)$ and $y_t \leq \phi \ln t$ holds true for $t = t_0, ..., T-1$ then $|y_{t+1} - y_t| \leq 1$ for all $t = t_0, \ldots, T-1$.

Finally, since we have assumed $\phi > 36$ and $z_{t_0} \leq \frac{\phi \ln^2 t_0}{2\delta^2 t_0}$ which implies $y_{t_0} \leq \frac{\phi \ln t_0}{2}$, we can apply martingale concentration inequality (c.f. [5, Theorem 18]):

$$
\begin{aligned}
\mathbf{Pr}\left[\exists t \geq t_0, y_t > \phi \ln t\right] &\leq \sum_{T=t_0+1}^{\infty} \mathbf{Pr}\left[y_T > \phi \ln T; \forall t \in \{t_0, ..., T-1\}, y_t \leq \phi \ln t\right] \\
&\leq \sum_{T=t_0+1}^{\infty} \mathbf{Pr}\left[y_T - y_{t_0} > \phi \ln T/2; \forall t \in \{t_0, ..., T-1\}, y_t \leq \phi \ln t\right] \\
&\leq \sum_{T=t_0+1}^{\infty} \exp\left\{\frac{-(\phi \ln T/2)^2}{2 \cdot 4(\phi+1)\ln(T-1) + \frac{2}{3}(\phi \ln T/2)}\right\} \\
&\leq \sum_{T=t_0+1}^{\infty} \exp\left\{-\ln T \frac{\phi^2/4}{8(\phi+1) + \phi/3}\right\} \\
&\leq \int_{T=t_0}^{\infty} \exp\left\{-\frac{\phi}{36}\ln T\right\} dT \quad \leq \quad \frac{\exp\{-(\frac{\phi}{36}-1)\ln t_0\}}{\frac{\phi}{36}-1} \ .
\end{aligned}
$$

□

*Proof of Corollary i.C.3.* Define $\phi \stackrel{\text{def}}{=} \frac{4\delta^2 t_0}{\ln^2 t_0} \geq 36 \ln \frac{1}{p} \geq 72$. It is easy to verify that $\frac{t_0}{\ln^2 t_0} \geq 7.5\kappa^2(\phi+1)$ (because $\kappa \leq 1/(\sqrt{2}\delta)$) and $z_{t_0} \leq \frac{\phi \ln^2 t_0}{2\delta^2 t_0} = 2$, so we can apply Lemma i.C.1:

$$
\mathbf{Pr}\left[\exists t \geq t_0, z_t > \frac{(\phi+1)\ln^2 t}{\delta^2 t}\right] \leq \frac{\exp\left\{-\left(\frac{\phi}{36}-1\right)\ln t_0\right\}}{\frac{\phi}{36}-1} \leq \exp\left\{-\left(\frac{\phi}{36}-1\right)\ln t_0\right\} \leq p \ ,
$$

where the last inequality uses $\ln t_0 \geq 2$ and $\left(\frac{\phi}{36}-1\right)\ln t_0 \geq \frac{\phi}{36}$. Therefore, we conclude that

$$
\mathbf{Pr}\left[\exists t \geq t_0, z_t > \frac{5(t_0/\ln^2 t_0)}{t/\ln^2 t}\right] \leq \mathbf{Pr}\left[\exists t \geq t_0, z_t > \frac{(\phi+1)\ln^2 t}{\delta^2 t}\right] \leq p \ .
$$

□

### i.C.3 Proofs for Multi-Dimensional Martingale

*Proof of Corollary i.C.4.* We apply Lemma i.C.2 with $\lambda = 2\max\left\{1, \max_{j \in [t+1]}\{[z_0]_j\}\right\} \geq 2$. Using the fact that $\beta_t \geq 0$, we have

$$
\mathbf{Pr}\left[[z_t]_1 \geq \lambda\right] = \mathbf{Pr}\left[[z_t]_1 \geq 2(\max_{j \in [t+1]}\{[z_0]_j\}+1)\right] \leq (1+1.4t)\exp\left\{-\ln(2^p) + 5p^2 \sum_{s=0}^{t-1} \tau_s^2\right\} \ .
$$

We can take $p = 4 \ln \frac{4t}{q} \leq \min_{s \in [t]}\{\frac{1}{6\kappa \tau_{s-1}}\}$ which satisfies the assumption of Lemma i.C.2. Therefore, denoting by $\alpha = \sum_{s=0}^{t-1} \tau_s^2$, we have

$$
\mathbf{Pr}\left[[z_t]_1 \geq \lambda\right] \leq 4t \exp\left\{-p \ln 2 + 5p^2 \alpha\right\} \leq q \ .
$$

Above, the last inequality is due to $-p \ln 2 + 5p^2 \alpha^2 \leq -2 \ln \frac{4t}{q} + 5\left(4 \ln \frac{4t}{q}\right)^2 \alpha \leq -\ln \frac{4t}{q}$ which holds for every $\alpha \leq \frac{1}{100} \ln^{-1} \frac{4t}{q}$.

□

*Proof of Corollary i.C.5.* We consider fixed $p = \frac{l}{5\gamma} = 2 \ln \frac{3t}{q}$. Let $y_t = \gamma \cdot z_t$, then $y_t$ satisfies (i.C.2) with (using the fact that $\gamma \geq 1$)

$$
\beta'_t = \beta_t, \quad \delta'_t = \delta_t, \quad (\tau'_t)^2 = \gamma \tau_t^2, \quad \kappa' = \kappa \ .
$$

We denote by $b \stackrel{\text{def}}{=} \sum_{s=0}^{t-1} \beta_s = b$ and $a \stackrel{\text{def}}{=} \sum_{s=0}^{t-1} \tau_s^2$, and apply Lemma i.C.2 on $y_t$ with $\lambda = 2$. Using the fact that $\beta_s \geq l\tau_s^2 = 5zp\tau_s^2$ we know $p\beta'_t \geq 5p^2(\tau'_t)^2$. Therefore, for all $s \in \{0, 1, \ldots, T-1\}$ we

have
$$\mathbf{Pr}\left[[y_t]_1 \geq 2\right] \leq \exp\left\{-pb + 5p^2\gamma a + p\ln\Xi - p\ln 2\right\} + 1.4t\exp\left\{-p\ln 2\right\} , \quad \text{(i.C.4)}$$
where we have denoted by $\Xi \stackrel{\text{def}}{=} \max_{j\in[t+1]}\{[z_0]_j\}$ for notational simplicity. Now, the choice $p = 2\ln\frac{3t}{q}$ satisfies the presumption of Lemma i.C.2 because we have assumed $\kappa\tau_s \leq \frac{1}{12\ln\frac{3t}{q}}$. Therefore, we have
$$-pb + 5p^2\gamma a + p\ln\Xi - p\ln 2 = p(-b + la + \ln\Xi - \ln 2) \leq \ln\frac{q}{2} \Longleftarrow b - la \geq \ln\Xi \bigwedge p \geq 2\ln\frac{3t}{q}$$
$$-p\ln 2 \leq \ln\frac{q}{3t} \Longleftarrow p \geq 2\ln\frac{3t}{q} .$$
Plugging them into (i.C.4) gives $\mathbf{Pr}\left[[z_t]_1 \geq \frac{2}{\gamma}\right] = \mathbf{Pr}\left[[y_t]_1 \geq 2\right] \leq \frac{q}{2} + \frac{q}{2} = q$ . $\square$

*Proof of Lemma i.C.2.* Define vector $s_t$ for every $t \in \{0, 1, \ldots, T-1\}$ and $i \in [D]$, it satisfies $[s_t]_i \stackrel{\text{def}}{=} \frac{[z_{t+1}]_i}{[z_t]_i} - 1$. We have
$$\mathbb{E}\left[[s_t]_i \mid \mathcal{F}_{\leq t}\right] \leq -(\delta_t + \beta_t - \tau_t^2) + \delta_t\frac{[z_t]_{i+1}}{[z_t]_i} + \frac{\tau_t^2}{[z_t]_i} . \quad \text{(i.C.5)}$$
In particular,
$$\text{if } [z_t]_i \geq 1, \text{ then } \quad \mathbb{E}\left[[s_t]_i^2 \mid \mathcal{F}_{\leq t}\right] \leq \tau_t^2 + \frac{\tau_t^2}{[z_t]_i} + \frac{\kappa^2\tau_t^4}{[z_t]_i^2} \leq (2 + (\tau_t\kappa)^2)\tau_t^2 \leq 3\tau_t^2 , \quad \text{(i.C.6)}$$
$$\left|[s_t]_i\right| \leq \kappa\tau_t + \frac{\kappa\tau_t}{\sqrt{[z_t]_i}} + \frac{\kappa^2\tau_t^2}{[z_t]_i} \leq \kappa\tau_t(2 + \kappa\tau_t) \leq 3\kappa\tau_t . \quad \text{(i.C.7)}$$
We consider $[z_{t+1}]_i^p$ for some fixed value $p \geq 1$ and derive that (using (i.C.7))
$$\text{if } (\kappa\tau_t)p \leq \frac{1}{6} \text{ and } [z_t]_i \geq 1, \text{ then } \quad [z_{t+1}]_i^p = [z_t]_i^p(1 + [s_t]_i)^p = [z_t]_i^p\left(\sum_{q=0}^{p}\binom{p}{q}[s_t]_i^q\right)$$
$$\leq [z_t]_i^p\left(1 + p[s_t]_i + p^2[s_t]_i^2\right) .$$
After taking expectation, we have if $(\kappa\tau_t)p \leq \frac{1}{6}$ and $[z_t]_i \geq 1$, then
$$\mathbb{E}\left[[z_{t+1}]_i^p \mid \mathcal{F}_{\leq t}\right] \stackrel{①}{\leq} [z_t]_i^p\left(1 + p\mathbb{E}\left[[s_t]_i \mid \mathcal{F}_{\leq t}\right] + 3p^2\tau_t^2\right)$$
$$\stackrel{②}{\leq} [z_t]_i^p\left(1 - p(\delta_t + \beta_t - \tau_t^2) + \delta_t p\frac{[z_t]_{i+1}}{[z_t]_i} + \frac{\tau_t^2 p}{[z_t]_i} + 3p^2\tau_t^2\right)$$
$$= [z_t]_i^p\left(1 - p(\delta_t + \beta_t - \tau_t^2) + 3p^2\tau_t^2\right) + \delta_t p[z_t]_i^{p-1}[z_t]_{i+1} + p\tau_t^2[z_t]_i^{p-1}$$
$$\stackrel{③}{\leq} [z_t]_i^p\left(1 - p(\delta_t + \beta_t - \tau_t^2) + 3p^2\tau_t^2 + p\tau_t^2\right) + \delta_t p\left(\frac{p-1}{p}[z_t]_i^p + \frac{1}{p}[z_t]_{i+1}^p\right)$$
$$= [z_t]_i^p\left(1 - \delta_t - p\beta_t + p\tau_t^2 + 3p^2\tau_t^2 + p\tau_t^2\right) + \delta_t[z_t]_{i+1}^p$$
$$\stackrel{④}{\leq} [z_t]_i^p\left(1 - \delta_t - p\beta_t + 5p^2\tau_t^2\right) + \delta_t[z_t]_{i+1}^p .$$
Above, ① uses (i.C.6); ② uses (i.C.5); ③ uses $[z_t]_i \geq 1$ and Young's inequality $ab \leq a^p/p + b^q/q$ for $1/p + 1/q = 1$; and ④ uses $p \geq 1$.

On the other hand, if $(\kappa\tau_t)p \leq \frac{1}{6}$ but $[z_t]_i < 1$, we have the following simple bound (using $\kappa\tau_t \leq 1/6$):
$$[z_{t+1}]_i \leq (1 + \kappa\tau_t)[z_t]_i + \kappa\tau_t\sqrt{[z_t]_i} + \kappa^2\tau_t^2 \leq (1 + \kappa\tau_t) + (\kappa\tau_t) + \kappa^2\tau_t^2 < 1.4 .$$
Therefore, as long as $(\kappa\tau_t)p \leq \frac{1}{6}$ we always have
$$\mathbb{E}\left[[z_{t+1}]_i^p \mid \mathcal{F}_{\leq t}\right] \leq [z_t]_i^p\left(1 - \delta_t - p\beta_t + 5p^2\tau_t^2\right) + \delta_t[z_t]_{i+1}^p + 1.4 =: (1 - \alpha_t)[z_t]_i^p + \delta_t[z_t]_{i+1}^p + 1.4 ,$$

and in the last inequality we have denoted by $\alpha_t \stackrel{\text{def}}{=} \delta_t + p\beta_t - 5p^2\tau_t^2$. Telescoping this expectation, and choosing $i = 1$, we have whenever $p \in [1, \min_{s \in [t]}\{\frac{1}{6\kappa\tau_{s-1}}\}]$, it satisfies

$$\begin{aligned}
\mathbb{E}\left[[z_{t+1}]_1^p\right] &\leq \prod_{s=1}^{t}(1 - \alpha_s + \delta_s)\left(\max_{j \in [t+2]}\{[z_0]_j^p\}\right) + 1.4\sum_{s=0}^{t}\left(\prod_{u=s+1}^{t}(1 - \alpha_u + \delta_u)\right) \\
&\leq \prod_{s=0}^{t}(1 - p\beta_s + 5p^2\tau_s^2)\left(\max_{j \in [t+2]}\{[z_0]_j^p\}\right) + 1.4\sum_{s=0}^{t}\left(\prod_{u=s+1}^{t}(1 - p\beta_u + 5p^2\tau_u^2)\right) \\
&\leq \max_{j \in [t+2]}\{[z_0]_j^p\}\exp\left\{-p\left(\sum_{s=0}^{t}\beta_s\right) + 5p^2\sum_{s=0}^{t}\tau_s^2\right\} + 1.4\sum_{s=0}^{t}\exp\left\{-p\left(\sum_{u=s+1}^{t}\beta_u\right) + 5p^2\sum_{u=s+1}^{t}\tau_u^2\right\} .
\end{aligned}$$

Finally, using Markov's inequality, we have for every $\lambda > 0$:

$$\begin{aligned}
\mathbf{Pr}\left[[z_{t+1}]_1 \geq \lambda\right] &\leq \lambda^{-p}\Big(\max_{j \in [t+2]}\{[z_0]_j^p\}\exp\left\{\sum_{s=0}^{t}5p^2\tau_s^2 - p\beta_s\right\} \\
&\quad + 1.4\sum_{s=0}^{t}\exp\left\{\sum_{u=s+1}^{t}5p^2\tau_u^2 - p\beta_u\right\}\Big) . \qquad \square
\end{aligned}$$

## i.D Decoupling Lemmas

We prove the following general lemma. Let $x_1, ..., x_T \in \Omega$ be random variables each i.i.d. drawn from some distribution $\mathcal{D}$. Let $\mathcal{F}_t$ be the sigma-algebra generated by $x_t$, and denote by $\mathcal{F}_{\leq t} = \vee_{s=1}^{t}\mathcal{F}_t$.[16]

**Lemma i.D.1** (decoupling lemma). *Consider a fixed value $q \in [0, 1)$. For every $t \in [T]$ and $s \in \{0, 1, ..., t-1\}$, let $y_{t,s} \in \mathbb{R}^D$ be an $\mathcal{F}_t \vee \mathcal{F}_{\leq s}$ measurable random vector and let $\phi_{t,s} \in \mathbb{R}^D$ be a fixed vector. Let $D' \in [D]$. Define events (we denote by $^{(i)}$ the i-th coordinate)*

$$\mathcal{C}'_t \stackrel{\text{def}}{=} \left\{(x_1, ..., x_{t-1}) \text{ satisfies } \mathbf{Pr}_{x_t}\left[\exists i \in [D']: y_{t,t-1}^{(i)} > \phi_{t,t-1}^{(i)} \mid \mathcal{F}_{t-1}\right] \leq q\right\}$$

$$\mathcal{C}''_t \stackrel{\text{def}}{=} \left\{(x_1, ..., x_t) \text{ satisfies } \forall i \in [D']: y_{t,t-1}^{(i)} \leq \phi_{t,t-1}^{(i)}\right\}$$

*and denote by $\mathcal{C}_t \stackrel{\text{def}}{=} \mathcal{C}'_t \wedge \mathcal{C}''_t$ and $\mathcal{C}_{\leq t} \stackrel{\text{def}}{=} \bigwedge_{s=1}^{t}\mathcal{C}_s$. Suppose the following three assumptions hold:*

(A1) *The random process $\{y_{t,s}\}_{t,s}$ satisfy that for every $i \in [D], t \in [T-1], s \in \{0, 1, ..., t-2\}$*

(a) $\mathbb{E}\left[y_{t,s+1}^{(i)} \mid \mathcal{F}_t, \mathcal{F}_{\leq s}, \mathcal{C}_{\leq s}\right] \leq f_s^{(i)}(y_{t,s}, q)$,

(b) $\mathbb{E}\left[|y_{t,s+1}^{(i)} - y_{t,s}^{(i)}|^2 \mid \mathcal{F}_t, \mathcal{F}_{\leq s}, \mathcal{C}_{\leq s}\right] \leq h_s^{(i)}(y_{t,s}, q)$, and

(c) $\left|y_{t,s+1}^{(i)} - y_{t,s}^{(i)}\right| \leq g_s^{(i)}(y_{t,s})$ whenever $\mathcal{C}_{\leq s}$ holds.

*Above, for each $i \in [D]$ and $s \in \{0, 1, ..., T-2\}$, we have $f_s, h_s : \mathbb{R}^d \times [0, 1] \to \mathbb{R}_{\geq 0}^D$, $g_s : \mathbb{R}^d \to \mathbb{R}_{\geq 0}^D$ are functions satisfying for every $x \in \mathbb{R}^d$,*

(d) $f_s^{(i)}(x, p), h_s^{(i)}(x, p)$ are monotone increasing in $p$, and

(e) $\left|x^{(i)} - f_s^{(i)}(x, 0)\right|^2 \leq h_s^{(i)}(x, 0)$ and $\left|x^{(i)} - f_s^{(i)}(x, 0)\right| \leq g_s^{(i)}(x)$ whenever $f_s^{(i)}(x, 0) \leq x^{(i)}$.

(A2) *Each $t \in [T]$ satisfies $\mathbf{Pr}_{x_t}[\overline{\mathcal{E}_t}] \leq q^2/2$ where event*

$$\mathcal{E}_t \stackrel{\text{def}}{=} \left\{x_t \text{ satisfies } \forall i \in [D]: y_{t,0}^{(i)} \leq \phi_{t,0}^{(i)}\right\} .$$

---
[16]For the purpose of this paper, one can feel free view $\Omega$ as $\mathbb{R}^d$, each $x_t$ as the $t$-th sample vector, and $\mathcal{D}$ as the distribution with covariance matrix $\Sigma$.

(A3) For every $t \in [T]$, letting $x_t$ be any vector satisfying $\mathcal{E}_t$, consider **any** random process $\{z_s\}_{s=0}^{t-1}$ where each $z_s \in \mathbb{R}_{\geq 0}^D$ is $\mathcal{F}_{\leq s}$ measurable with $z_0 = y_{t,0}$ as the starting vector. Suppose that **whenever** $\{z_s\}_{s=0}^{t-1}$ satisfies

$$\forall i \in [D], \forall s \in \{0, 1, \ldots, t-2\}: \left\{ \begin{array}{r} \mathbb{E}\left[z_{s+1}^{(i)} \mid \mathcal{F}_{\leq s}\right] \leq f_s^{(i)}(z_s, q) \\ \mathbb{E}\left[|z_{s+1}^{(i)} - z_s^{(i)}|^2 \mid \mathcal{F}_{\leq s}\right] \leq h_s^{(i)}(z_s, q) \\ |z_{s+1}^{(i)} - z_s^{(i)}| \leq g_s^{(i)}(z_s) \end{array} \right\} \quad \text{(i.D.1)}$$

then it holds $\mathbf{Pr}_{x_1,\ldots,x_{t-1}}[\exists i \in [D'] : z_{t-1}^{(i)} > \phi_{t,t-1}^{(i)}] \leq q^2/2$.

Under the above two assumptions, we have for every $t \in [T]$, it satisfies $\mathbf{Pr}[\overline{\mathcal{C}_t}] \leq 2tq$.

*Proof of Lemma i.D.1.* We prove the lemma by induction. For the base case, by applying assumption (A2) we know that $\mathbf{Pr}_{x_1}\left[\exists i \in [D'] : y_{1,0}^{(i)} > \phi_{1,0}^{(i)}\right] \leq \mathbf{Pr}[\overline{\mathcal{E}_1}] \leq q^2/2 \leq q$ so event $\mathcal{C}_1$ holds with probability at least $1 - q$. In other words, $\mathbf{Pr}[\overline{\mathcal{C}_{\leq 1}}] = \mathbf{Pr}[\overline{\mathcal{C}_1}] \leq q < 2q$.

Suppose $\mathbf{Pr}[\overline{\mathcal{C}_{\leq t-1}}] \leq 2(t-1)q$ is true for some $t \geq 2$, we will prove $\mathbf{Pr}[\overline{\mathcal{C}_{\leq t}}] \leq 2tq$. Since it satisfies $\mathbf{Pr}[\overline{\mathcal{C}_{\leq t}}] \leq \mathbf{Pr}[\overline{\mathcal{C}_{\leq t-1}}] + \mathbf{Pr}[\overline{\mathcal{C}_t}]$, it suffices to prove that $\mathbf{Pr}[\overline{\mathcal{C}_t}] \leq 2q$.

Note also $\mathbf{Pr}[\overline{\mathcal{C}_t}] \leq \mathbf{Pr}[\overline{\mathcal{C}_t'}] + \mathbf{Pr}[\overline{\mathcal{C}_t''} \mid \mathcal{C}_t']$ but the second quantity $\mathbf{Pr}[\overline{\mathcal{C}_t''} \mid \mathcal{C}_t']$ is no more than $q$ according to our definition of $\mathcal{C}_t'$ and $\mathcal{C}_t''$. Therefore, in the rest of the proof, it suffices to show $\mathbf{Pr}[\overline{\mathcal{C}_t'}] \leq q$.

We use $y_{t,s}(x_t, x_{\leq s})$ to emphasize that $y_{t,s}$ is an $\mathcal{F}_t \times \mathcal{F}_{\leq s}$ measurable random vector. Let us now fix $x_t$ to be a vector satisfying $\mathcal{E}_t$. Define $\{z_s\}_{s=0}^{t-1}$ to be a random process where each $z_s \in \mathbb{R}^D$ is $\mathcal{F}_{\leq s}$ measurable:

$$z_s^{(i)} = z_s^{(i)}(x_{\leq s}) \stackrel{\text{def}}{=} \begin{cases} y_{t,s}^{(i)}(x_t, x_{\leq s}) & \text{if } x_{\leq s} \text{ satisfies } \mathcal{C}_{\leq s}; \\ \min\left\{f_{s-1}^{(i)}(z_{s-1}(x_{\leq s-1}), 0), z_{s-1}^{(i)}(x_{\leq s-1})\right\} & \text{if } x_{\leq s} \text{ satisfies } \overline{\mathcal{C}_{\leq s}}. \end{cases} \quad \text{(i.D.2)}$$

Then $z_s^{(i)}$ satisfies for every $i \in [D], s \leq \{0, 1, \ldots, t-2\}$,

$$\mathbb{E}\left[z_{s+1}^{(i)} \mid \mathcal{F}_{\leq s}\right] = \mathbf{Pr}[\mathcal{C}_{\leq s+1} \mid \mathcal{F}_{\leq s}] \cdot \mathbb{E}\left[z_{s+1}^{(i)} \mid \mathcal{C}_{\leq s+1}, \mathcal{F}_{\leq s}\right] + \mathbf{Pr}\left[\overline{\mathcal{C}_{\leq s+1}} \mid \mathcal{F}_{\leq s}\right] \cdot \mathbb{E}\left[z_{s+1}^{(i)} \mid \overline{\mathcal{C}_{\leq s+1}}, \mathcal{F}_{\leq s}\right]$$

$$\stackrel{①}{\leq} \mathbf{Pr}[\mathcal{C}_{\leq s+1} \mid \mathcal{F}_{\leq s}] \cdot \mathbb{E}\left[y_{t,s+1}^{(i)} \mid \mathcal{C}_{\leq s+1}, \mathcal{F}_{\leq s}\right] + \mathbf{Pr}\left[\overline{\mathcal{C}_{\leq s+1}} \mid \mathcal{F}_{\leq s}\right] \cdot f_s^{(i)}(z_s, 0)$$

$$\stackrel{②}{\leq} \mathbf{Pr}[\mathcal{C}_{\leq s+1} \mid \mathcal{F}_{\leq s}] \cdot f_s^{(i)}(y_{t,s}, q) + \mathbf{Pr}\left[\overline{\mathcal{C}_{\leq s+1}} \mid \mathcal{F}_{\leq s}\right] \cdot f_s^{(i)}(z_s, q)$$

$$\stackrel{③}{\leq} \mathbf{Pr}[\mathcal{C}_{\leq s+1} \mid \mathcal{F}_{\leq s}] \cdot f_s^{(i)}(y_{t,s}, q) + \mathbf{Pr}\left[\overline{\mathcal{C}_{\leq s+1}} \mid \mathcal{F}_{\leq s}\right] \cdot f_s^{(i)}(y_{t,s}, q) \quad \text{(i.D.3)}$$

$$= f_s^{(i)}(z_s, q) \quad \text{(i.D.4)}$$

Above, ① is because whenever $\mathcal{C}_{\leq s+1}$ holds it satisfies $z_{s+1}^{(i)} = y_{t,s+1}^{(i)}$, as well as whenever $\overline{\mathcal{C}_{\leq s+1}}$ holds it satisfies $z_{s+1}^{(i)} \leq f_s^{(i)}(z_s, 0)$; ② uses assumptions (A1a) and (A1d) as well as the fact that we have fixed $x_t$; ③ uses the fact that whenever $\mathbf{Pr}\left[\mathcal{C}_{\leq s+1} \mid \mathcal{F}_{\leq s}\right] > 0$ it must hold that $\mathcal{C}_{\leq s}$ is satisfied, and therefore it satisfies $y_{t,s} = z_s$.

Similarly, we can also show for every $i \in [D], s \leq \{0, 1, \ldots, t-2\}$,

$$\mathbb{E}\left[|z_{s+1}^{(i)} - z_s^{(i)}|^2 \mid \mathcal{F}_{\leq s}\right]$$

$$= \mathbf{Pr}[\mathcal{C}_{\leq s+1} \mid \mathcal{F}_{\leq s}] \cdot \mathbb{E}\left[|z_{s+1}^{(i)} - z_s^{(i)}|^2 \mid \mathcal{C}_{\leq s+1}, \mathcal{F}_{\leq s}\right] + \mathbf{Pr}\left[\overline{\mathcal{C}_{\leq s+1}} \mid \mathcal{F}_{\leq s}\right] \cdot \mathbb{E}\left[|z_{s+1}^{(i)} - z_s^{(i)}|^2 \mid \overline{\mathcal{C}_{\leq s+1}}, \mathcal{F}_{\leq s}\right]$$

$$\stackrel{①}{\leq} \mathbf{Pr}[\mathcal{C}_{\leq s+1} \mid \mathcal{F}_{\leq s}] \cdot \mathbb{E}\left[|y_{t,s+1}^{(i)} - y_{t,s}^{(i)}|^2 \mid \mathcal{C}_{\leq s+1}, \mathcal{F}_{\leq s}\right] + \mathbf{Pr}\left[\overline{\mathcal{C}_{\leq s+1}} \mid \mathcal{F}_{\leq s}\right] \cdot h_s^{(i)}(z_s, 0)$$

$$\stackrel{②}{\leq} \mathbf{Pr}[\mathcal{C}_{\leq s+1} \mid \mathcal{F}_{\leq s}] \cdot h_s^{(i)}(y_{t,s}, q) + \mathbf{Pr}\left[\overline{\mathcal{C}_{\leq s+1}} \mid \mathcal{F}_{\leq s}\right] \cdot h_s^{(i)}(z_s, q)$$

$$\stackrel{③}{\leq} h_s^{(i)}(z_s, q) . \quad \text{(i.D.5)}$$

Above, ① is because whenever $\mathcal{C}_{\leq s+1}$ holds it satisfies $z_{s+1}^{(i)} = y_{t,s+1}^{(i)}$ and $z_s^{(i)} = y_{t,s}^{(i)}$, together with whenever $\overline{\mathcal{C}_{\leq s+1}}$ holds it satisfies $|z_{s+1}^{(i)} - y_s^{(i)}|^2$ either equal zero or equal $|f_s^{(i)}(z_s, 0) - z_s^{(i)}|^2$, but in the latter case we must have $f_s^{(i)}(z_s, 0) < z_s^{(i)}$ (owing to (i.D.2)) and therefore it holds $|f_s^{(i)}(z_s, 0) - z_s^{(i)}|^2 \leq h_s^{(i)}(z_s, 0)$ using assumption (A1e). ② uses assumptions (A1b) and (A1d) as well as the fact that we have fixed $x_t$. ③ uses the fact that whenever $\mathbf{Pr}\left[\mathcal{C}_{\leq s+1} \mid \mathcal{F}_{\leq s}\right] > 0$ then $\mathcal{C}_{\leq s}$ must hold, and therefore it satisfies $y_{t,s} = z_s$.

Finally, we also have
$$|z_{s+1}^{(i)} - z_s^{(i)}| \leq g_s^{(i)}(z_s^{(i)}) \ . \tag{i.D.6}$$
This is so because whenever $\mathcal{C}_{\leq s+1}$ holds it satisfies $|z_{s+1}^{(i)} - z_s^{(i)}| = |y_{t,s+1}^{(i)} - y_{t,s}^{(i)}|$ so we can apply assumption (A1c). Otherwise, $\overline{\mathcal{C}_{\leq s+1}}$ holds we either have $|z_{s+1}^{(i)} - z_s^{(i)}| = 0$ (so (i.D.6) trivially holds) or $|z_{s+1}^{(i)} - z_s^{(i)}| = |f_s^{(i)}(z_s, 0) - z_s^{(i)}|$, but in the latter case we must have $f_s^{(i)}(z_s, 0) < z_s^{(i)}$ (owing to (i.D.2)) so it must satisfy $|f_s^{(i)}(z_s, 0) - z_s^{(i)}| \leq g_s^{(i)}(z_s)$ using assumption (A1e).

We are now ready to apply assumption (A3), which together with (i.D.4), (i.D.5), (i.D.6), implies that (recalling we have fixed $x_t$ to be any vector satisfying $\mathcal{E}_t$)
$$\mathop{\mathbf{Pr}}_{x_1,\ldots,x_{t-1}}\left[\exists i \in [D'] : z_{t-1}^{(i)} > \phi_{t,t-1}^{(i)} \mid \mathcal{E}_t\right] \leq q^2/2 \ .$$
This implies, after translating back to the random process $\{y_{t,s}\}$, we have
$$\mathop{\mathbf{Pr}}_{x_1,\ldots,x_t}\left[\exists i \in [D'] : y_{t,t-1}^{(i)} > \phi_{t,t-1}^{(i)}\right] \leq \mathop{\mathbf{Pr}}_{x_1,\ldots,x_t}\left[\exists i \in [D'] : y_{t,t-1}^{(i)} > \phi_{t,t-1}^{(i)} \mid \mathcal{E}_t\right] + \mathbf{Pr}[\overline{\mathcal{E}_t}]$$
$$\leq \mathop{\mathbf{Pr}}_{x_1,\ldots,x_{t-1}}\left[\exists i \in [D'] : z_{t-1}^{(i)} > \phi_{t,t-1}^{(i)} \mid \mathcal{E}_t\right] + q^2/2$$
$$\leq q^2/2 + q^2/2 = q^2 \ .$$
where the last inequality uses (A2). Finally, using Markov's inequality,
$$\mathop{\mathbf{Pr}}_{x_1,\ldots,x_{t-1}}\left[\overline{\mathcal{C}_t'}\right] = \mathop{\mathbf{Pr}}_{x_1,\ldots,x_{t-1}}\left[\mathop{\mathbf{Pr}}_{x_t}\left[\exists i \in [D'] : y_{t,t-1}^{(i)} > \phi_{t,t-1}^{(i)} \mid \mathcal{F}_{\leq t-1}\right] > q\right]$$
$$\leq \frac{1}{q} \cdot \mathop{\mathbb{E}}_{x_1,\ldots,x_{t-1}}\left[\mathop{\mathbf{Pr}}_{x_t}[\exists i \in [D'] : y_{t,t-1}^{(i)} > \phi_{t,t-1}^{(i)} \mid \mathcal{F}_{\leq t-1}]\right]$$
$$= \frac{1}{q} \cdot \mathop{\mathbf{Pr}}_{x_1,\ldots,x_t}\left[\exists i \in [D'] : y_{t,t-1}^{(i)} > \phi_{t,t-1}^{(i)}\right] \leq q \ .$$
Therefore, we finish proving $\mathbf{Pr}[\overline{\mathcal{C}_t'}] \leq q$ which implies $\mathbf{Pr}[\overline{\mathcal{C}_{\leq t}}] \leq 2tq$ as desired. This finishes the proof of Lemma i.D.1. $\square$

## i.E  Main Lemmas (for Section 7)

### i.E.1  Before Warm Start

*Proof of Lemma Main 1.* For every $t \in [T]$ and $s \in \{0, 1, \ldots, t-1\}$, consider random vectors $y_{t,s} \in \mathbb{R}^{T+2}$ defined as:
$$y_{t,s}^{(1)} \stackrel{\text{def}}{=} \|\mathbf{Z}^\top \mathbf{P}_s \mathbf{Q} (\mathbf{V}^\top \mathbf{P}_s \mathbf{Q})^{-1}\|_F^2 \ ,$$
$$y_{t,s}^{(2)} \stackrel{\text{def}}{=} \|\mathbf{W}^\top \mathbf{P}_s \mathbf{Q} (\mathbf{V}^\top \mathbf{P}_s \mathbf{Q})^{-1}\|_F^2 \ ,$$
$$y_{t,s}^{(3+j)} \stackrel{\text{def}}{=} \begin{cases} \left\|x_t^\top \mathbf{Z}\mathbf{Z}^\top (\mathbf{\Sigma}/\lambda_{k+1})^j \mathbf{P}_s \mathbf{Q} (\mathbf{V}^\top \mathbf{P}_s \mathbf{Q})^{-1}\right\|_2^2, & \text{for } j \in \{0, 1, \ldots, t-s-1\}; \\ (1 - \eta_s \lambda_k) \cdot y_{t,s-1}^{(3+j)}, & \text{for } j \in \{t-s, \ldots, T-1\}. \end{cases}$$

(In fact, we are only interested in $y_{t,s}^{(3+j)}$ for $j \leq t - s - 1$, and can "almost" define $y_{t,s}^{(3+j)} = +\infty$ whenever $j \geq t - s$. However, we still decide to give such out-of-boundary variables meaningful values in order to make all of our vectors $y_{t,s}$ (and functions $f, g, h$ defined later) to be of the same dimension $T + 2$. This allows us to greatly simplify our notations.)

We consider upper bounds

$$\phi_{t,s}^{(1)} \stackrel{\text{def}}{=} 2\Xi_{\mathbf{Z}}, \quad \phi_{t,s}^{(2)} \stackrel{\text{def}}{=} \begin{cases} 2\Xi_{\mathbf{Z}} & s < T_0; \\ 2 & \text{otherwise.} \end{cases}, \quad \text{and } \phi_{t,s}^{(3+j)} \stackrel{\text{def}}{=} 2\Xi_x^2 .$$

For each $t \in [T]$, define event $\mathcal{C}_t'$ and $\mathcal{C}_t''$ in the same way as decoupling Lemma i.D.1 (with $D' = 3$):

$$\mathcal{C}_t' \stackrel{\text{def}}{=} \left\{ (x_1, ..., x_{t-1}) \text{ satisfies } \Pr_{x_t}\left[ \exists i \in [3] : y_{t,t-1}^{(i)} > \phi_{t,t-1}^{(i)} \,\big|\, \mathcal{F}_{t-1} \right] \leq q \right\}$$

$$\mathcal{C}_t'' \stackrel{\text{def}}{=} \left\{ (x_1, ..., x_t) \text{ satisfies } \forall i \in [3] : y_{t,t-1}^{(i)} \leq \phi_{t,t-1}^{(i)} \right\}$$

and denote by $\mathcal{C}_t \stackrel{\text{def}}{=} \mathcal{C}_t' \wedge \mathcal{C}_t''$ and $\mathcal{C}_{\leq t} \stackrel{\text{def}}{=} \bigwedge_{s=1}^{t} \mathcal{C}_s$.

As a result, if $\mathcal{C}_{\leq s+1}$ holds, then we always have

$$\|x_{s+1}^\top \mathbf{P}_s \mathbf{Q} (\mathbf{V}^\top \mathbf{P}_s \mathbf{Q})^{-1}\|_2^2 \leq \left( \|x_{s+1}^\top \mathbf{V}\mathbf{V}^\top \mathbf{P}_s \mathbf{Q} (\mathbf{V}^\top \mathbf{P}_s \mathbf{Q})^{-1}\|_2 + \|x_{s+1}^\top \mathbf{Z}\mathbf{Z}^\top \mathbf{P}_s \mathbf{Q} (\mathbf{V}^\top \mathbf{P}_s \mathbf{Q})^{-1}\|_2 \right)^2$$

$$\leq (1 + \phi_{s+1,s}^{(3)})^2 = (\sqrt{2}\Xi_x + 1)^2 \leq 4\Xi_x^2 ,$$

where last inequality uses $\Xi_x \geq 2$. This allows us to later apply Lemma i.B.1 and Lemma 6.1 with $\phi_t = 2\Xi_x$.

**Verification of Assumption (A1) in Lemma i.D.1.**

Suppose $\mathbb{E}[x_s x_s^\top \mid \mathcal{C}_{\leq s}, \mathcal{F}_{\leq s-1}] = \mathbf{\Sigma} + \mathbf{\Delta}$, and we want to bound $\|\mathbf{\Delta}\|_2$. Defining $q_1 \stackrel{\text{def}}{=} \Pr[\overline{\mathcal{C}_s''} \mid \mathcal{C}_s', \mathcal{C}_{\leq s-1}, \mathcal{F}_{\leq s-1}]$, then we must have $q_1 \leq q$ according to the definition of $\mathcal{C}_s'$ and $\mathcal{C}_s''$. Using law of total expectation:

$$\mathbb{E}[x_s x_s^\top \mid \mathcal{C}_s', \mathcal{C}_{\leq s-1}, \mathcal{F}_{\leq s-1}] = \mathbb{E}[x_s x_s^\top \mid \mathcal{C}_{\leq s}, \mathcal{F}_{\leq s-1}] \cdot (1 - q_1) + \mathbb{E}[x_s x_s^\top \mid \overline{\mathcal{C}_s''}, \mathcal{C}_s', \mathcal{C}_{\leq s-1}, \mathcal{F}_{\leq s-1}] \cdot q_1 ,$$

and combining it with the fact that $0 \preceq x_s x_s^\top \preceq \mathbf{I}$ and $\mathbb{E}[x_s x_s^\top \mid \mathcal{C}_s', \mathcal{C}_{\leq s-1}, \mathcal{F}_{\leq s-1}] = \mathbb{E}[x_s x_s^\top] = \mathbf{\Sigma}$, we have[17]

$$\mathbf{\Sigma} \preceq (\mathbf{\Sigma} + \mathbf{\Delta})(1 - q_1) + q_1 \cdot \mathbf{I} \quad \text{and} \quad \mathbf{\Sigma} \succeq (\mathbf{\Sigma} + \mathbf{\Delta})(1 - q_1) .$$

After rearranging, these two properties imply $\|\mathbf{\Delta}\|_2 \leq \frac{q_1}{1-q_1} \leq \frac{q}{1-q}$.

Now, we can apply Lemma 6.1 and obtain for every $t \in [T]$, $s \in \{0, 1, \ldots, t-2\}$, and every $j \in \{0, 1, \ldots, T-1\}$, it satisfies[18]

$$\mathbb{E}[y_{t,s+1}^{(1)} \mid \mathcal{F}_t, \mathcal{F}_{\leq s}, \mathcal{C}_{\leq s+1}] \leq (1 + 56\eta_{s+1}^2 \Xi_x^2) y_{t,s}^{(1)} + 40\eta_{s+1}^2 \Xi_x^2 + 20\eta_{s+1} \frac{q\Xi_{\mathbf{Z}}^{3/2}}{1-q} ,$$

$$\mathbb{E}[y_{t,s+1}^{(2)} \mid \mathcal{F}_t, \mathcal{F}_{\leq s}, \mathcal{C}_{\leq s+1}] \leq (1 - 2\eta_{s+1}\rho + 56\eta_{s+1}^2 \Xi_x^2) y_{t,s}^{(2)} + 40\eta_{s+1}^2 \Xi_x^2 + 20\eta_{s+1} \frac{q\Xi_{\mathbf{Z}}^{3/2}}{1-q} , \text{ and}$$

$$\mathbb{E}[y_{t,s+1}^{(3+j)} \mid \mathcal{F}_t, \mathcal{F}_{\leq s}, \mathcal{C}_{\leq s+1}] \leq (1 - \eta_{s+1}\lambda_k + 56\eta_{s+1}^2 \Xi_x^2) y_{t,s}^{(3+j)} + \eta_{s+1}\lambda_k y_{t,s}^{(3+j+1)} + 40\eta_{s+1}^2 \Xi_x^2 + 20\eta_{s+1} \frac{q\Xi_{\mathbf{Z}}^{3/2}}{1-q}.$$

Moreover, for every $i \in [T+2]$, using Lemma i.B.1-(c) with $\phi_t = 2\Xi_x$ we have whenever $\mathcal{C}_{\leq s+1}$

---

[17]Here, we use notation $\mathbf{A} \preceq \mathbf{B}$ to indicate spectral dominance: that is, $\mathbf{B} - \mathbf{A}$ is positive semidefinite.

[18]We make a few comments regarding how to derive these upper bounds.
- Event $\mathcal{C}_{\leq s+1}$ implies we can safely apply Lemma 6.1 with $\phi_t = 2\Xi_x$.
- To obtain the third inequality, one needs to use the fact that when $w = x_t \mathbf{Z}\mathbf{Z}^\top$ the quantity $\frac{\eta_t}{\lambda_k}\|w^\top \mathbf{\Sigma} \mathbf{L}_{t-1}\|_2^2$ that appeared in Corollary 6.1-(c) can be upper bounded by $\frac{\eta_t \lambda_{k+1}^2}{\lambda_k}\|w^\top \mathbf{L}_{t-1}\|_2^2 \leq \eta_t \lambda_k \|w^\top \mathbf{L}_{t-1}\|_2^2$.
- Whenever $j \geq t - s$ we have $y_{t,s}^{(3+j)}$ is "out of boundary" which means $y_{t,s}^{(3+j)} \stackrel{\text{def}}{=} (1 - \eta_s \lambda_k) \cdot y_{t,s-1}^{(3+j)}$ according to our definition. In this case it easily satisfies the third inequality.
- When $j = T - 1$, in the third inequality we should replace $y_{t,s}^{(3+j+1)}$ with zero because it is out of bound.

holds it satisfies
$$\left|y_{t,s+1}^{(i)} - y_{t,s}^{(i)}\right| \leq 18\eta_{s+1}\Xi_x \cdot y_{t,s}^{(i)} + 4\eta_{s+1}\Xi_x \cdot \sqrt{y_{t,s}^{(i)}} + 40\eta_{s+1}^2\Xi_x^2 \leq 20\eta_{s+1}\Xi_x \cdot y_{t,s}^{(i)} + 42\eta_{s+1}^2\Xi_x^2 \ .$$

Putting the above bounds together, one can verify that the random process $\{y_{t,s}\}_{t\in[T], s\leq t-1}$ satisfy assumption (A1) of Lemma i.D.1 with[19]

$$f_s^{(1)}(y,q) = (1 + 56\eta_{s+1}^2\Xi_x^2)y^{(1)} + 40\eta_{s+1}^2\Xi_x^2 + 20\eta_{s+1}\frac{q\Xi_{\mathbf{Z}}^{3/2}}{1-q} \ ,$$

$$f_s^{(2)}(y,q) = (1 - 2\eta_{s+1}\rho + 56\eta_{s+1}^2\Xi_x^2)y^{(2)} + 40\eta_{s+1}^2\Xi_x^2 + 20\eta_{s+1}\frac{q\Xi_{\mathbf{Z}}^{3/2}}{1-q} \ ,$$

$$f_s^{(3+j)}(y,q) = (1 - \eta_{s+1}\lambda_k + 56\eta_{s+1}^2\Xi_x^2)y^{(3+j)} + \eta_{s+1}\lambda_k y^{(3+j+1)} + 40\eta_{s+1}^2\Xi_x^2 + 20\eta_{s+1}\frac{q\Xi_{\mathbf{Z}}^{3/2}}{1-q} \ ,$$

$$g_s^{(i)}(y) = 20\eta_{s+1}\Xi_x \cdot y^{(i)} + 42\eta_{s+1}^2\Xi_x^2 \ , \text{ and}$$

$$h_s^{(i)}(y,q) = \left(g_s^{(i)}(y)\right)^2$$

**Verification of Assumption (A2) of Lemma i.D.1.**

For coordinates $i = 1$ and $i = 2$, our assumption $\|\mathbf{Z}^\top \mathbf{Q}(\mathbf{V}^\top \mathbf{Q})^{-1}\|_F^2 \leq \Xi_{\mathbf{Z}}$ implies $y_{t,0}^{(i)} \leq \Xi_{\mathbf{Z}} < \phi_{t,0}^{(i)}$. For coordinates $i \geq 3$, we have assumption $\mathbf{Pr}_{x_t}\left[\forall j \in [T], \|x_t^\top \mathbf{Z}\mathbf{Z}^\top (\mathbf{\Sigma}/\lambda_{k+1})^{j-1}\mathbf{Q}(\mathbf{V}^\top \mathbf{Q})^{-1}\|_2 \leq \Xi_x\right] \geq 1 - q^2/2$. Together, event $\mathcal{E}_t$ (recall $\mathcal{E}_t \stackrel{\text{def}}{=} \{x_t \text{ satisfies } \forall i \in [D]: y_{t,0}^{(i)} \leq \phi_{t,0}^{(i)}\}$) holds for all $t \in [T]$ with probability at least $1 - q^2/2$. In sum, assumption (A2) is satisfied in Lemma i.D.1.

**Verification of Assumption (A3) of Lemma i.D.1.**

For every $t \in [T]$, at a high level assumption (A3) is satisfied once we plug in the following three sets of parameter choices to Corollary i.C.4 and Corollary i.C.5: for every $s \in [T-1]$, define

$$\begin{aligned}
\beta_{s,1} &= 0, & \delta_{s,1} &= 0, & \tau_{s,1} &= 20\eta_{s+1}\Xi_x \\
\beta_{s,2} &= 2\eta_{s+1}\rho, & \delta_{s,2} &= 0, & \tau_{s,2} &= 20\eta_{s+1}\Xi_x \\
\beta_{s,3} &= 0, & \delta_{s,3} &= \eta_{s+1}\lambda_k & \tau_{s,3} &= 20\eta_{s+1}\Xi_x
\end{aligned}$$

More specifically, for every $t \in [T]$, let $\{z_s\}_{s=0}^{t-1}$ be the *arbitrary* random vector satisfying (i.D.1) of Lemma i.D.1. Define $q_2 = q^2/8$.

- For coordinate $i = 1$ of $\{z_s\}_{s=0}^{t-1}$,
  - apply Corollary i.C.4 with $\{\beta_{s,1}, \delta_{s,1}, \tau_{s,1}\}_{s=0}^{t-2}$, $q = q_2$, $D = 1$, and $\kappa = 1$;

- For coordinate $i = 2$ of $\{z_s\}_{s=0}^{t-1}$,
  - if $t < T_0$, apply Corollary i.C.4 with $\{\beta_{s,2}, \delta_{s,2}, \tau_{s,2}\}_{s=0}^{t-2}$, $q = q_2$, $D = 1$, and $\kappa = 1$;
  - if $t \geq T_0$, apply Corollary i.C.5 with $\{\beta_{s,2}, \delta_{s,2}, \tau_{s,2}\}_{s=0}^{t-2}$, $q = q_2$, $D = 1$, $\gamma = 1$, and $\kappa = 1$;

- For coordinates $i = 3, 4, \ldots, T+2$ of $\{z_s\}_{s=0}^{t-1}$,
  - apply Corollary i.C.4 with $\{\beta_{s,3}, \delta_{s,3}, \tau_{s,3}\}_{s=0}^{t-2}$, $q = q_2$, $D = T$, and $\kappa = 1$.

---

[19] The only part of (A1) that is non-trivial to verify is (A1e) for $g_s^{(i)}$. Whenever $f_s^{(i)}(x,0) \leq x^{(i)}$, it satisfies

$$\left|f_s^{(i)}(x,0) - x^{(i)}\right| \leq \begin{cases} 0, & \text{if } i = 1; \\ 2\eta_{s+1}\rho \cdot x^{(2)}, & \text{if } i = 2; \\ \eta_{s+1}\lambda_k \cdot x^{(i)}, & \text{if } i \geq 3. \end{cases} \leq 2\eta_{s+1} \cdot x^{(i)} \leq g_s^{(i)}(x) \ ,$$

where the second inequality uses $\rho, \lambda_k \leq 1$ and the last inequality uses $\Xi_x \geq 2$.

One needs to verify that the assumptions of Corollary i.C.4 and i.C.5 are satisfied as follows. First of all, one can carefully check that our parameters $\beta, \delta, \tau$ satisfy (i.C.2) with $\kappa = 1$ and this needs our assumption $q \leq \frac{\eta_{s+1}}{\Xi_{\mathbf{Z}}^{3/2}}$. Next, we can apply Corollary i.C.4 because our assumptions on $\eta_s$ imply $\sum_{s=0}^{T-1} \tau_{s,i}^2 \leq \frac{1}{100} \ln^{-1} \frac{4T}{q_2}$ and $\tau_{s,i} \leq \frac{1}{24 \log(4T/q_2)}$ for $i = 1, 2, 3$. To verify the presumption of Corollary i.C.5 with $\gamma = 1$, we notice that

- our assumption $\eta_s \leq \frac{\rho}{4000 \cdot \Xi_x^2 \ln \frac{3T}{q_2}}$ implies $\beta_{s,2} \geq 10 \ln \frac{3T}{q_2} \cdot \tau_{s,2}^2$ and $\kappa \tau_s \leq \frac{1}{12 \ln \frac{3T}{q_2}}$ for every $s$,

- our assumption $\sum_{s=0}^{T_0-1} \beta_{s,2} \geq 1 + \ln \Xi_{\mathbf{Z}}$ implies $\sum_{s=0}^{t-1} \beta_s - 10 \ln \frac{3t}{q_2} \tau_s^2 \geq \ln \Xi_{\mathbf{Z}} + 1 - 1 = \ln \Xi_{\mathbf{Z}}$ whenever $t > T_0$,

Therefore, the conclusion of Corollary i.C.4 and Corollary i.C.5 imply that

$$\mathbf{Pr}[\exists i \in [3] : z_{t-1}^{(i)} > \phi_{t,t-1}^{(i)}] \leq 3q_2 < q^2/2$$

so assumption (A3) of Lemma i.D.1 holds.

**Application of Lemma i.D.1.** Applying Lemma i.D.1, we have $\mathbf{Pr}[\overline{\mathcal{C}_T}] \leq 2qT$ which implies our desired bounds and this finishes the proof of Lemma Main 1. $\square$

### i.E.2 After Warm Start

*Proof of Lemma Main 2.* For every $t \in [T]$ and $s \in \{0, 1, ..., t-1\}$, consider *the same* random vectors $y_{t,s} \in \mathbb{R}^{T+2}$ defined in the proof of Lemma Main 1:

$$y_{t,s}^{(1)} \stackrel{\text{def}}{=} \|\mathbf{Z}^\top \mathbf{P}_s \mathbf{Q} (\mathbf{V}^\top \mathbf{P}_s \mathbf{Q})^{-1}\|_F^2 \ ,$$

$$y_{t,s}^{(2)} \stackrel{\text{def}}{=} \|\mathbf{W}^\top \mathbf{P}_s \mathbf{Q} (\mathbf{V}^\top \mathbf{P}_s \mathbf{Q})^{-1}\|_F^2 \ ,$$

$$y_{t,s}^{(3+j)} \stackrel{\text{def}}{=} \begin{cases} \left\| x_t^\top \mathbf{Z} \mathbf{Z}^\top \left( \mathbf{\Sigma}/\lambda_{k+1} \right)^j \mathbf{P}_s \mathbf{Q} (\mathbf{V}^\top \mathbf{P}_s \mathbf{Q})^{-1} \right\|_2^2, & \text{for } j \in \{0, 1, \ldots, t-s-1\}; \\ (1 - \eta_s \lambda_k) \cdot y_{t,s-1}^{(3+j)}, & \text{for } j \in \{t-s, \ldots, T-1\}. \end{cases}$$

This time, we consider slightly different upper bounds

$$\phi_{t,s}^{(1)} \stackrel{\text{def}}{=} 2\Xi_{\mathbf{Z}}, \quad \phi_{t,s}^{(2)} \stackrel{\text{def}}{=} \begin{cases} 2\Xi_{\mathbf{Z}} & \text{if } s < T_0; \\ 2 & \text{if } s = T_0; \\ \frac{5T_0/\ln^2(T_0)}{s/\ln^2 s} & \text{if } s > T_0. \end{cases}, \quad \text{and } \phi_{t,s}^{(3+j)} \stackrel{\text{def}}{=} 2\Xi_x^2 \ .$$

We stress that the only difference between the above upper bounds and the ones we used in the proof of Lemma Main 1 is the choice of $\phi_{t,s}^{(2)}$ for $s > T_0$. Instead of setting it to be constant 2 for all such $s$, we make it decrease almost linearly with respect to index $s$.

Again, define event

$$\mathcal{C}'_t \stackrel{\text{def}}{=} \left\{ (x_1, ..., x_{t-1}) \text{ satisfies } \mathbf{Pr}_{x_t} \left[ \exists i \in [3] : y_{t,t-1}^{(i)} > \phi_{t,t-1}^{(i)} \,\middle|\, \mathcal{F}_{t-1} \right] \leq q \right\}$$

$$\mathcal{C}''_t \stackrel{\text{def}}{=} \left\{ (x_1, ..., x_t) \text{ satisfies } \forall i \in [3] : y_{t,t-1}^{(i)} \leq \phi_{t,t-1}^{(i)} \right\}$$

and denote by $\mathcal{C}_t \stackrel{\text{def}}{=} \mathcal{C}'_t \wedge \mathcal{C}''_t$ and $\mathcal{C}_{\leq t} \stackrel{\text{def}}{=} \bigwedge_{s=1}^t \mathcal{C}_s$.

We next want to apply the decoupling Lemma i.D.1.

**Verification of Assumption (A1) in Lemma i.D.1.**

The same functions $f_s^{(i)}$, $g_s^{(i)}$, and $h_s^{(i)}$ used in the proof of Lemma Main 1 still apply here. However, we want to make a minor change on $g_s^{(i)}$ whenever $s \geq T_0$.

Applying Lemma i.B.1-(c) with $\phi_t = 2\Xi_x$, we have whenever $\mathcal{C}_{\leq s+1}$ holds for some $s \geq T_0$ (which implies $y_{t,s}^{(2)} \leq 5$),

$$|y_{t,s+1}^{(2)} - y_{t,s}^{(2)}| \leq 18\eta_{s+1}\Xi_x y_{t,s}^{(2)} + 4\eta_{s+1}\Xi_x\sqrt{y_{t,s}^{(2)}} + 40\eta_{s+1}^2\Xi_x^2 \leq 45\eta_{s+1}\Xi_x\sqrt{y_{t,s}^{(2)}} + 40\eta_{s+1}^2\Xi_x^2 \ .$$

Therefore, we can choose

$$g_s^{(2)}(y) = 45\eta_{s+1}\Xi_x\sqrt{y^{(2)}} + 40\eta_{s+1}^2\Xi_x^2$$

for all $s \geq T_0$ and this still satisfies assumption (A1) of Lemma i.D.1.[20]

**Verification of Assumption (A2) of Lemma i.D.1.**

This is the same as the proof of Lemma Main 1.

**Verification of Assumption (A3) in Lemma i.D.1.**

Again, for every $t \in [T]$, let $\{z_s\}_{s=0}^{t-1}$ be the *arbitrary* random vector satisfying (i.D.1) of Lemma i.D.1. Choosing $q_2 = q^2/8$ again, the same proof of Lemma Main 1 shows that

$$\mathbf{Pr}[\exists i \in \{1,3\} : z_{t-1}^{(i)} > \phi_{t,t-1}^{(i)}] \leq 2q_2 \ .$$

Therefore, it suffices to prove that $\mathbf{Pr}[z_{t-1}^{(2)} > \phi_{t,t-1}^{(2)}] \leq 2q_2$.

We only need to focus on the case $t \geq T_0 + 2$, because otherwise if $t \leq T_0 + 1$ then $g_s^{(2)}$ is not changed for all $s \in \{0, \ldots, t-2\}$ so the same proof of Lemma Main 1 also shows $\mathbf{Pr}[z_{t-1}^{(2)} > \phi_{t,t-1}^{(2)}] \leq q_2$.

When $t \geq T_0 + 2$, we can first apply the same proof of Lemma Main 1 (for $t = T_0 + 1$) to show that $\mathbf{Pr}[z_{T_0}^{(2)} > \phi_{T_0+1,T_0}^{(2)} = 2] \leq q_2$. Next, conditioning on $z_{T_0}^{(2)} \leq 2$ which happens with probability at least $1 - q_2$, we want to apply Corollary i.C.3 with $\kappa = 2$ and $\tau_s = \frac{1}{\delta s}$.

More specifically, for every $t \in \{T_0 + 2, \ldots, T\}$, we have shown that the random sequence $\{z_s^{(2)}\}_{s=T_0}^{t-1}$ satisfies (i.D.1) with

$$f_s^{(2)}(y, q) = (1 - 2\eta_{s+1}\rho + 56\eta_{s+1}^2\Xi_x^2)y^{(2)} + 40\eta_{s+1}^2\Xi_x^2 + 20\eta_{s+1}\frac{q\Xi_\mathbf{Z}^{3/2}}{1-q}$$

$$g_s^{(2)}(y) = 45\eta_{s+1}\Xi_x\sqrt{y^{(2)}} + 40\eta_{s+1}^2\Xi_x^2$$

$$h_s^{(2)}(y, q) = \left(g_s^{(2)}(y)\right)^2$$

Therefore, $\{z_s^{(2)}\}_{s=T_0}^{t-1}$ also satisfies (i.C.1) with $\kappa = 2$ and $\tau_s = \frac{1}{\delta \tau_s}$ because the following holds from our assumptions:

$$q\Xi_\mathbf{Z}^{3/2} \leq \eta_{s+1} \qquad\qquad \delta\tau_s = \frac{1}{s} \leq 2\eta_{s+1}\rho - 56\eta_{s+1}^2\Xi_x^2$$

$$\tau_s^2 = \frac{1}{\delta^2 s^2} \geq 60\eta_{s+1}^2\Xi_x^2 \geq 40\eta_{s+1}^2\Xi_x^2 + 20\eta_{s+1}\frac{q\Xi_\mathbf{Z}^{3/2}}{1-q} \qquad \kappa\tau_s = \frac{2}{\delta s} \geq 40\eta_{s+1}\Xi_x$$

Now, we are ready to apply Corollary i.C.3 with $q = q_2$, $t_0 = T_0$, and $\kappa = 2$. Because $q_2 \leq e^{-2}$, $z_{T_0}^{(2)} \leq 2$, $\delta \leq 1/\sqrt{8}$ and $\frac{T_0}{\ln^2 T_0} \geq \frac{9\ln(1/q_2)}{\delta^2}$, the conclusion of Corollary i.C.3 tells us

$$\mathbf{Pr}[z_{t-1}^{(2)} > \phi_{t,t-1}^{(2)} \mid z_{T_0}^{(2)} \leq 2] \leq q_2 \ .$$

By union bound, we have $\mathbf{Pr}[z_{t-1}^{(2)} > \phi_{t,t-1}^{(2)}] \leq q_2 + q_2 = 2q_2$ as desired.

Finally, we conclude (for every $t \geq T_0 + 2$) that

$$\mathbf{Pr}[\exists i \in [3] : z_{t-1}^{(i)} > \phi_{t,t-1}^{(i)}] \leq 4q_2 < q^2/2$$

---

[20]Similar to Footnote 19, we also need to verify (A1e) for $g_s^{(i)}$. Whenever $f_s^{(2)}(x, 0) \leq x^{(2)}$, it satisfies $|f_s^{(2)}(x, 0) - x^{(2)}| \leq 2\eta_{s+1} \cdot x^{(2)} \leq g_s^{(2)}(x)$, where the first inequality uses $\rho \leq 1$ and the second uses $\Xi_x \geq 2$.

so assumption (A3) of Lemma i.D.1 holds.

**Application of Lemma i.D.1.** Applying Lemma i.D.1, we have $\mathbf{Pr}[\overline{\mathcal{C}_T}] \leq 2qT$ which implies our desired bounds and this finishes the proof of Lemma Main 2. $\square$

## i.F  Proof of Theorem 1' and 2'

*Proof of Theorem 2'.* First for a sufficiently large constant $C$, we can apply Lemma 5.1 with $p' = \frac{p}{6}$ and $q = \min\left\{\frac{1}{CT^2d^2}, \frac{p}{4T}\right\}$ and obtain: with probability at least $1 - p' - q^2 \geq 1 - p/2$ over the random choice of $\mathbf{Q}$, the following holds:

$$\begin{cases} \left\|(\mathbf{Z}^\top\mathbf{Q})(\mathbf{V}^\top\mathbf{Q})^{-1}\right\|_F^2 \leq \frac{20736dk}{p^2}\ln\frac{6d}{p} \text{ , and} \\ \mathbf{Pr}_{x_1,\ldots,x_T}\left[\exists i \in [T], \exists t \in [T], \left\|x_t^\top \mathbf{ZZ}^\top (\mathbf{\Sigma}/\lambda_{k+1})^{i-1}\mathbf{Q}(\mathbf{V}^\top\mathbf{Q})^{-1}\right\|_2 \geq \frac{216\sqrt{k\ln\frac{2T}{q}}}{p}\right] \leq \frac{q^2}{2} \text{ .} \end{cases}$$

Denote by $\mathcal{C}_1$ the union of the above two events, and we have $\mathbf{Pr}_{\mathbf{Q}}[\mathcal{C}_1] \geq 1 - p/2$.

Now, for every fixed $\mathbf{Q}$, whenever $\mathcal{C}_1$ holds, we can let

$$\Xi_{\mathbf{Z}} = \frac{20736dk}{p^2}\ln\frac{6d}{p}, \quad \Xi_x = \frac{216\sqrt{2k\ln\frac{2T}{p}}}{p} \text{ ,}$$

so the initial conditions in Lemma Main 1 (and thus Lemma Main 2) is satisfied. Also, according to Parameter 7.1, our parameter choices satisfy the assumptions in Lemma Main 2. Finally, the conclusion of Lemma Main 2 immediately implies for every $T \geq T_0$

$$\mathbf{Pr}_{x_1,\ldots,x_T}\left[\forall t = \{T_0, \ldots, T\} : \|\mathbf{W}^\top\mathbf{P}_t\mathbf{Q}(\mathbf{V}^\top\mathbf{P}_t\mathbf{Q})^{-1}\|_F^2 \leq \widetilde{O}\left(\frac{T_0}{t}\right) \bigg| \mathcal{C}_1\right] \geq 1 - 2qT \geq 1 - \frac{p}{2} \text{ .}$$

Union bounding this with event $\overline{\mathcal{C}_1}$, we have

$$\mathbf{Pr}_{\mathbf{Q},x_1,\ldots,x_T}\left[\forall t = \{T_0, \ldots, T\} : \|\mathbf{W}^\top\mathbf{P}_t\mathbf{Q}(\mathbf{V}^\top\mathbf{P}_t\mathbf{Q})^{-1}\|_F^2 \leq \widetilde{O}\left(\frac{T_0}{t}\right)\right] \geq 1 - p \text{ .}$$

Combining this with Lemma 2.2 completes the proof. $\square$

Finally, Theorem 1' is a direct corollary of Theorem 2' by setting $\rho \leftarrow \mathsf{gap}$.

# APPENDIX (PART II)

In this part II of the appendix, we provide our proofs for the lower bound Theorem 6, as well as that for the Rayleigh quotient Theorem 3.

- Appendix ii.G extends our main lemmas to better serve for the rayleigh quotient setting;
- Appendix ii.H provides the final proof for our Rayleigh Quotient Theorem 3;
- Appendix ii.I includes a three-paged proof of our lower bound Theorem 6.

Below we address, at a high level, the main ideas needed behind our main lemma extension as well as the proof of Theorem 3.

**Additional Ideas Needed for Theorem 3.** In order to prove Theorem 3 which is the rayleigh-quotient guarantee in gap-free streaming PCA, we want to strengthen Lemma Main 1 so that it provides guarantee essentially of the form:

$$\text{for every } \gamma \geq 1: \quad \left\|\mathbf{W}_\gamma^\top \mathbf{P}_t \mathbf{Q}(\mathbf{V}^\top \mathbf{P}_t \mathbf{Q})^{-1}\right\|_F^2 \leq 2/\gamma \ , \tag{ii.F.1}$$

where $\mathbf{W}_\gamma$ is the column orthonormal matrix consisting of all eigenvectors of $\mathbf{\Sigma}$ with eigenvalues $\leq \lambda_k - \gamma \cdot \rho$. For obvious reason Lemma Main 1 is a special case of (ii.F.1) when restricting only to $\gamma = 1$. It is a simple exercise to show that (ii.F.1) implies our desired rayleigh-quotient guarantee (via an Abel transformation and an integral computation, see Appendix ii.H).

Therefore, it suffices to prove (ii.F.1). If one were allowed to magically change learning rates and apply Lemma Main 1 multiple times, then (ii.F.1) would be trivial to prove: just replace $\mathbf{W}$ with $\mathbf{W}_\gamma$ and replacing $\rho$ with $\gamma \cdot \rho$ and repeatedly apply Lemma Main 1. Unfortunately, the difficulty arises because want to prove (ii.F.1) for all $\gamma \geq 1$ but with a *fixed* set of learning rates $\eta_t$.

In Appendix ii.G, we showed that the same learning rates in Parameter 7.1, together with a more general martingale concentration lemma (i.e., Corollary i.C.5 with $\gamma \geq 1$), one can obtain (ii.F.1) and we call it Lemma Main 3. This proof follows from the same structure as that of Lemma Main 1 except for the change in how we apply Corollary i.C.5.

Finally, Theorem 3 follows from Lemma Main 3 for a one-paged reason, see Appendix ii.H.

## ii.G  Improved Main Lemma

In this section we also sketch the proof to obtain Rayleigh quotient result. We will prove the following lemma which is a strengthened version of Lemma Main 1.

**Lemma Main 3** (before warm start). *In the same setting as Lemma Main 1, suppose we redefine $\mathbf{W} = \mathbf{W}_\gamma$ to be the column orthonormal matrix consisting of eigenvectors of $\mathbf{\Sigma}$ with values $\leq \lambda_k - \gamma \cdot \rho$.*

*Then, for every $\gamma \in [1, 1/\rho]$, with probability at least $1 - 2qT$:*

$$\forall t \in \{T_0, \ldots, T\}, \left\|\mathbf{W}_\gamma^\top \mathbf{P}_t \mathbf{Q} (\mathbf{V}^\top \mathbf{P}_t \mathbf{Q})^{-1}\right\|_F^2 \leq \frac{2}{\gamma} \ .$$

*Proof of Lemma Main 3.* The proof is a non-trivial adaption of the proof of Lemma Main 1.

We redefine $\mathbf{W} = \mathbf{W}_\gamma$ and consider random vectors $y_{t,s} \in \mathbb{R}^{T+2}$ defined in the same way as the proof of Lemma Main 1. This time, we consider upper bounds

$$\phi_{t,s}^{(1)} \stackrel{\text{def}}{=} 2\Xi_\mathbf{Z}, \quad \phi_{t,s}^{(2)} \stackrel{\text{def}}{=} \begin{cases} 2\Xi_\mathbf{Z} & s < T_0; \\ 2/\gamma & s \geq T_0. \end{cases}, \quad \text{and } \phi_{t,s}^{(3+j)} \stackrel{\text{def}}{=} 2\Xi_x^2$$

so the only difference we make here is on coordinate $i = 2$ for $s \geq T_0$. For each $t \in [T]$, we also consider events $\mathcal{C}'_t$, $\mathcal{C}''_t$, $\mathcal{C}_t \stackrel{\text{def}}{=} \mathcal{C}'_t \wedge \mathcal{C}''_t$, and $\mathcal{C}_{\leq t} \stackrel{\text{def}}{=} \bigwedge_{s=1}^t \mathcal{C}_s$ defined in the same way as before.

**Verification of Assumption (A1) in Lemma i.D.1.**

We consider the same functions $f_s$, $g_s$, $h_s$ as defined in the proof of Lemma Main 1, except that we replace $\rho$ with $\gamma \cdot \rho$ because this time we have redefined $\mathbf{W} = \mathbf{W}_\gamma$ so that it consists of eigenvectors with values $\leq \lambda_k - \gamma \cdot \rho$. In other words, we redefine

$$f_s^{(2)}(y, q) = (1 - 2\eta_{s+1}\gamma\rho + 56\Lambda\eta_{s+1}^2\Xi_x^2)y^{(2)} + 40\eta_{s+1}^2\Lambda\Xi_x^2 + 20\eta_{s+1}\frac{q\Xi_\mathbf{Z}^{3/2}}{1-q} \ .$$

In the same way we can verify that these functions satisfy assumption (A1) of Lemma i.D.1.

**Verification of Assumption (A2) of Lemma i.D.1.**

This step is exactly the same as the proof of Lemma Main 1 so ignored here.

**Verification of Assumption (A3) of Lemma i.D.1.**

We consider the same parameters $\{\beta_s, \delta_s, \tau_s\}_s$ as Lemma Main 1 except that at coordinate $i = 2$ we replace $\rho$ with $\gamma \cdot \rho$:

$$\beta_{s,2} = 2\eta_{s+1}\gamma\rho, \qquad \delta_{s,2} = 0, \qquad \tau_{s,2} = 20\Lambda\eta_{s+1}\Xi_x \ .$$

Now, for every $t \in [T]$, let $\{z_s\}_{s=0}^{t-1}$ be the *arbitrary* random vector satisfying (i.D.1) of Lemma i.D.1. Letting $q_2 = q^2/8$, we can handle coordinates $i = 1$ and $i \geq 3$ in the same way as before. As for coordinate $i = 2$ of $\{z_s\}_{s=0}^{t-1}$,

- if $t < T_0$, apply Corollary i.C.4 with $\{\beta_{s,2}, \delta_{s,2}, \tau_{s,2}\}_{s=0}^{t-2}$, $q = q_2$, $D = 1$, and $\kappa = \frac{1}{\sqrt{\Lambda}}$;

- if $t \geq T_0$, apply Corollary i.C.5 with $\{\beta_{s,2}, \delta_{s,2}, \tau_{s,2}\}_{s=0}^{t-2}$, $q = q_2$, $D = 1$, $\gamma = \gamma$, and $\kappa = \frac{1}{\sqrt{\Lambda}}$;

Note that the $t < T_0$ case is exactly the same as before. When $t \geq T_0$, we again apply Corollary i.C.5 but this time with value $\gamma \geq 1$ rather than $\gamma = 1$. Since this is the only difference here, we only need to verify the the presumptions of Corollary i.C.5:

- our assumption $\eta_s \leq \frac{\rho}{4000\Lambda \cdot \Xi_x^2 \ln \frac{3T}{q_2}}$ implies $\beta_{s,2} \geq 20\gamma \ln \frac{3T}{q_2} \cdot \tau_{s,2}^2$ and $\kappa\tau_s \leq \frac{1}{12 \ln \frac{3T}{q_2}}$ for every $s$,

- our assumption $\sum_{s=0}^{T_0-1} \beta_{s,2} \geq 1 + \ln \Xi_\mathbf{Z}$ implies $\sum_{s=0}^{t-1} \beta_{s,2} - 10\gamma \ln \frac{3t}{q_2} \tau_{s,2}^2 \geq \frac{1}{2}\sum_{s=0}^{t-1} \beta_{s,2} \geq \ln \Xi_\mathbf{Z}$ whenever $t > T_0$.

Therefore, in the same way as the old proof in Lemma Main 1, we can conclude using Corollary i.C.4 and Corollary i.C.5 that
$$\mathbf{Pr}[\exists i \in [3] : z_{t-1}^{(i)} > \phi_{t,t-1}^{(i)}] \leq 3q_2 < q^2/2 \ .$$
This verifies assumption (A3) of Lemma i.D.1.

**Application of Lemma i.D.1.** Applying Lemma i.D.1, we have $\mathbf{Pr}[\overline{\mathcal{C}_T}] \leq 2qT$ which implies our desired bounds and this finishes the proof of Lemma Main 3. □

## ii.H Proof of Theorem 3: Rayleigh Quotient

**Theorem 3** (restated). *In the same setting as Theorem 2', we have for every $T = \widetilde{\Theta}(\frac{k}{\rho^2 \cdot p^2})$, letting $q_i$ be the $i$-th column of the output matrix $\mathbf{Q}_T$, then*
$$\mathbf{Pr}\left[\forall i \in [k], \quad q_i^\top \mathbf{\Sigma} q_i \geq \lambda_i - \widetilde{\Theta}(\rho)\right] \geq 1 - p \ .$$
*Again, $\widetilde{\Theta}$ hides poly-log factors in $\frac{1}{p}, \frac{1}{\rho}$ and $d$.*

*Proof of Theorem 3.* Since the statement of Theorem 3 uses the same learning rates[21] as in Theorem 2', the same proof of Theorem 2' ensures that the initialization assumptions in Lemma Main 3 are satisfied and thus we can apply Lemma Main 3.

We want to prove next the output matrix $\mathbf{Q}_T = [q_1, \ldots, q_k] \in \mathbb{R}^{d \times k}$ satisfies

$$\text{with probability at least } 1 - (2kdT)q, \quad \forall i \in [k]\colon \quad q_i^\top \mathbf{\Sigma} q_i \geq \lambda_i - 3\rho \ln \frac{1}{\rho} \ . \tag{ii.H.1}$$

For every $i \in [k]$, let $\mathbf{Q}_T^i \in \mathbb{R}^{d \times i}$ denote the first $i$-columns of $\mathbf{Q}_T$. By the property of Oja's algorithm, the same $\mathbf{Q}_T^i$ would have been the output if we started from an $\mathbb{R}^{d \times i}$ random matrix $\mathbf{Q}_0$ for streaming $i$-PCA. In other words, we can write $\mathbf{Q}_T^i = [q_1, \ldots, q_i]$.

Letting $\mathbf{W}_\gamma^i$ be the column orthonormal matrix consisting of all eigenvectors of $\mathbf{\Sigma}$ with eigenvalue $\leq \lambda_i - \gamma \cdot \rho$, we applying Lemma Main 3 (with $k = i$) and obtain:

$$\text{w.p. at least } 1 - 2qT\colon \quad \|(\mathbf{W}_\gamma^i)^\top \mathbf{Q}_T^i\|_F^2 \leq \|(\mathbf{W}_\gamma^i)^\top \mathbf{P}_T \mathbf{Q}^i (\mathbf{V}^\top \mathbf{P}_T \mathbf{Q}^i)^{-1}\|_F^2 \leq 2/\gamma \ .$$

(Above, the first inequality uses Lemma 2.2.) This in particular implies $\|(\mathbf{W}_\gamma^i)^\top q_i\|_2^2 \leq \frac{2}{\gamma}$.

Let us define for each $i \in [k]$,
$$\Gamma_i \stackrel{\text{def}}{=} \left\{\frac{\lambda_i - \lambda_j}{\rho} \,\middle|\, \lambda_i - \lambda_j \geq \rho\right\} \subseteq \mathbb{R}_{\geq 1} \quad \text{and} \quad \gamma_{i,j} \stackrel{\text{def}}{=} \frac{\lambda_i - \lambda_j}{\rho} \in \left[1, \frac{1}{\rho}\right] \ .$$

By union bound,
$$\text{w.p. at least } 1 - 2qkdT, \quad \forall i \in [k], \forall \gamma \in \Gamma_i\colon \quad \|(\mathbf{W}_\gamma^i)^\top q_i\|_2^2 \leq 2/\gamma \ . \tag{ii.H.2}$$

We are now ready to bound Rayleigh quotient. For each $i \in [k]$, let $i_0$ be the index of the first (i.e., the largest) eigenvector with eigenvalue $\leq \lambda_i - \rho$ and define $b_{i,j} \stackrel{\text{def}}{=} \sum_{s=j}^d \langle q_i, \nu_j \rangle^2$ where $\nu_j$ is the $j$-th largest eigenvector of $\mathbf{\Sigma}$. It satisfies $b_{i,1} = 1$. By Abel's formula,

$$q_i^\top \mathbf{\Sigma} q_i = \sum_{j=1}^d \lambda_j \langle q_i, \nu_j \rangle^2 \geq (\lambda_i - \rho) - \sum_{j=i_0+1}^d b_{i,j}(\lambda_{j-1} - \lambda_j) \ .$$

---
[21] More precisely, we can use the same Parameter 7.1 together with the same values of $\Xi_x$ and $\Xi_\mathbf{Z}$ used in Section i.F.

Note that for every $j \geq i_0 + 1$, we have $b_{i,j} \leq \|\mathbf{W}^i_{\gamma_{i,j}} q_i\|_2^2 \leq \frac{2}{\gamma_{i,j}}$ according to (ii.H.2). Therefore,

$$\sum_{j=i_0+1}^{d} b_{i,j}(\lambda_{j-1} - \lambda_j) \leq \sum_{j=i_0+1}^{d} \frac{2}{\gamma_{i,j}} \rho(\gamma_{i,j} - \gamma_{i,j-1}) \leq 2\rho \int_{1}^{\frac{1}{\rho}} \frac{1}{\gamma} dz \leq 2\rho \ln \frac{1}{\rho} ,$$

which implies $q_i^\top \mathbf{\Sigma} q_i \geq \lambda_i - 3\rho \ln \frac{1}{\rho}$ so (ii.H.1) holds. $\square$

## ii.I Proof of Theorem 6: Lower Bound

In this section we prove the following lower bound. Without loss of generality, we only focus on the case when $d = 2(k+m)$ and the case for $d > 2(k+m)$ can be done by padding zeros. Denoting by $\mathcal{S}_{(2(k+m))\times k}$ the set of all (column) orthonormal matrix in $\mathbb{R}^{(2(k+m))\times k}$,

> **Theorem 6** (lower bound, restated)**.** *For every $k \in \mathbb{N}^*$, every $m \geq 0$ that is an integral multiple of $k$, for every $\lambda \in (0, \frac{1}{4(k+m)}]$, every $\delta \in (0, \frac{\lambda}{2}]$, every $T \in \mathbb{N}^*$ satisfying $T \geq \Omega(\frac{\lambda}{\delta^2})$, every algorithm $\mathcal{A} : (\mathbb{R}^{2(k+m)})^{\otimes T} \to \mathcal{S}_{(2(k+m))\times k}$, there exists a distribution $\mu$ over vectors in $\mathbb{R}^{2(k+m)}$ such that*
>
> *all $x \sim \mu$ satisfies $\|x\|_2 \leq 1$, $\lambda_k \left( \mathbb{E}_{x \sim \mu} [xx^\top] \right) \geq \lambda, \quad \lambda_{k+m+1} \left( \mathbb{E}_{x \sim \mu} [xx^\top] \right) \leq \lambda - \delta$ .*
>
> *Furthermore, let $\mathbf{Q}_T \overset{\text{def}}{=} \mathcal{A}(x_1, \ldots, x_T)$ be the output with respect to $T$ i.i.d. random inputs $x_1, \ldots, x_T$ from $\mu$, we have*
>
> $$\mathbb{E}_{x_1,\ldots,x_T,\mathcal{A}} \left[ \|\mathbf{Q}_T^\top \mathbf{W}\|_F^2 \right] = \Omega \left( \frac{k\lambda}{\delta^2 T} \right) ,$$
>
> *where $\mathbf{W}$ consists of all the last $d - (k+m)$ eigenvectors of $\mathbb{E}_{x \sim \mu}[xx^\top]$.*

We shall just prove this theorem for gap case, i.e. when $m = 0$. For the gap free case, the theorem follows by first obtaining $\mu$ for the $m = 0$ case and for $\lambda' = \lambda \cdot \frac{m+k}{k}$, and then modifying it into some $\mu'$ as follows: first draw $x^{(0)}, \ldots, x^{(m/k)}$ independent copies from $\mu$, and then concatenate them into $x' = (x^{(0)}, \ldots, x^{(m/k)}) \cdot \frac{\sqrt{k}}{\sqrt{m+k}}$ and output this $x'$ instead.

Throughout this section, we use the phrase $i$-th eigenvalue to denote the $i$-th largest eigenvalue of a matrix (tie breaking arbitrarily) and the $i$-th eigenvector to denote that corresponding to the $i$-th eigenvector (with unit norm). We first state a lemma regarding $2 \times 2$ matrices:

**Lemma ii.I.1** ($2 \times 2$ matrix)**.** *For every $\beta \in \left( \frac{\sqrt{2}}{2}, \frac{\sqrt{3}}{2} \right)$ and $\varepsilon \in (0, 2\beta^2 - 1]$, choose $a = \left[ \beta, \sqrt{1 - \beta^2} \right]^\top, b = \left[ \beta, -\sqrt{1 - \beta^2} \right]^\top \in \mathbb{R}^2$ and define*

$$\mathbf{A} = \frac{1}{2} aa^\top + \frac{1}{2} bb^\top, \quad \mathbf{B} = \left( \frac{1}{2} + \varepsilon \right) aa^\top + \left( \frac{1}{2} - \varepsilon \right) bb^\top .$$

*Then, letting $\nu_1^a$ and $\nu_1^b$ be the top eigenvectors of $\mathbf{A}$ and $\mathbf{B}$ respectively, letting $\lambda_1, \lambda_2$ be the eigenvalues of $\mathbf{A}$, and $\lambda_1^b, \lambda_2^b$ be the eigenvalues of $\mathbf{B}$, we have*

$$|\langle \nu_1^a, \nu_1^b \rangle|^2 \leq 1 - \frac{\varepsilon^2}{16(\lambda_1 - \lambda_2)^2} \quad \text{and} \quad \lambda_1^b > \lambda_1 = \beta^2 > 1 - \beta^2 = \lambda_2 > \lambda_2^b .$$

*Proof of Lemma ii.I.1.* We can calculate that $\mathbf{A} = \begin{pmatrix} \beta^2 & 0 \\ 0 & 1 - \beta^2 \end{pmatrix}$ and therefore $\lambda_1 = \beta^2$ and $\lambda_2 = 1 - \beta^2$. We can also compute $\mathbf{B} = \begin{pmatrix} \beta^2 & 2\varepsilon\beta\sqrt{1-\beta^2} \\ 2\varepsilon\beta\sqrt{1-\beta^2} & 1 - \beta^2 \end{pmatrix}$ and the eigenvectors of $\mathbf{B}$ are $\nu_i^b = \frac{1}{\sqrt{s_i^2+1}}(1, s_i)^\top$ with eigenvalues $\lambda_i^b = \beta^2 + 2\varepsilon\beta\sqrt{1-\beta^2}s_i$ for $i = 1, 2$. Here, $s_1 \geq s_2$ are the two roots of equation $s^2 + \left(\frac{2\beta^2-1}{2\varepsilon\beta\sqrt{1-\beta^2}}\right)s - 1 = 0$. Since exactly one of the roots is positive (denoted by $s_1$), we immediately have that $\lambda_1^b = \beta^2 + 2\varepsilon\beta\sqrt{1-\beta^2}s_1 > \lambda_1$. This also implies $\lambda_2^b < \lambda_2$ since $\lambda_1^b + \lambda_2^b = \mathbf{Tr}(\mathbf{B}) = 1 = \lambda_1 + \lambda_2$.

We now turn to the eigenvectors. Denote by $\alpha = \frac{2\beta^2-1}{2\varepsilon\beta\sqrt{1-\beta^2}}$ and it satisfies $\alpha \leq \frac{2(\lambda_1-\lambda_2)}{\varepsilon}$ because $\beta \in \left(\frac{\sqrt{2}}{2}, \frac{\sqrt{3}}{2}\right)$. Therefore,

$$s_1 = \frac{-\alpha + \sqrt{\alpha^2+4}}{2} = \frac{2}{\alpha + \sqrt{\alpha^2+4}} \geq \frac{1}{\sqrt{\alpha^2+4}} \geq \frac{1}{\sqrt{\frac{4(\lambda_1-\lambda_2)^2}{\varepsilon^2}+4}} \geq \frac{\varepsilon}{\sqrt{8}(\lambda_1-\lambda_2)} \ .$$

Above, the last inequality uses our assumption that $\varepsilon \leq 2\beta^2 - 1 = \lambda_1 - \lambda_2$. On the other hand, $s_1 = \frac{2}{\alpha+\sqrt{\alpha^2+4}} \leq 1$. Therefore, we have:

$$|\langle \nu_1^a, \nu_1^b \rangle|^2 = \frac{1}{1+s_1^2} \leq 1 - \frac{s_1^2}{2} \leq 1 - \frac{\varepsilon^2}{16(\lambda_1-\lambda_2)^2} \ . \qquad \square$$

**Corollary ii.I.2** (Corollary of Lemma ii.I.1)**.** *For every $\beta \in \left(\frac{\sqrt{2}}{2}, \frac{\sqrt{3}}{2}\right)$ and every $\varepsilon \in (0, \mathsf{gap}]$ where $\mathsf{gap} \stackrel{\text{def}}{=} 2\beta^2 - 1$, choose $a$ and $b$ in $\mathbb{R}^2$ as defined as in Lemma ii.I.1. Define distributions over $\{a, b\}$:*

$$\mu_1 \colon \mathbf{Pr}[x=a] = \mathbf{Pr}[x=b] = \frac{1}{2} \quad \text{and} \quad \mu_2 \colon \mathbf{Pr}[x=a] = \frac{1}{2}+\varepsilon, \mathbf{Pr}[x=b] = \frac{1}{2}-\varepsilon \ .$$

*Suppose we are given an arbitrary (possibly randomized) algorithm $\mathcal{A} : (\mathbb{R}^2)^{\otimes T} \to \mathbb{R}^2$ for $T < O(\frac{1}{\varepsilon^2})$. Then, if we pick $i \in \{1, 2\}$ each with probability $1/2$ and sample $T$ i.i.d. vectors $x_1, \ldots, x_T$ from distribution $\mu_i$, it satisfies*

$$\mathbb{E}[\langle q, \nu_2^+ \rangle^2] \geq \frac{\|q\|_2^2 \varepsilon^2}{512 \mathsf{gap}^2} \ ,$$

*where $q = \mathcal{A}(x_1, \ldots, x_T)$ is the output of $\mathcal{A}$, $\nu_2^+$ is the second eigenvector of $\mathbb{E}_{\mu_i}[xx^\top]$, and the expectation is taken over the randomness of $i, x_1, \ldots, x_T$, and $\mathcal{A}$.*

*Proof of Corollary ii.I.2.* Without lose of generality, let us assume that $\|q\|_2 = 1$. Let $\nu_2^a, \nu_2^b$ be the second eigenvectors of $\mathbb{E}_{\mu_1}[xx^\top]$ and $\mathbb{E}_{\mu_2}[xx^\top]$ respectively. Let $\nu_1^+, \nu_2^+$ be the two eigenvectors of $\mathbb{E}_{\mu_i}[xx^\top]$, and $\nu_1^-, \nu_2^-$ be the two eigenvectors of $\mathbb{E}_{\mu_{3-i}}[xx^\top]$.

Suppose by way of contradiction that $\mathbb{E}[\langle q, \nu_2^+ \rangle^2] < \frac{\varepsilon^2}{512 \mathsf{gap}^2}$. Then, we shall construct a protocol that can, given samples $x_1, \ldots, x_T$, with success probability at least $3/4$ to tell if the distribution $\mu_i$ is $\mu_1$ or $\mu_2$. However, we cannot distinguish between a fair coin and a $\frac{1}{2} \pm \varepsilon$ biased coin with probability $\geq 3/4$ with fewer than $O(1/\varepsilon^2)$ samples. This will give a contradiction and finish the proof that $\mathbb{E}[\langle q, \nu_2^+ \rangle^2] \geq \frac{\varepsilon^2}{512 \mathsf{gap}^2}$.

We define the following protocol: on input samples $x_1, \ldots, x_T$, run algorithm $\mathcal{A}$ and get output $q = \mathcal{A}(x_1, \ldots, x_T)$; then, declare distribution $\mu_1$ if $\langle q, \nu_2^a \rangle^2 < \frac{\varepsilon^2}{128 \mathsf{gap}^2}$, or distribution $\mu_2$ otherwise.

Using Markov's inequality, with probability at least $\frac{3}{4}$, it satisfies $\langle q, \nu_2^+ \rangle^2 < \frac{\varepsilon^2}{128 \mathsf{gap}^2}$. In this case, we have $\langle q, \nu_1^+ \rangle^2 \geq 1 - \frac{\varepsilon^2}{128 \mathsf{gap}^2} \geq \frac{1}{2}$. By Lemma ii.I.1, we also have $\langle \nu_1^+, \nu_2^- \rangle^2 \geq \frac{\varepsilon^2}{16 \mathsf{gap}^2}$. Using these two together we can derive that

$$\begin{aligned}
\langle q, \nu_2^- \rangle^2 &= (\langle q, \nu_1^+ \rangle \langle \nu_1^+, \nu_2^- \rangle + \langle q, \nu_2^+ \rangle \langle \nu_2^+, \nu_2^- \rangle)^2 \\
&\geq \frac{1}{2} \langle q, \nu_1^+ \rangle^2 \langle \nu_1^+, \nu_2^- \rangle^2 - \langle q, \nu_2^+ \rangle^2 \langle \nu_2^+, \nu_2^- \rangle^2 \geq \frac{\varepsilon^2}{64 \mathsf{gap}^2} - \frac{\varepsilon^2}{128 \mathsf{gap}^2} \geq \frac{\varepsilon^2}{128 \mathsf{gap}^2} \ .
\end{aligned}$$

This means, the protocol we just defined can declare the correct $\mu_i$: indeed, among the two vectors $\{\nu_2^a, \nu_2^b\} = \{\nu_2^+, \nu_2^-\}$, the incorrect one (namely, $\nu_2^-$) will not satisfy $\langle q, \nu_2^- \rangle^2 < \frac{\varepsilon^2}{128 \mathsf{gap}^2}$. $\square$

*Proof of Theorem 6.* Choose $\beta = \sqrt{\frac{1+\delta/\lambda}{2}}$ so $2\beta^2 - 1 = \frac{\delta}{\lambda}$. Let $\varepsilon \in (0, \frac{\delta}{\lambda}]$ be defined such that $2\lambda T/\beta^2 = \Theta(\frac{1}{\varepsilon^2})$. (We can do so because $T \geq \Omega(\frac{\lambda}{\delta^2})$.) We let $a$ and $b$ be the two vectors defined in Lemma ii.I.1 with parameters $\beta, \varepsilon$.

We also denote by $\mathbf{A}_0$ and $\mathbf{A}_1$ respectively the two $2 \times 2$ matrices from Lemma ii.I.1:

$$\mathbf{A}_0 = \begin{pmatrix} \beta^2 & 0 \\ 0 & 1-\beta^2 \end{pmatrix} \quad \text{and} \quad \mathbf{A}_1 = \begin{pmatrix} \beta^2 & 2\varepsilon\beta\sqrt{1-\beta^2} \\ 2\varepsilon\beta\sqrt{1-\beta^2} & 1-\beta^2 \end{pmatrix} .$$

Now, consider the following procedure to generate $T$ random vectors $x_1, \ldots, x_T \in \mathbb{R}^{2k}$. At the beginning, pick a vector $z \in \{0,1\}^k$ uniformly at random from the $2^k$ choices. Then, in each of the $T$ rounds,

1. Pick $y \in [0,1]$ uniformly at random.
2. If $y > k\lambda/\beta^2$, output $x_t = 0$; otherwise continue to the next step. (Note $k\lambda/\beta^2 \in [0,1]$.)
3. Pick $i \in [k]$ uniformly at random.
4. If $z_i = 0$, then pick

$$x_t = 0^{\oplus 2(i-1)} \oplus a \oplus 0^{\oplus 2(k-i)} \text{ w.p. } 1/2 \text{ and } x_t = 0^{\oplus 2(i-1)} \oplus b \oplus 0^{\oplus 2(k-i)} \text{ w.p. } 1/2$$

If $z_i = 1$, then pick

$$x_t = 0^{\oplus(i-1)} \oplus a \oplus 0^{\oplus(k-i)} \text{ w.p. } 1/2 + \varepsilon \text{ and } x_t = 0^{\oplus(i-1)} \oplus b \oplus 0^{\oplus(k-i)} \text{ w.p. } 1/2 - \varepsilon$$

It is clear from the above definition that $x_1, \ldots, x_T$ are generated i.i.d. from a distribution $\mathcal{D}_z$ which is characterized by vector $z$. Since $\|a\|_2 = \|b\|_2 = 1$, we also have $\|x_t\|_2 \leq 1$. By the construction,

$$\mathbb{E}_{x \sim \mathcal{D}_z}[xx^\top] = \bigoplus_{i=1}^{k} \left(\frac{\lambda}{\beta^2} \mathbf{A}_{z_i}\right) ,$$

which implies, according to Lemma ii.I.1, for every possible $z \in \{0,1\}^k$

$$\lambda_k\left(\mathbb{E}_{x \sim \mathcal{D}_z}[xx^\top]\right) \geq \frac{\lambda}{\beta^2}\beta^2 = \lambda \quad \text{and} \quad \lambda_{k+1}\left(\mathbb{E}_{x \sim \mathcal{D}_z}[xx^\top]\right) \leq \frac{\lambda}{\beta^2}(1-\beta^2) \leq \lambda - \delta .$$

Now, in the total number of $T$ rounds, with probability $\geq \frac{1}{2}$, there are less than $2T \cdot \frac{k\lambda}{\beta^2}$ rounds that $x_t \neq 0$. Denote this event as $\mathcal{E}$. If this happens, then at least $7/8$ of $i \in [k]$ are picked less than $16T \cdot \frac{k\lambda}{\beta^2} \cdot \frac{1}{k} = 16\frac{\lambda T}{\beta^2}$ times. Denote this set of indices as $\mathcal{S}_1 \subset [k]$, and let $\mathbf{Q}_T = \mathcal{A}(x_1, \ldots, x_T) \in \mathbb{R}^{2k \times k}$ be the output of the algorithm on this random input sequence $x_1, \ldots, x_T$. (Recall that the algorithm does not know the vector $z \in \{0,1\}^k$). Let $[\mathbf{Q}_T]_i$ be the $i$-th row of $\mathbf{Q}_T$, it holds:

1. $\forall i \in [2k], \|[\mathbf{Q}_T]_i\|_2^2 \leq 1$.
2. $\sum_{i=1}^{2k} \|[\mathbf{Q}_T]_i\|_2^2 = \|\mathbf{Q}_T\|_F^2 = k$.

Therefore, at least $1/4$ of $j \in [k]$ satisfies $\|[\mathbf{Q}_T]_{2j-1}\|_2^2 + \|[\mathbf{Q}_T]_{2j}\|_2^2 \geq \frac{1}{2}$. We denote this set of $j$ as $\mathcal{S}_2 \subset [k]$. Let $\mathcal{S} \stackrel{\text{def}}{=} \mathcal{S}_1 \cap \mathcal{S}_2$, then it holds that $|\mathcal{S}_2| \geq \frac{k}{8}$.

Now we apply Corollary ii.I.2 to each index $i \in \mathcal{S}$ under event $\mathcal{E}$. We can apply Corollary ii.I.2 because $16\lambda T/\beta^2 \leq O(\frac{1}{\varepsilon^2})$ and the algorithm $\mathcal{A}$ has received only $16\lambda T/\beta^2$ nonzero samples under event $\mathcal{E}$. The conclusion of Corollary ii.I.2 implies: for every $i \in \mathcal{S}$:

$$\mathbb{E}[(\nu_{i,2}^+)^\top \mathbf{Q}_T \mathbf{Q}_T^\top \nu_{i,2}^+ \mid \mathcal{E}] \geq \left(\|[\mathbf{Q}_T]_{2i-1}\|_2^2 + \|[\mathbf{Q}_T]_{2i}\|_2^2\right) \cdot \frac{\varepsilon^2}{512(2\beta^2-1)^2} \geq \frac{\varepsilon^2}{1024(2\beta^2-1)^2} .$$

Above, we have denoted by $\nu_{i,2}^+ = 0^{\oplus 2(i-1)} \oplus \nu_2^+ \oplus 0^{\oplus 2(k-i)}$ where $\nu_2^+$ is the second eigenvector of $\mathbf{A}_{z_i}$. It is clear that $\nu_{i,2}^+ \mathbf{Q}_T = (\nu_2^+)^\top (q_1, \ldots, q_k)$ where each $q_j = ([\mathbf{Q}_T]_{2i-1,j}, [\mathbf{Q}_T]_{2i,j})^\top \in \mathbb{R}^2$ so we can apply Corollary ii.I.2 on each of them. Thus, using $\mathbf{Pr}[\mathcal{E}] \geq \frac{1}{2}$ we conclude:

$$\mathbb{E}\left[\sum_{i\in\mathcal{S}}\|(\nu_{i,2}^+)^\top \mathbf{Q}_T\|_2^2\right] \geq \frac{1}{2}\cdot\frac{k}{8}\cdot\frac{\varepsilon^2}{1024(2\beta^2-1)^2} = \frac{\Theta(k)}{\lambda(2\beta^2-1)^2 T/\beta^2} = \Theta\left(\frac{k\lambda}{\delta^2 T}\right) \ .$$

Finally, let $\mathbf{W}$ be the column orthonormal matrix consisting of the $(k+1)$-th till the $2k$-th eigenvectors of matrix $\mathbb{E}_{x\sim\mathcal{D}_z}[xx^\top]$. Since each $\nu_{i,2}^+$ for $i \in [k]$ is in one of the columns of $\mathbf{W}$, the above lower bound implies

$$\mathop{\mathbb{E}}_{z,x_1,\ldots,x_T,\mathbf{Q}_T}\left[\|\mathbf{Q}_T^\top \mathbf{W}\|_F^2\right] \geq \Omega\left(\frac{k\lambda}{\delta^2 T}\right) \ ,$$

where the expectation is over all possible choices $z \in \{0,1\}^k$, the random vectors $x_1,\ldots,x_T$ generated from $\mathcal{D}_z$, and the randomness of the algorithm. By an averaging argument, there exists some $z \in \{0,1\}^k$ such that

$$\mathop{\mathbb{E}}_{x_1,\ldots,x_T\sim\mathcal{D}_z,\mathbf{Q}_T}\left[\|\mathbf{Q}_T^\top \mathbf{W}\|_F^2\right] \geq \Omega\left(\frac{k\lambda}{\delta^2 T}\right) \ . \qquad \square$$

# APPENDIX (PART III)

In this part III of the appendix, we make non-trivial modifications to Appendix I and derive our final theorems. More specifically,

- Appendix iii.J extends our initialization lemma in Appendix i.A to stronger settings;
- Appendix iii.K extends our expectation lemmas in Appendix i.B to stronger settings;
- Appendix iii.L extends our main lemmas in Appendix i.E to stronger settings;
- Appendix iii.M provides the final proofs for our Theorem 1 and Theorem 2;
- Appendix iii.N provides the final proofs for our Theorem 4 and Theorem 5.

**Main Ideas Behind Main Lemmas.**

Our New Main Lemmas in Appendix iii are extensions of Main Lemmas in Appendix i essentially with two additional ingredients.

- UNDER SAMPLING

    This we have already discussed in Section 3 of the main body.

- VARIANCE BOUND

    Recall that our weaker theorems (such as Theorem 1') do not have the factor $\Lambda = \lambda_1 + \cdots + \lambda_k \in [0,1]$ show up in the complexity. This speed-up factor $1/\Lambda$ could be large and has been argued as a very important factor in the total running time by [10]. To achieve this, we need tighter martingale concentrations on our random variables and below we discuss the main intuition.

    Recall that all martingale concentrations for a random process $\{z_t\}_t$ require some upper bound between consecutive variables $|z_t - z_{t+1}|$. If this upper bound holds with probability 1, that is, $|z_t - z_{t+1}|^2 \leq M$, then an Azuma-type of concentration can be naturally used. However, Azuma concentration is not tight: if one knows a better bound on $\mathbb{E}\left[|z_{t+1} - z_t|^2 \mid z_t\right]$, the $M$ term can be replaced with this expected bound a tighter concentration bound can be implied. See for instance the survey [5].

    This same issue also shows up in streaming PCA. Our Lemma Main 1 and Main 2 have adopted a probability-one absolute bound on $|z_t - z_{t+1}|$ because that gives the simplest proof. If one replaces it with a tighter (but very sophisticated) expected bound, the concentration result can be further improved and this improvement translates to a factor $1/\Lambda$ speed-up in the running time on Oja's algorithm in the gap-dependent case (and similarly in the gap-free case). We present such expected bounds in Appendix iii.K.

**Additional Ideas for Theorem 2.**

While the extended main lemmas are sufficient to prove Theorem 1, they lead again to a weaker version of Theorem 2 where $\Lambda_1, \Lambda_2$ are replaced by $(k\Lambda, k\Lambda_2)$.[22] The reason that the factor $k$ shows up is because in *gap-free* case, there is no "local convergence" (i.e. convergence when $\|\mathbf{P}_t\mathbf{Q}(\mathbf{V}^\top\mathbf{P}_t\mathbf{Q})^{-1}\|_2 = O(1)$) as opposite to the gap-case. Therefore, the quantity $\|x_t^\top \mathbf{P}_t\mathbf{Q}(\mathbf{V}^\top\mathbf{P}_t\mathbf{Q})^{-1}\|_2^2$ will stay at $\Omega(k)$ during the entire execution of Oja's algorithm, and never decreased to $O(1)$.

---

[22]This new weaker version is still weaker than Theorem 2 but already stronger than Theorem 2'. We refrain from stating it formally in this paper.

To overcome this issue, we consider an *auxiliary* "under-sampled" objective that has a "local convergence" behavior. As illustrated by Figure 2 on page 57, we define $\mathbf{W}_1$ to consist of all the eigenvectors of eigenvalue $\leq \lambda - \rho/2$, and $\mathbf{V}_1 = \mathbf{W}_1^\perp$ contains all eigenvectors of eigenvalue $> \lambda - \rho/2$. Now, we define the auxiliary objective to be $\mathbf{W}_1^\top \mathbf{P}_t \mathbf{Q} (\mathbf{V}^\top \mathbf{P}_t \mathbf{Q})^{-1}$.

Unlike $\mathbf{Z}$, this new matrix $\mathbf{W}_1$ has an eigengap $\rho/2$ as compared with $\mathbf{V}$, and thus after certain number of iterations, the quantity $\|\mathbf{W}_1^\top \mathbf{P}_t \mathbf{Q} (\mathbf{V}^\top \mathbf{P}_t \mathbf{Q})^{-1}\|_F^2$ can drop to a constant, say $\leq 2$.

Now, let us suddenly "*shift*" our objective to " $\|\mathbf{W}^\top \mathbf{P}_t \mathbf{Q} (\mathbf{V}_1^\top \mathbf{P}_t \mathbf{Q})^{-1}\|_F^2$ ".[23] The crucial observation is that

$$\|x_t^\top \mathbf{P}_t \mathbf{Q} (\mathbf{V}_1^\top \mathbf{P}_t \mathbf{Q})^{-1}\|_2 \leq 1 + \|\mathbf{W}_1^\top \mathbf{P}_t \mathbf{Q} (\mathbf{V}_1^\top \mathbf{P}_t \mathbf{Q})^{-1}\|_2 \leq 1 + \|\mathbf{W}_1^\top \mathbf{P}_t \mathbf{Q} (\mathbf{V}^\top \mathbf{P}_t \mathbf{Q})^{-1}\|_2 \leq 3 \ ,$$

and this is $\sqrt{k}$ times smaller than the previous bound

$$\|x_t^\top \mathbf{P}_t \mathbf{Q} (\mathbf{V}^\top \mathbf{P}_t \mathbf{Q})^{-1}\|_2 \leq O(\sqrt{k}) \ .$$

Therefore, we can now focus on this random process " $\|\mathbf{W}^\top \mathbf{P}_t \mathbf{Q} (\mathbf{V}_1^\top \mathbf{P}_t \mathbf{Q})^{-1}\|_F^2$ ", and it converges a factor $\Omega(k)$ faster than the old quantity $\|\mathbf{W}^\top \mathbf{P}_t \mathbf{Q} (\mathbf{V}^\top \mathbf{P}_t \mathbf{Q})^{-1}\|_F^2$. Finally, Lemma 2.2 tells us that bounding the former (with respect to $\mathbf{V}_1$) as opposed to the old quantity (with respect to $\mathbf{V}$) also implies the convergence of Oja's algorithm, so we are done.

**Additional Ideas for Oja$^{++}$ Algorithm.**

The analysis behind our Oja$^{++}$ requires a more fine-grind argument regarding the alternation between "under-sampling" and "objective shift". As illustrated by Figure 4 on page 60, we consider a logarithmic long sequence $(\mathbf{V}_1, \mathbf{W}_1), (\mathbf{V}_2, \mathbf{W}_2), \ldots$ and shift carefully and only when needed. Since our main goal of Oja$^{++}$ is to remove the factor $O(k)$ in the running time, the spirit behind these shiftings is the same as the "under-sampled" objective as we discussed above. However, in this case we need to shift our objective once per epoch, so we need in total logarithmic many of them.

Also, we need to establish a stronger initialization lemma. Recall that at the beginning of each epoch of Oja$^{++}$, we insert random columns to the working matrix from the previous epoch. Therefore, the initial matrix $\mathbf{Q}_0$ for every epoch (except the first epoch) is not completely fresh random, but consists of two parts $\mathbf{Q}_0 = [\widetilde{\mathbf{Q}}, \mathbf{Q}]$ where $\widetilde{\mathbf{Q}}$ is a (fixed) column orthonormal matrix and $\mathbf{Q}$ is a fresh random gaussian matrix. We therefore need to re-design our initialization lemma for this particular case, especially to bound the quantity $\|x_t^\top \mathbf{Q}_0 (\mathbf{V}^\top \mathbf{Q}_0)^{-1}\|_2 = \|x_t^\top [\widetilde{\mathbf{Q}}, \mathbf{Q}] (\mathbf{V}^\top [\widetilde{\mathbf{Q}}, \mathbf{Q}])^{-1}\|_2$. This will be our main focus in the next sub-section.

---

[23]The matrix $\mathbf{V}_1^\top \mathbf{P}_t \mathbf{Q}$ is not a square matrix anymore because of the shifting, so by inverse we actually mean the Moore-Penrose pseudo-inverse. $\|\mathbf{W}^\top \mathbf{P}_t \mathbf{Q} (\mathbf{V}_1^\top \mathbf{P}_t \mathbf{Q})^{-1}\|_F^2$ is also equal to $\|\mathbf{W}^\top \mathbf{P}_t \mathbf{Q} (\mathbf{Q}^\top \mathbf{P}_t^\top \mathbf{V}_1 \mathbf{V}_1^\top \mathbf{P}_t \mathbf{Q})^{-1/2}\|_F^2$.

### iii.J  Final Initialization Lemmas

We state and prove the following lemma using random matrix theory:

**Lemma iii.J.1.** *There exists constants $C$ such that the following holds: Let $\mathbf{Q}$ in $\mathbb{R}^{N \times n}$, $N \geq n$ be a random matrix with each entry i.i.d. $\mathcal{N}(0,1)$, for every $\varepsilon > 0$,*

$$\mathbf{Pr}[\lambda_{\min}(\mathbf{Q}^\top \mathbf{Q}) \leq \varepsilon^2(\sqrt{N} - \sqrt{n-1})^2] \leq C\varepsilon \log \frac{1}{\varepsilon}$$

*Proof of Lemma iii.J.1.* Theorem 1.1 of [15] states that there exist constants $C', c > 0$ such that
$$\mathbf{Pr}[\lambda_{\min}(\mathbf{Q}^\top \mathbf{Q}) \leq \varepsilon^2(\sqrt{N} - \sqrt{n-1})^2] \leq (C'\varepsilon)^{N-n+1} + e^{-cN} \ .$$
Therefore, if $N > \frac{\log \frac{1}{\varepsilon}}{2c}$ then we are done. In the case when $N \leq \frac{\log \frac{1}{\varepsilon}}{2c}$, we view $\mathbf{Q}$ as a submatrix of an $N \times N$ random matrix $[\mathbf{Q}, \widetilde{\mathbf{Q}}]$ where the entries of $\widetilde{\mathbf{Q}}$ are also i.i.d. generated from $\mathcal{N}(0,1)$. In such a case, we can apply Lemma i.A.2 on $[\mathbf{Q}, \widetilde{\mathbf{Q}}]$ and conclude that for every $\alpha \leq 0$,

$$\mathbf{Pr}\left[\lambda_{\min}([\mathbf{Q}, \widetilde{\mathbf{Q}}]^\top [\mathbf{Q}, \widetilde{\mathbf{Q}}]) \leq \frac{\alpha^2}{N}\right] \leq C'' \cdot \alpha \ .$$

Since $\mathbf{Q}^\top \mathbf{Q}$ is a sub matrix of $[\mathbf{Q}, \widetilde{\mathbf{Q}}]^\top [\mathbf{Q}, \widetilde{\mathbf{Q}}]$, we know that $\lambda_{\min}(\mathbf{Q}^\top \mathbf{Q}) \geq \lambda_{\min}([\mathbf{Q}, \widetilde{\mathbf{Q}}]^\top [\mathbf{Q}, \widetilde{\mathbf{Q}}])$. Choosing $\alpha = \varepsilon N$, we conclude that
$$\mathbf{Pr}[\lambda_{\min}(\mathbf{Q}^\top \mathbf{Q}) \leq \varepsilon^2(\sqrt{N} - \sqrt{n-1})^2] \leq \mathbf{Pr}[\lambda_{\min}(\mathbf{Q}^\top \mathbf{Q}) \leq \varepsilon^2 N] \leq O(\varepsilon \log(1/\varepsilon)) \ . \quad \square$$

We prove the following the following initialization lemma that extends Lemma 5.1. Note that the matrix we are interested now is $[\widetilde{\mathbf{Q}}, \mathbf{Q}]$ where $\mathbf{Q}$ is a random matrix but $\widetilde{\mathbf{Q}}$ is a fixed one. This will allow us to perform "under-sampling" and "objective shift" properly.

**Lemma iii.J.2** (initialization). *Let $\widetilde{\mathbf{Q}} \in \mathbb{R}^{d \times \alpha}$, $\mathbf{Q} \in \mathbb{R}^{d \times \beta}$, $\mathbf{V} \in \mathbb{R}^{d \times k}$ be three matrices such that*

- $\alpha, \beta, k \in \mathbb{N}$ where $\alpha \geq 0$, $\beta \geq 1$, and $\alpha + \beta \leq k \leq d$;
- *each entry of $\mathbf{Q}$ is i.i.d. random from $\mathcal{N}(0,1)$;*
- $\widetilde{\mathbf{Q}}$ *and $\mathbf{V}$ are (column) orthonormal;*
- *either $\widetilde{\mathbf{Q}} = []$ has zero column or $\widetilde{\mathbf{Q}}^\top \mathbf{V} \mathbf{V}^\top \widetilde{\mathbf{Q}} \succeq \frac{1}{2}\mathbf{I}$ holds.*

*Denote by $\mathbf{Z} \in \mathbb{R}^{d \times (d-k)}$ the orthogonal complement of $\mathbf{V}$.*

*Then, for every $p, q \in (0,1)$, every $T \in \mathbb{N}^*$, every set of random vectors $\{x_t\}_{t=1}^T$ with $\|x_t\|_2 \leq 1$, with probability at least $1 - p - q$ over the random choice of $\mathbf{Q}$, the following holds:*

- $\left\| \mathbf{Z}^\top [\widetilde{\mathbf{Q}}, \mathbf{Q}] \left([\widetilde{\mathbf{Q}}, \mathbf{Q}]^\top \mathbf{V}\mathbf{V}^\top [\widetilde{\mathbf{Q}}, \mathbf{Q}]\right)^{-1/2} \right\|_F^2 \leq d \cdot \clubsuit$ ;

- $\displaystyle \Pr_{x_1, \ldots, x_T}\left[\exists i, t \in [T], \left\| x_t^\top \mathbf{Z}\mathbf{Z}^\top \left(\mathbf{\Sigma}/\lambda_{k+1}\right)^{i-1} [\widetilde{\mathbf{Q}}, \mathbf{Q}] \left([\widetilde{\mathbf{Q}}, \mathbf{Q}]^\top \mathbf{V}\mathbf{V}^\top [\widetilde{\mathbf{Q}}, \mathbf{Q}]\right)^{-1/2} \right\|_2^2 \geq \clubsuit \right] \leq q$ ;

- $\left\| \nu_j^\top \mathbf{Z}\mathbf{Z}^\top \mathbf{Q}(\mathbf{Q}^\top \mathbf{V}\mathbf{V}^\top \mathbf{Q})^{-1/2} \right\|_2 \leq \clubsuit$ *for every $j \in [d]$* .

*where*

$$\clubsuit \stackrel{\text{def}}{=} C \cdot \frac{\beta \cdot \log^2(1/p) \cdot \log(Td/q)}{p^2(\sqrt{k-\alpha} - \sqrt{\beta-1})^2} = \widetilde{\Theta}\left(\frac{\beta}{p^2(\sqrt{k-\alpha} - \sqrt{\beta-1})^2}\right) \ .$$

*Proof of Lemma iii.J.2.* Let $\mathbf{A}^2 \stackrel{\text{def}}{=} \widetilde{\mathbf{Q}}^\top \mathbf{V}\mathbf{V}^\top \widetilde{\mathbf{Q}} \succeq \frac{1}{2}\mathbf{I} \in \mathbb{R}^{\alpha \times \alpha}$ and by $\mathbf{B}^\top \stackrel{\text{def}}{=} \mathbf{A}^{-1}\widetilde{\mathbf{Q}}^\top \mathbf{V}\mathbf{V}^\top \in \mathbb{R}^{\alpha \times d}$. We have

$$\mathbf{B}^\top \mathbf{B} = \mathbf{I} \quad \text{and} \quad \mathbf{B}\mathbf{B}^\top = \mathbf{V}^\top \left(\mathbf{V}^\top \mathbf{Q}[\widetilde{\mathbf{Q}}^\top \mathbf{V}\mathbf{V}^\top \widetilde{\mathbf{Q}}]^{-1}\mathbf{Q}^\top \mathbf{V}\right) \mathbf{V} \ .$$

Since $\mathbf{V}^\top\mathbf{Q}[\widetilde{\mathbf{Q}}^\top\mathbf{V}\mathbf{V}^\top\widetilde{\mathbf{Q}}]^{-1/2}$ is a column orthonormal matrix, we can always find $\mathbf{C}$ such that
$$\mathbf{V}\mathbf{V}^\top - \mathbf{B}\mathbf{B}^\top = \mathbf{C}\mathbf{C}^\top \succeq 0, \quad \mathbf{C}^\top\mathbf{B} = 0, \quad \mathbf{C} \in \mathbb{R}^{d\times(k-\alpha)} .$$
Therefore,
$$\left([\widetilde{\mathbf{Q}},\mathbf{Q}]^\top\mathbf{V}\mathbf{V}^\top[\widetilde{\mathbf{Q}},\mathbf{Q}]\right) = \begin{pmatrix} \widetilde{\mathbf{Q}}^\top\mathbf{V}\mathbf{V}^\top\widetilde{\mathbf{Q}} & \widetilde{\mathbf{Q}}^\top\mathbf{V}\mathbf{V}^\top\mathbf{Q} \\ \mathbf{Q}^\top\mathbf{V}\mathbf{V}^\top\widetilde{\mathbf{Q}} & \mathbf{Q}^\top\mathbf{V}\mathbf{V}^\top\mathbf{Q} \end{pmatrix}$$
$$= \begin{pmatrix} \mathbf{A} & \\ & \mathbf{I} \end{pmatrix}\begin{pmatrix} \mathbf{I} & \mathbf{B}^\top\mathbf{Q} \\ \mathbf{Q}^\top\mathbf{B} & \mathbf{Q}^\top(\mathbf{B}\mathbf{B}^\top + \mathbf{C}\mathbf{C}^\top)\mathbf{Q} \end{pmatrix}\begin{pmatrix} \mathbf{A} & \\ & \mathbf{I} \end{pmatrix}$$

Let us denote $\mathbf{Q}_b = \mathbf{B}^\top\mathbf{Q} \in \mathbb{R}^{\alpha\times\beta}$, $\mathbf{Q}_c = \mathbf{C}^\top\mathbf{Q} \in \mathbb{R}^{(k-\alpha)\times\beta}$, $\sigma = \sigma_{\min}(\mathbf{Q}_c^\top\mathbf{Q}_c)$, we then have:
$$\left([\widetilde{\mathbf{Q}},\mathbf{Q}]^\top\mathbf{V}\mathbf{V}^\top[\widetilde{\mathbf{Q}},\mathbf{Q}]\right) \succeq \begin{pmatrix} \mathbf{A} & \\ & \mathbf{I} \end{pmatrix}\begin{pmatrix} \mathbf{I} & \mathbf{Q}_b \\ \mathbf{Q}_b^\top & \mathbf{Q}_b^\top\mathbf{Q}_b + \sigma\mathbf{I} \end{pmatrix}\begin{pmatrix} \mathbf{A} & \\ & \mathbf{I} \end{pmatrix}$$

Which implies that
$$\left([\widetilde{\mathbf{Q}},\mathbf{Q}]^\top\mathbf{V}\mathbf{V}^\top[\widetilde{\mathbf{Q}},\mathbf{Q}]\right)^{-1} \preceq \begin{pmatrix} \mathbf{A}^{-1} & \\ & \mathbf{I} \end{pmatrix}\begin{pmatrix} \mathbf{I}+\sigma^{-1}\mathbf{Q}_b\mathbf{Q}_b^\top & -\sigma^{-1}\mathbf{Q}_b \\ -\sigma^{-1}\mathbf{Q}_b^\top & \sigma^{-1}\mathbf{I} \end{pmatrix}\begin{pmatrix} \mathbf{A}^{-1} & \\ & \mathbf{I} \end{pmatrix}$$
$$\preceq 2\begin{pmatrix} \mathbf{I}+\sigma^{-1}\mathbf{Q}_b\mathbf{Q}_b^\top & -\sigma^{-1}\mathbf{Q}_b \\ -\sigma^{-1}\mathbf{Q}_b^\top & \sigma^{-1}\mathbf{I} \end{pmatrix}$$
$$= 2\begin{pmatrix} \mathbf{I} & 0 \\ 0 & 0 \end{pmatrix} + 2\sigma^{-1}[\mathbf{Q}_b^\top,-\mathbf{I}]^\top[\mathbf{Q}_b^\top,-\mathbf{I}] . \qquad \text{(iii.J.1)}$$

Above, the first spectral dominance uses $\begin{pmatrix} \mathbf{I} & \mathbf{Q}_b \\ \mathbf{Q}_b^\top & \mathbf{Q}_b^\top\mathbf{Q}_b + \sigma\mathbf{I} \end{pmatrix}^{-1} = \begin{pmatrix} \mathbf{I}+\sigma^{-1}\mathbf{Q}_b\mathbf{Q}_b^\top & -\sigma^{-1}\mathbf{Q}_b \\ -\sigma^{-1}\mathbf{Q}_b^\top & \sigma^{-1}\mathbf{I} \end{pmatrix}$;
the second spectral dominance uses $\mathbf{A}^2 \preceq \frac{1}{2}\mathbf{I}$

Now consider any fixed $y \in \mathbb{R}^d$ with $\|y\|_2 \leq 1$. We have
$$y^\top[\widetilde{\mathbf{Q}},\mathbf{Q}]\left([\widetilde{\mathbf{Q}},\mathbf{Q}]^\top\mathbf{V}\mathbf{V}^\top[\widetilde{\mathbf{Q}},\mathbf{Q}]\right)^{-1}[\widetilde{\mathbf{Q}},\mathbf{Q}]^\top y$$
$$= \left\| y^\top\mathbf{V}\mathbf{V}^\top[\widetilde{\mathbf{Q}},\mathbf{Q}]\left([\widetilde{\mathbf{Q}},\mathbf{Q}]^\top\mathbf{V}\mathbf{V}^\top[\widetilde{\mathbf{Q}},\mathbf{Q}]\right)^{-1/2} + y^\top\mathbf{Z}\mathbf{Z}^\top[\widetilde{\mathbf{Q}},\mathbf{Q}]\left([\widetilde{\mathbf{Q}},\mathbf{Q}]^\top\mathbf{V}\mathbf{V}^\top[\widetilde{\mathbf{Q}},\mathbf{Q}]\right)^{-1/2} \right\|_2^2$$
$$\leq 2\left\| y^\top\mathbf{V} \right\|_2^2 \left\| \mathbf{V}^\top[\widetilde{\mathbf{Q}},\mathbf{Q}]\left([\widetilde{\mathbf{Q}},\mathbf{Q}]^\top\mathbf{V}\mathbf{V}^\top[\widetilde{\mathbf{Q}},\mathbf{Q}]\right)^{-1/2} \right\|_2^2 + 2\left\| y^\top\mathbf{Z}\mathbf{Z}^\top[\widetilde{\mathbf{Q}},\mathbf{Q}]\left([\widetilde{\mathbf{Q}},\mathbf{Q}]^\top\mathbf{V}\mathbf{V}^\top[\widetilde{\mathbf{Q}},\mathbf{Q}]\right)^{-1/2} \right\|_2^2$$
$$\leq 2 + 2y^\top\mathbf{Z}\mathbf{Z}^\top[\widetilde{\mathbf{Q}},\mathbf{Q}]\left([\widetilde{\mathbf{Q}},\mathbf{Q}]^\top\mathbf{V}\mathbf{V}^\top[\widetilde{\mathbf{Q}},\mathbf{Q}]\right)^{-1}[\widetilde{\mathbf{Q}},\mathbf{Q}]^\top\mathbf{Z}\mathbf{Z}^\top y . \qquad \text{(iii.J.2)}$$

Above, the last inequality uses $\|y^\top\mathbf{V}\|_2 \leq 1$ and the fact that $\|\mathbf{A}(\mathbf{A}^\top\mathbf{A})^{-1/2}\|_2 \leq 1$ for every matrix $\mathbf{A}$ that has full column rank. Next, using inequality (iii.J.1), we can bound
$$\frac{1}{2}y^\top\mathbf{Z}\mathbf{Z}^\top[\widetilde{\mathbf{Q}},\mathbf{Q}]\left([\widetilde{\mathbf{Q}},\mathbf{Q}]^\top\mathbf{V}\mathbf{V}^\top[\widetilde{\mathbf{Q}},\mathbf{Q}]\right)^{-1}[\widetilde{\mathbf{Q}},\mathbf{Q}]^\top\mathbf{Z}\mathbf{Z}^\top y$$
$$\leq \|y^\top\mathbf{Z}\mathbf{Z}^\top\widetilde{\mathbf{Q}}\|_2^2 + \sigma^{-1}\|y^\top\mathbf{Z}\mathbf{Z}^\top(\widetilde{\mathbf{Q}}\mathbf{Q}_b^\top - \mathbf{Q})]\|_2^2$$
$$\leq 1 + \frac{2\|\mathbf{Q}_b^\top\widetilde{\mathbf{Q}}^\top\mathbf{Z}\mathbf{Z}^\top y\|_2^2 + 2\|\mathbf{Q}^\top\mathbf{Z}\mathbf{Z}^\top y\|_2^2}{\sigma} . \qquad \text{(iii.J.3)}$$

Note that $\mathbf{Q}_b = \mathbf{B}^\top\mathbf{Q} \in \mathbb{R}^{\alpha\times\beta}$, $\mathbf{Q}_c = \mathbf{C}^\top\mathbf{Q} \in \mathbb{R}^{(k-\alpha)\times\beta}$, and $\mathbf{Z}^\top\mathbf{Q} \in \mathbb{R}^{(d-k)\times\beta}$ are independent of each other, because the entries of $\mathbf{Q}$ remain independent after rotation and the columns of $\mathbf{B},\mathbf{C},\mathbf{Z}$ are pairwise orthogonal. Therefore, letting $u_1 \stackrel{\text{def}}{=} \mathbf{Q}_b^\top\widetilde{\mathbf{Q}}^\top\mathbf{Z}\mathbf{Z}^\top y$ and $u_2 \stackrel{\text{def}}{=} \mathbf{Q}^\top\mathbf{Z}\mathbf{Z}^\top y \in \mathbb{R}^\beta$, we have that

- $u_1 \in \mathbb{R}^\beta$ is a random vector with entries i.i.d. in $\mathcal{N}(0,\|\mathbf{B}\widetilde{\mathbf{Q}}^\top\mathbf{Z}\mathbf{Z}^\top y\|_2^2)$;
- $u_2 \in \mathbb{R}^\beta$ is a random vector with entries i.i.d. in $\mathcal{N}(0,\|\mathbf{Z}\mathbf{Z}^\top y\|_2^2)$;

- $u_1$, $u_2$ and $\sigma = \sigma_{\min}(\mathbf{Q}_c^\top \mathbf{Q}_c)$ are independent variables.

Since $\|\mathbf{B}\widetilde{\mathbf{Q}}^\top \mathbf{Z}\mathbf{Z}^\top y\|_2 \leq 1$ and $\|\mathbf{Z}\mathbf{Z}^\top y\|_2^2 \leq 1$, it implies (using tail bound for chi-squared distribution) that
$$\forall s \geq 4 \colon \mathbf{Pr}[\|u_1\|_2^2 \geq s\beta], \mathbf{Pr}[\|u_2\|_2^2 \geq s\beta] \leq (se^{1-s})^{\beta/2} \leq e^{-s/6} \ .$$
Putting them into (iii.J.3) and then altogether to (iii.J.2), we claim that under event $\mathcal{C}_\sigma \stackrel{\text{def}}{=} \{\sigma_{\min}(\mathbf{Q}^\top \mathbf{C}\mathbf{C}^\top \mathbf{Q}) \geq \sigma\}$, we have
$$\mathbf{Pr}\left[ y^\top [\widetilde{\mathbf{Q}}, \mathbf{Q}] \left( [\widetilde{\mathbf{Q}}, \mathbf{Q}]^\top \mathbf{V}\mathbf{V}^\top [\widetilde{\mathbf{Q}}, \mathbf{Q}] \right)^{-1} [\widetilde{\mathbf{Q}}, \mathbf{Q}]^\top y \geq 6 + \frac{16\beta s}{\sigma} \ \Big| \ \mathcal{C}_\sigma \right] \leq 2e^{-s/6} \ . \qquad \text{(iii.J.4)}$$

**Final Probability Arguments.** Let us now take $\sigma = \frac{p^2(\sqrt{k-\alpha} - \sqrt{\beta-1})^2}{c_2 \cdot \log^2(1/p)}$ for a large enough constant $c_2$, we have $\mathbf{Pr}[C_\sigma] \geq 1 - p/2$ according to Lemma iii.J.1. We also take $s = \Theta(\log(Td/q))$, and it satisfies
$$6 + \frac{16\beta s}{\sigma} = \Theta\left(\frac{\beta \cdot \log^2(1/p) \cdot \log(Td/q)}{p^2(\sqrt{k-\alpha} - \sqrt{\beta-1})^2}\right) = \clubsuit \ .$$
We now apply (iii.J.4) multiple times in two different manners

- We can apply (iii.J.4) with $y = \nu_1, \nu_2, \ldots, \nu_d$ which are the eigenvector of $\mathbf{\Sigma}$. By union bound, we have $\mathbf{Pr}[\mathcal{C}_1 | \mathcal{C}_\sigma] \geq 1 - q$ where $\mathcal{C}_1$ is the following event
$$C_1 \stackrel{\text{def}}{=} \left\{ \begin{array}{l} \left\| \mathbf{Z}^\top [\widetilde{\mathbf{Q}}, \mathbf{Q}] \left( [\widetilde{\mathbf{Q}}, \mathbf{Q}]^\top \mathbf{V}\mathbf{V}^\top [\widetilde{\mathbf{Q}}, \mathbf{Q}] \right)^{-1/2} \right\|_F^2 \leq d \cdot \clubsuit \ , \text{ and} \\ \forall j \in [d] \colon \left\| \nu_j \mathbf{Z}\mathbf{Z}^\top [\widetilde{\mathbf{Q}}, \mathbf{Q}] \left( [\widetilde{\mathbf{Q}}, \mathbf{Q}]^\top \mathbf{V}\mathbf{V}^\top [\widetilde{\mathbf{Q}}, \mathbf{Q}] \right)^{-1/2} \right\|_2^2 \leq \clubsuit \end{array} \right\}$$

- Define event
$$\mathcal{C}_2 = \left\{ \exists i \in [T], \exists t \in [T], \left\| x_t^\top \mathbf{Z}\mathbf{Z}^\top (\mathbf{\Sigma}/\lambda_{k+1})^{i-1} [\widetilde{\mathbf{Q}}, \mathbf{Q}] \left( [\widetilde{\mathbf{Q}}, \mathbf{Q}]^\top \mathbf{V}\mathbf{V}^\top [\widetilde{\mathbf{Q}}, \mathbf{Q}] \right)^{-1/2} \right\|_2 \geq \clubsuit \right\} \ .$$

Then, we apply (iii.J.4) multiple times each with $y = x_t^\top \mathbf{Z}\mathbf{Z}^\top (\mathbf{\Sigma}/\lambda_{k+1})^{i-1}$ (and we can do so because $\|y\|_2 \leq 1$). By union bound, we have for every fixed $x_1, \ldots, x_T$, it satisfies $\mathbf{Pr}_{\mathbf{Q}}[\mathcal{C}_2 | \mathcal{C}_\sigma, x_1, \ldots, x_T] \leq q^2$. Denoting by $\mathbb{1}_{\mathcal{C}_2}$ the indicator function of event $\mathcal{C}_2$, then
$$\mathbf{Pr}_{\mathbf{Q}}\left[ \mathbf{Pr}_{x_1,\ldots,x_T}[\mathcal{C}_2 | \mathbf{Q}] \geq q \ \Big| \ \mathcal{C}_\sigma \right] \leq \frac{1}{q} \mathbb{E}_{\mathbf{Q}}\left[ \mathbf{Pr}_{x_1,\ldots,x_T}[\mathcal{C}_2 | \mathbf{Q}] \ \Big| \ \mathcal{C}_\sigma \right]$$
$$= \frac{1}{q} \mathbb{E}_{\mathbf{Q}}\left[ \mathbb{E}_{x_1,\ldots,x_T}[\mathbb{1}_{\mathcal{C}_2} | \mathbf{Q}] \ \Big| \ \mathcal{C}_\sigma \right]$$
$$= \frac{1}{q} \mathbb{E}_{x_1,\ldots,x_T}\left[ \mathbb{E}_{\mathbf{Q}}[\mathbb{1}_{\mathcal{C}_2} | \mathcal{C}_\sigma, x_1, \ldots, x_T] \right]$$
$$= \frac{1}{q} \mathbb{E}_{x_1,\ldots,x_T}\left[ \mathbf{Pr}_{\mathbf{Q}}[\mathcal{C}_2 | \mathcal{C}_\sigma, x_1, \ldots, x_T] \right] \leq q \ .$$

Above, the first inequality uses Markov's bound.

In sum, we just derived that $\mathbf{Pr}_{\mathbf{Q}}[\mathcal{C}_1 | \mathcal{C}_\sigma] \geq 1 - q$ and $\mathbf{Pr}_{\mathbf{Q}}\left[ \mathbf{Pr}_{x_1,\ldots,x_T}[\mathcal{C}_2 | \mathbf{Q}] \geq q \right] \leq q$. Applying union bound again we have
$$\mathbf{Pr}_{\mathbf{Q}}\left[ \mathcal{C}_1 \wedge \mathbf{Pr}_{x_1,\ldots,x_T}[\mathcal{C}_2 | \mathbf{Q}] \leq q \right] \geq 1 - p - 2q \ .$$
This finishes the proof of Lemma iii.J.2. $\square$

## iii.K  Final Expectation Lemmas

This section extends Appendix i.B to provide better bounds regarding the expected behaviors of the random process we are interested. Recall that if $\mathbf{X} \in \mathbb{R}^{d \times r}$ is a generic matrix (either $\mathbf{X} = \mathbf{W}$, $\mathbf{X} = \mathbf{Z}$, or $\mathbf{X} = [w]$ for some vector $w$), we have the following notions from Section 6.

$$\begin{array}{ll}
\mathbf{L}_t = \mathbf{P}_t \mathbf{Q} (\mathbf{V}^\top \mathbf{P}_t \mathbf{Q})^{-1} \in \mathbb{R}^{d \times k} & \mathbf{R}'_t = \mathbf{X}^\top x_t x_t^\top \mathbf{L}_{t-1} \in \mathbb{R}^{r \times k} \\
\mathbf{S}_t = \mathbf{X}^\top \mathbf{L}_t \in \mathbb{R}^{r \times k} & \mathbf{H}'_t = \mathbf{V}^\top x_t x_t^\top \mathbf{L}_{t-1} \in \mathbb{R}^{k \times k}
\end{array}$$

Our next lemma is a stronger version of Lemma 6.1. It provide tighter expected bounds by introducing an additional factor $\Gamma$, and introduce new bounds regarding the variance terms $\mathbb{E}\left[\left|\mathbf{Tr}(\mathbf{S}_t^\top \mathbf{S}_t) - \mathbf{Tr}(\mathbf{S}_{t-1}\mathbf{S}_{t-1}^\top)\right|^2\right]$.

**Lemma iii.K.1.** *For every $t \in [T]$, For every $t \in [T]$, let $\mathcal{C}_{\leq t}$ be any event that depends on random $x_1, \ldots, x_t$ and implies*

$$\forall j \in [d], \|\nu_j^\top \mathbf{L}_{t-1}\|_2 \leq \phi_t, \quad \|x_t^\top \mathbf{L}_{t-1}\|_2 = \|x_t^\top \mathbf{P}_{t-1}\mathbf{Q}(\mathbf{V}^\top \mathbf{P}_{t-1}\mathbf{Q})^{-1}\|_2 \leq \phi_t \quad \text{where} \quad \eta_t \phi_t \leq \frac{1}{2},$$

*and $\mathbb{E}[x_t x_t^\top \mid \mathcal{F}_{\leq t-1}, \mathcal{C}_{\leq t}] = \mathbf{\Sigma} + \mathbf{\Delta}$. Let $m$ be the largest integer such that $\lambda_{k+m} > \lambda_k - \rho$, and*

$$\Gamma \stackrel{\text{def}}{=} \min\left\{\sum_{i=1}^k \lambda_i + \phi_t^2 \sum_{j=k+1}^{k+m} \lambda_j + \|\mathbf{\Delta}\|_2, \; 1\right\} \quad \text{and} \quad \chi \stackrel{\text{def}}{=} k + m\phi_t^2$$

*. We have:*

*(a) If $\mathbf{X} = [w] \in \mathbb{R}^{d \times 1}$ where $w$ is a vector with Euclidean norm at most 1,*

$$\mathbb{E}\left[\mathbf{Tr}(\mathbf{S}_t^\top \mathbf{S}_t) \mid \mathcal{F}_{\leq t-1}, \mathcal{C}_{\leq t}\right] \leq (1 - \eta_t \lambda_k + 14\Gamma \eta_t^2 \phi_t^2) \mathbf{Tr}(\mathbf{S}_{t-1}\mathbf{S}_{t-1}^\top) + 10\Gamma \eta_t^2 \phi_t^2 + \frac{\eta_t}{\lambda_k} \|w^\top \mathbf{\Sigma} \mathbf{L}_{t-1}\|_2^2$$
$$+ 2\eta_t \|\mathbf{\Delta}\|_2 \left([\mathbf{Tr}(\mathbf{S}_{t-1}^\top \mathbf{S}_{t-1})]^{1/2} + \mathbf{Tr}(\mathbf{S}_{t-1}^\top \mathbf{S}_{t-1})\right)\left(1 + [\mathbf{Tr}(\mathbf{Z}^\top \mathbf{L}_{t-1}\mathbf{L}_{t-1}^\top \mathbf{Z})]^{1/2}\right)$$

*(b) If $\mathbf{X} = [w] \in \mathbb{R}^{d \times 1}$ where $w$ is a vector with Euclidean norm at most 1,*

$$\mathbb{E}\left[\left|\mathbf{Tr}(\mathbf{S}_t^\top \mathbf{S}_t) - \mathbf{Tr}(\mathbf{S}_{t-1}\mathbf{S}_{t-1}^\top)\right|^2 \mid \mathcal{F}_{t-1}, \mathcal{C}_{\leq t}\right]$$
$$\leq 243\Gamma \eta_t^2 \phi_t^2 \mathbf{Tr}(\mathbf{S}_{t-1}^\top \mathbf{S}_{t-1})^2 + 12\Gamma \eta_t^2 \phi_t^2 \mathbf{Tr}(\mathbf{S}_{t-1}^\top \mathbf{S}_{t-1}) + 300\Gamma \eta_t^4 \phi_t^4$$

*(c) If $\mathbf{X} = \mathbf{W}$,*

$$\mathbb{E}\left[\mathbf{Tr}(\mathbf{S}_t^\top \mathbf{S}_t) \mid \mathcal{F}_{t-1}, \mathcal{C}_{\leq t}\right]$$
$$\leq \left(1 - 2\eta_t \rho + 12\Gamma \eta_t^2 \phi_t^2 + \eta_t^2(6\phi_t + 8)\lambda_{k+1}\right) \mathbf{Tr}(\mathbf{S}_{t-1}\mathbf{S}_{t-1}^\top) + 10\Gamma \eta_t^2(2\phi_t + 8)$$
$$+ 2\eta_t \|\mathbf{\Delta}\|_2 \Big(\eta_t(4 + \phi_t)\chi + [\mathbf{Tr}(\mathbf{L}_{t-1}^\top \mathbf{Z}\mathbf{Z}^\top \mathbf{L}_{t-1})]^{3/2} + (5 + 4\eta_t)\mathbf{Tr}(\mathbf{L}_{t-1}^\top \mathbf{Z}\mathbf{Z}^\top \mathbf{L}_{t-1})$$
$$+ [\mathbf{Tr}(\mathbf{L}_{t-1}^\top \mathbf{Z}\mathbf{Z}^\top \mathbf{L}_{t-1})]^{1/2}\Big) \; .$$

*(d) If $\mathbf{X} = \mathbf{W}$,*

$$\mathbb{E}[|\mathbf{Tr}(\mathbf{S}_t^\top \mathbf{S}_t) - \mathbf{Tr}(\mathbf{S}_{t-1}\mathbf{S}_{t-1}^\top)|_2^2 \mid \mathcal{F}_{t-1}, \mathcal{C}_{\leq t}]$$
$$\leq 192\Gamma \eta_t^2 \phi_t^2 \mathbf{Tr}(\mathbf{S}_{t-1}\mathbf{S}_{t-1}^\top)^2 + 4\eta_t^2 (4\phi_t + 8)^2 \lambda_{k+1} \mathbf{Tr}(\mathbf{S}_{t-1}\mathbf{S}_{t-1}^\top)$$
$$+ 192\Gamma \eta_t^4 \phi_t^2 + \|\mathbf{\Delta}\|_2 \cdot 4\eta_t^2 (4\phi_t + 8)^2 (\chi + \mathbf{Tr}(\mathbf{S}_{t-1}^\top \mathbf{S}_{t-1})) \; .$$

*Proof.* The proof of the first two cases rely on the following tighter upper bounds when $\mathbf{X} = [w]$:

$$\mathbb{E}\left[\|\mathbf{H}'_t\|_2^2 \mid \mathcal{F}_{\leq t-1}, \mathcal{C}_{\leq t}\right] \leq \phi_t^2 \,\mathbb{E}\left[\|x_t \mathbf{V}\|_2^2 \mid \mathcal{F}_{\leq t-1}, \mathcal{C}_{\leq t}\right] \leq \min\left\{\left(\sum_{i=1}^k \lambda_i + \|\mathbf{\Delta}\|_2\right), 1\right\} \phi_t^2 \leq \Gamma \phi_t^2$$

$$\mathbb{E}\left[\|\mathbf{R}'_t\|_2^2 \mid \mathcal{F}_{\leq t-1}, \mathcal{C}_{\leq t}\right] \leq \phi_t^2 \,\mathbb{E}\left[\|x_t \mathbf{X}\|_2^2 \mid \mathcal{F}_{\leq t-1}, \mathcal{C}_{\leq t}\right] \leq \min\left\{\left(\sum_{i=1}^k \lambda_i + \|\mathbf{\Delta}\|_2\right), 1\right\} \phi_t^2 \leq \Gamma \phi_t^2$$
(iii.K.1)

as opposed to $\phi_t^2$ that we have used in the past.

The proof of the last two cases rely on different upper bounds for $\mathbf{X} = \mathbf{W}$. We introduce some notations that shall be only used in this proof. Let $\mathbf{Y} \stackrel{\text{def}}{=} [\nu_{k+1}, \ldots, \nu_{k+m}] \in \mathbb{R}^{d \times m}$ be the matrix consisting of all the eigenvectors of $\mathbf{\Sigma}$ with eigenvalues $\lambda_{k+1}, \ldots, \lambda_{k+m}$. In this notation, $\mathbf{Z}\mathbf{Z}^\top = \mathbf{Y}\mathbf{Y}^\top + \mathbf{W}\mathbf{W}^\top$. We also denote by $\widetilde{\mathbf{V}} = [\mathbf{V}, \mathbf{Y}]$. Let $\Lambda_1 = \sum_{j=1}^k \lambda_j, \Lambda_2 = \sum_{j=k+1}^{k+m} \lambda_j$, $\Lambda = \Lambda_1 + \Lambda_2 \phi_t^2 \leq \Gamma$.

We make three quick observations:

$$\begin{cases} \mathbf{Tr}(\mathbf{L}_{t-1}^\top \widetilde{\mathbf{V}} \mathbf{\Sigma}_{\leq k+m} \widetilde{\mathbf{V}}^\top \mathbf{L}_t) = \sum_{j=1}^k \lambda_j + \sum_{j=k+1}^{k+m} \lambda_j \|\nu_j^\top \mathbf{L}_t\|_2^2 \leq \Lambda_1 + \Lambda_2 \phi_t^2 = \Lambda \ , \\ \mathbf{Tr}(\mathbf{L}_{t-1}^\top \widetilde{\mathbf{V}} \widetilde{\mathbf{V}}^\top \mathbf{L}_t) = k + \sum_{j=k+1}^{k+m} \|\nu_j \mathbf{L}_t\|_2^2 \leq k + m\phi_t^2 = \chi \ , \text{ and} \\ \mathbf{Tr}(\mathbf{L}_{t-1}^\top \mathbf{L}_t) = \mathbf{Tr}(\mathbf{L}_{t-1}^\top (\widetilde{\mathbf{V}} \widetilde{\mathbf{V}}^\top + \mathbf{W}\mathbf{W}^\top)\mathbf{L}_t) \leq \chi + \mathbf{Tr}(\mathbf{S}_{t-1}^\top \mathbf{S}_{t-1}) \ . \end{cases}$$
(iii.K.2)

Therefore, we have:

$$\mathbb{E}[\|\mathbf{H}'_t\|_2^2 \mid \mathcal{F}_{\leq t-1}, \mathcal{C}_{\leq t}] \leq \phi_t^2 \cdot \mathbb{E}[\|x_t^\top \mathbf{V}\|_F^2 \mid \mathcal{F}_{\leq t-1}, \mathcal{C}_{\leq t}] \leq \Gamma \phi_t^2$$

$$\mathbb{E}[\|\mathbf{R}'_t\|_2^2 \mid \mathcal{F}_{\leq t-1}, \mathcal{C}_{\leq t}] \leq \mathbb{E}[\|x_t^\top \mathbf{L}_{t-1}\|_2^2 \mid \mathcal{F}_{\leq t-1}, \mathcal{C}_{\leq t}] = \mathbf{Tr}(\mathbf{L}_{t-1}^\top \mathbf{\Sigma} \mathbf{L}_{t-1}) + \mathbf{Tr}(\mathbf{L}_{t-1}^\top \mathbf{\Delta} \mathbf{L}_{t-1})$$
$$\leq \mathbf{Tr}(\mathbf{L}_{t-1}^\top (\widetilde{\mathbf{V}} \mathbf{\Sigma}_{\leq k+m} \widetilde{\mathbf{V}}^\top + \mathbf{W} \mathbf{\Sigma}_{>k+m} \mathbf{W}^\top) \mathbf{L}_{t-1}) + \|\mathbf{\Delta}\|_2 \mathbf{Tr}(\mathbf{L}_{t-1}^\top \mathbf{L}_{t-1})$$
$$\leq \Lambda + \lambda_{k+1} \|\mathbf{W}^\top \mathbf{L}_{t-1}\|_F^2 + \|\mathbf{\Delta}\|_2 \cdot (\chi + \mathbf{Tr}(\mathbf{S}_{t-1}^\top \mathbf{S}_{t-1}))$$
$$= \Lambda + \lambda_{k+1} \mathbf{Tr}(\mathbf{S}_{t-1}^\top \mathbf{S}_{t-1}) + \|\mathbf{\Delta}\|_2 \cdot (\chi + \mathbf{Tr}(\mathbf{S}_{t-1}^\top \mathbf{S}_{t-1})) \ . \quad \text{(iii.K.3)}$$

We are now ready to prove our four cases individually.

(a) This follows from almost the same proof of Corollary 6.1-(c), except that one can replace the use of (i.B.10) with the following (owing to (iii.K.1))

$$\mathbb{E}\left[\mathbf{Tr}(\mathbf{S}_t^\top \mathbf{S}_t) \mid \mathcal{F}_{\leq t-1}, \mathcal{C}_{\leq t}\right] \leq (1 - 2\eta_t \lambda_k + 14\Gamma \eta_t^2 \phi_t^2) \mathbf{Tr}(\mathbf{S}_{t-1} \mathbf{S}_{t-1}^\top) + 10\Gamma \eta_t^2 \phi_t^2$$
$$- 2\eta_t \mathbf{Tr}(\mathbf{S}_{t-1}^\top \mathbf{S}_{t-1} \mathbf{V}^\top \mathbf{\Delta} \mathbf{L}_{t-1}) + 2\eta_t \mathbf{Tr}(\mathbf{S}_{t-1}^\top w^\top \mathbf{\Delta} \mathbf{L}_{t-1}) + 2\eta_t \mathbf{Tr}(\mathbf{S}_{t-1}^\top w^\top \mathbf{\Sigma} \mathbf{L}_{t-1}) \ .$$

(b) This follows directly from Lemma i.B.1-(b) and (iii.K.1).

(c) The exact same first half of the proof of Lemma i.B.1 gives (i.B.1) which using $\|\mathbf{H}'_t\|_2 \leq \phi_t$ gives

$$\mathbf{Tr}(\mathbf{S}_t^\top \mathbf{S}_t) \leq \mathbf{Tr}(\mathbf{S}_{t-1} \mathbf{S}_{t-1}^\top) - 2\eta_t \mathbf{Tr}(\mathbf{S}_{t-1}^\top \mathbf{S}_{t-1} \mathbf{H}'_t) + 2\eta_t \mathbf{Tr}(\mathbf{S}_{t-1}^\top \mathbf{R}'_t)$$
$$+ 4\eta_t^2 \phi_t \left|\mathbf{Tr}(\mathbf{S}_{t-1}^\top \mathbf{R}'_t)\right| + 12\eta_t^2 \|\mathbf{H}'_t\|_2^2 \mathbf{Tr}(\mathbf{S}_{t-1} \mathbf{S}_{t-1}^\top) + 8\eta_t^2 \|\mathbf{R}'_t\|_2^2 \ . \text{(iii.K.4)}$$

This time, we upper bound

$$|\mathbf{Tr}(\mathbf{S}_{t-1}^\top \mathbf{R}'_t)| = |\mathbf{Tr}(\mathbf{S}_{t-1}^\top \mathbf{W}^\top x_t x_t^\top \mathbf{L}_{t-1})| = |\mathbf{Tr}(\mathbf{S}_{t-1}^\top \mathbf{W}^\top x_t x_t^\top (\widetilde{\mathbf{V}}\widetilde{\mathbf{V}}^\top + \mathbf{W}\mathbf{W}^\top) \mathbf{L}_{t-1})|$$
$$= |\mathbf{Tr}(\mathbf{S}_{t-1}^\top \mathbf{W}^\top x_t x_t^\top \mathbf{W} \mathbf{S}_{t-1}) + \mathbf{Tr}(\mathbf{S}_{t-1}^\top \mathbf{W}^\top x_t x_t^\top \widetilde{\mathbf{V}}\widetilde{\mathbf{V}}^\top \mathbf{L}_{t-1})| \quad \text{(iii.K.5)}$$
$$\stackrel{①}{\leq} \frac{3}{2} \mathbf{Tr}(\mathbf{S}_{t-1}^\top \mathbf{W}^\top x_t x_t^\top \mathbf{W} \mathbf{S}_{t-1}) + \frac{1}{2} \|x_t^\top \widetilde{\mathbf{V}}\widetilde{\mathbf{V}}^\top \mathbf{L}_{t-1}\|_2^2 \ .$$

Above, inequality ① is because $2\mathbf{Tr}(\mathbf{A}^\top \mathbf{B}) \leq \mathbf{Tr}(\mathbf{A}^\top \mathbf{A}) + \mathbf{Tr}(\mathbf{B}^\top \mathbf{B})$ which is Young's inequality in the matrix case. We take expectation and get:

$$\mathbb{E}\left[|\mathbf{Tr}(\mathbf{S}_{t-1}^\top \mathbf{R}_t')| \mid \mathcal{F}_{t-1}, \mathcal{C}_{\leq t}\right]$$
$$\leq \frac{3}{2}\mathbf{Tr}(\mathbf{S}_{t-1}^\top \mathbf{W}^\top (\mathbf{\Sigma} + \mathbf{\Delta})\mathbf{W}\mathbf{S}_{t-1}) + \frac{1}{2}\mathbf{Tr}(\mathbf{L}_{t-1}^\top \widetilde{\mathbf{V}}\widetilde{\mathbf{V}}^\top (\mathbf{\Sigma} + \mathbf{\Delta})\widetilde{\mathbf{V}}\widetilde{\mathbf{V}}^\top \mathbf{L}_{t-1})$$
$$\leq \frac{3}{2}\lambda_{k+1}\mathbf{Tr}(\mathbf{S}_{t-1}^\top \mathbf{S}_{t-1}) + \frac{1}{2}\mathbf{Tr}(\mathbf{L}_{t-1}^\top \widetilde{\mathbf{V}}\mathbf{\Sigma}_{\leq k+m}\widetilde{\mathbf{V}}^\top \mathbf{L}_t) + \|\mathbf{\Delta}\|_2 \cdot \left(\frac{3}{2}\mathbf{Tr}(\mathbf{S}_{t-1}^\top \mathbf{S}_{t-1}) + \frac{1}{2}\mathbf{Tr}(\mathbf{L}_{t-1}^\top \widetilde{\mathbf{V}}\widetilde{\mathbf{V}}^\top \mathbf{L}_{t-1})\right)$$
$$\overset{①}{\leq} \frac{3}{2}\lambda_{k+1}\mathbf{Tr}(\mathbf{S}_{t-1}^\top \mathbf{S}_{t-1}) + \frac{1}{2}\Lambda + \|\mathbf{\Delta}\|_2 \cdot \left(\frac{3}{2}\mathbf{Tr}(\mathbf{S}_{t-1}^\top \mathbf{S}_{t-1}) + \frac{\chi}{2}\right) \qquad \text{(iii.K.6)}$$

Above, inequality ① has relied on our earlier observations (iii.K.2). At this point, plugging (iii.K.6), (iii.K.3) into (iii.K.4) and using the assumption $\eta_t \phi_t \leq 1/2$, we have

$$\begin{aligned}\mathbb{E}[\mathbf{Tr}(\mathbf{S}_t^\top \mathbf{S}_t) \mid \mathcal{F}_{\leq t-1}, \mathcal{C}_{\leq t}] &\leq \left(1 + 12\Gamma\eta_t^2 \phi_t^2 + \eta_t^2(6\phi_t + 8)\lambda_{k+1}\right)\mathbf{Tr}(\mathbf{S}_{t-1}\mathbf{S}_{t-1}^\top) + \eta_t^2(2\phi_t + 8)\Lambda \\ &\quad + \mathbb{E}\left[-2\eta_t \mathbf{Tr}(\mathbf{S}_{t-1}^\top \mathbf{S}_{t-1}\mathbf{H}_t') + 2\eta_t \mathbf{Tr}(\mathbf{S}_{t-1}^\top \mathbf{R}_t') \mid \mathcal{F}_{\leq t-1}, \mathcal{C}_{\leq t}\right] \\ &\quad + 2\eta_t \|\mathbf{\Delta}\|_2 \cdot \left(\eta_t(4 + \phi_t)\chi + \eta_t(4 + 3\phi_t)\mathbf{Tr}(\mathbf{S}_{t-1}^\top \mathbf{S}_{t-1})\right) \ . \end{aligned}$$

Finally, inequality (i.B.8) (from the proof of Corollary 6.1-(a)) gives an upper bound on the expected value of $-\mathbf{Tr}(\mathbf{S}_{t-1}^\top \mathbf{S}_{t-1}\mathbf{H}_t') + \mathbf{Tr}(\mathbf{S}_{t-1}^\top \mathbf{R}_t')$. Putting it in gives us the desired bound.

(d) This time we compute a slightly different upper bound from (iii.K.5)

$$\begin{aligned}|\mathbf{Tr}(\mathbf{S}_{t-1}^\top \mathbf{R}_t')| &= |\mathbf{Tr}(\mathbf{S}_{t-1}^\top \mathbf{W}^\top x_t x_t^\top \mathbf{L}_{t-1})| = |\mathbf{Tr}(\mathbf{S}_{t-1}^\top \mathbf{W}^\top x_t x_t^\top (\mathbf{V}\mathbf{V}^\top + \mathbf{Z}\mathbf{Z}^\top)\mathbf{L}_{t-1})| \\ &= |\mathbf{Tr}(\mathbf{S}_{t-1}^\top \mathbf{W}^\top x_t x_t^\top \mathbf{V}) + \mathbf{Tr}(\mathbf{S}_{t-1}^\top \mathbf{W}^\top x_t x_t^\top \mathbf{Z}\mathbf{Z}^\top \mathbf{L}_{t-1})| \\ &\leq |\mathbf{Tr}(\mathbf{S}_{t-1}^\top \mathbf{W}^\top x_t x_t^\top \mathbf{Z}\mathbf{Z}^\top \mathbf{L}_{t-1})| + \|\mathbf{S}_{t-1}^\top \mathbf{W}^\top x_t\|_2 \ . \end{aligned}$$

Plugging this into (i.B.1), we obtain

$$\begin{aligned}|\mathbf{Tr}(\mathbf{S}_t^\top \mathbf{S}_t) - \mathbf{Tr}(\mathbf{S}_{t-1}\mathbf{S}_{t-1}^\top)| &\leq 2\eta_t |\mathbf{Tr}(\mathbf{S}_{t-1}^\top \mathbf{S}_{t-1}\mathbf{H}_t')| + 2\eta_t |\mathbf{Tr}(\mathbf{S}_{t-1}^\top \mathbf{R}_t')| + 4\eta_t^2 \|\mathbf{H}_t'\|_2 \left|\mathbf{Tr}(\mathbf{S}_{t-1}^\top \mathbf{R}_t')\right| \\ &\quad + 12\eta_t^2 \|\mathbf{H}_t'\|_2^2 \mathbf{Tr}(\mathbf{S}_{t-1}\mathbf{S}_{t-1}^\top) + 8\eta_t^2 \|\mathbf{R}_t'\|_2^2 \ . \\ &\overset{①}{\leq} 8\eta_t \|\mathbf{H}_t'\|_2 \mathbf{Tr}(\mathbf{S}_{t-1}\mathbf{S}_{t-1}^\top) + 4\eta_t |\mathbf{Tr}(\mathbf{S}_{t-1}^\top \mathbf{R}_t')| + 8\eta_t^2 \|\mathbf{R}_t'\|_2^2 \\ &\leq 8\eta_t \|\mathbf{H}_t'\|_2 \mathbf{Tr}(\mathbf{S}_{t-1}\mathbf{S}_{t-1}^\top) + 4\eta_t |\mathbf{Tr}(\mathbf{S}_{t-1}^\top \mathbf{W}^\top x_t x_t^\top \mathbf{Z}\mathbf{Z}^\top \mathbf{L}_{t-1})| \\ &\quad + 4\eta_t \|\mathbf{S}_{t-1}^\top \mathbf{W}^\top x_t\|_2 + 8\eta_t^2 \|\mathbf{R}_t'\|_2^2 \\ &\overset{②}{\leq} 8\eta_t \|\mathbf{H}_t'\|_2 \mathbf{Tr}(\mathbf{S}_{t-1}\mathbf{S}_{t-1}^\top) + \eta_t(4\phi_t + 8)\mathbf{Tr}(\mathbf{S}_{t-1}^\top \mathbf{W}^\top x_t x_t^\top \mathbf{W}\mathbf{S}_{t-1})^{1/2} \\ &\quad + 8\eta_t^2 \phi_t \|\mathbf{R}_t'\|_2 \end{aligned}$$

Above, ① uses the fact that $\eta_t \|\mathbf{H}_t'\|_2 \leq \eta_t \phi_t \leq 1/2$, and ② uses $\|\mathbf{R}_t'\|_2 \leq \phi_t$ as well as

$$\|x_t^\top \mathbf{Z}\mathbf{Z}^\top \mathbf{L}_{t-1}\|_2 \leq \|x_t^\top \mathbf{L}_{t-1}\|_2 + \|x_t^\top \mathbf{V}\mathbf{V}^\top \mathbf{L}_{t-1}\|_2 \leq (\phi_t + 1)$$

Taking square on both sides, we have

$$\begin{aligned}|\mathbf{Tr}(\mathbf{S}_t^\top \mathbf{S}_t) - \mathbf{Tr}(\mathbf{S}_{t-1}\mathbf{S}_{t-1}^\top)|_2^2 &\leq 192\eta_t^2 \|\mathbf{H}_t'\|_2^2 \mathbf{Tr}(\mathbf{S}_{t-1}\mathbf{S}_{t-1}^\top)^2 + 3\eta_t^2(4\phi_t + 8)^2 \mathbf{Tr}(\mathbf{S}_{t-1}^\top \mathbf{W}^\top x_t x_t^\top \mathbf{W}\mathbf{S}_{t-1}) \\ &\quad + 192\eta_t^4 \phi_t^2 \|\mathbf{R}_t'\|_2^2\end{aligned}$$

Finally, taking expectation and using (iii.K.3), we have (noticing that $\eta_t \phi_t \leq 1/2$)

$$\mathbb{E}[|\mathbf{Tr}(\mathbf{S}_t^\top \mathbf{S}_t) - \mathbf{Tr}(\mathbf{S}_{t-1}\mathbf{S}_{t-1}^\top)|_2^2 \mid \mathcal{F}_{\leq t-1}, \mathcal{C}_{\leq t}]$$
$$\leq 192 \Gamma \eta_t^2 \phi_t^2 \mathbf{Tr}(\mathbf{S}_{t-1}\mathbf{S}_{t-1}^\top)^2 + 3\eta_t^2 (4\phi_t + 8)^2 \mathbf{Tr}(\mathbf{S}_{t-1}^\top \mathbf{W}^\top (\mathbf{\Sigma} + \mathbf{\Delta})\mathbf{W}\mathbf{S}_{t-1})$$
$$+ 192\eta_t^4 \phi_t^2 \left( \Lambda + \lambda_{k+1}\mathbf{Tr}(\mathbf{S}_{t-1}^\top \mathbf{S}_{t-1}) + \|\mathbf{\Delta}\|_2 \cdot (\chi + \mathbf{Tr}(\mathbf{S}_{t-1}^\top \mathbf{S}_{t-1})) \ . \right)$$
$$\leq 192 \Gamma \eta_t^2 \phi_t^2 \mathbf{Tr}(\mathbf{S}_{t-1}\mathbf{S}_{t-1}^\top)^2 + 3\eta_t^2 (4\phi_t + 8)^2 \left( (\lambda_{k+1} + \|\mathbf{\Delta}\|_2) \mathbf{Tr}(\mathbf{S}_{t-1}^\top \mathbf{S}_{t-1}) \right)$$
$$+ 192\eta_t^4 \phi_t^2 \left( \Lambda + \lambda_{k+1}\mathbf{Tr}(\mathbf{S}_{t-1}^\top \mathbf{S}_{t-1}) + \|\mathbf{\Delta}\|_2 \cdot (\chi + \mathbf{Tr}(\mathbf{S}_{t-1}^\top \mathbf{S}_{t-1})) \ . \right)$$
$$\leq 192 \Gamma \eta_t^2 \phi_t^2 \mathbf{Tr}(\mathbf{S}_{t-1}\mathbf{S}_{t-1}^\top)^2 + 4\eta_t^2 (4\phi_t + 8)^2 \lambda_{k+1} \mathbf{Tr}(\mathbf{S}_{t-1}\mathbf{S}_{t-1}^\top)$$
$$+ 192 \Gamma \eta_t^4 \phi_t^2 + \|\mathbf{\Delta}\|_2 \cdot 4\eta_t^2 (4\phi_t + 8)^2 (\chi + \mathbf{Tr}(\mathbf{S}_{t-1}^\top \mathbf{S}_{t-1})) \ .$$

□

## iii.L Final Main Lemmas

In this section, we extend our main lemmas in Section 7 to their strong forms. Specifically,

- Lemma Main 4 is an extension of Lemma Main 1;
- Lemma Main 5 is an extension of Lemma Main 2;
- Lemma Main 6 is a new main lemma that takes into account "under-sampling".

Recall that given parameter $\rho \in (0, \lambda_k)$,

- $\mathbf{V}$ is a matrix consisting of all the eigenvectors of $\mathbf{\Sigma}$ with eigenvalue $\geq \lambda_k$,
- $\mathbf{Z}$ is a matrix consisting of all the eigenvectors of $\mathbf{\Sigma}$ with eigenvalue $< \lambda_k$,
- $\mathbf{W}$ is a matrix consisting of all the eigenvectors of $\mathbf{\Sigma}$ with eigenvalue $\leq \lambda_k - \rho$.

When we apply these main lemmas in later sections, we may redefine the meaning of $(\mathbf{V}, \mathbf{Z}, \mathbf{W})$ to be with respect to some other $\lambda$ that is not necessarily $\lambda_k$.

### iii.L.1 Before Warm Start

**Lemma Main 4** (before warm start). *For every $q \in \left(0, \frac{1}{2}\right]$, $\Xi_{\mathbf{Z}} \geq 2$, $\Xi_x \geq 2$, and fixed matrix $\mathbf{Q} \in \mathbb{R}^{d \times k}$, suppose the initial matrix $\mathbf{Q}$ satisfies*

- $\|\mathbf{Z}^\top \mathbf{Q}(\mathbf{V}^\top \mathbf{Q})^{-1}\|_F^2 \leq \Xi_{\mathbf{Z}}$,
- $\mathbf{Pr}_{x_t}\left[\forall j \in [T], \left\|x_t^\top \mathbf{Z}\mathbf{Z}^\top (\mathbf{\Sigma}/\lambda_{k+1})^{j-1} \mathbf{Q}(\mathbf{V}^\top \mathbf{Q})^{-1}\right\|_2 \leq \Xi_x\right] \geq 1 - q^2/2$ *for every $t \in [T]$,*
- $\left\|\nu_j^\top \mathbf{Z}\mathbf{Z}^\top \mathbf{Q}(\mathbf{V}^\top \mathbf{Q})^{-1}\right\|_2 \leq \Xi_x$ *for every $j \in [d]$.*[24]

---
[24] This assumption is redundant for $j \in [k]$ because $\nu_j^\top \mathbf{Z} = 0$ for $j \in [k]$.

Let $m$ be the number of eigenvalues of $\boldsymbol{\Sigma}$ in $(\lambda_k - \rho, \lambda_{k+1}]$. Let $\Lambda = \sum_{i=1}^{k} \lambda_i + \Xi_x^2 \sum_{j=k+1}^{k+m} \lambda_j$. Suppose also the learning rates $\{\eta_s\}_{s \in [T]}$ satisfy

$$(1): \; \forall s \in [T], \frac{2q(\Xi_{\mathbf{Z}}^{3/2} + d\Xi_x^2)}{\Lambda} \leq \eta_s \leq O\Big(\frac{\rho}{\Lambda \cdot \Xi_x^2 \ln \frac{dT}{q}}\Big) \quad (2): \; \sum_{t=1}^{T} \Lambda \eta_t^2 \Xi_x^2 \leq O\Big(\frac{1}{\ln \frac{dT}{q}}\Big) \;.$$

$$(3): \; \exists T_0 \in [T] \text{ such that } \sum_{t=1}^{T_0} \eta_t \geq \Omega\Big(\frac{\ln(\Xi_{\mathbf{Z}})}{\rho}\Big)$$

(iii.L.1)

Then, for every $t \in [T-1]$, we have with probability at least $1 - 2qT$ (over the randomness of $x_1, \ldots, x_t$):

- if $t \geq T_0$ then $\|\mathbf{W}^\top \mathbf{P}_t \mathbf{Q} (\mathbf{V}^\top \mathbf{P}_t \mathbf{Q})^{-1}\|_F^2 \leq 2$.

*Proof of Lemma Main 4.* The proof is a non-trivial adaption of the proof of Lemma Main 1.

This time we consider random vectors $y_{t,s} \in \mathbb{R}^{2+d+T}$ defined as:

$$y_{t,s}^{(1)} \stackrel{\text{def}}{=} \|\mathbf{Z}^\top \mathbf{P}_s \mathbf{Q} (\mathbf{V}^\top \mathbf{P}_s \mathbf{Q})^{-1}\|_F^2 \;,$$

$$y_{t,s}^{(2)} \stackrel{\text{def}}{=} \|\mathbf{W}^\top \mathbf{P}_s \mathbf{Q} (\mathbf{V}^\top \mathbf{P}_s \mathbf{Q})^{-1}\|_F^2 \;,$$

$$y_{t,s}^{(2+T+j)} \stackrel{\text{def}}{=} \left\|\nu_j^\top \mathbf{Z}\mathbf{Z}^\top \mathbf{P}_s \mathbf{Q} (\mathbf{V}^\top \mathbf{P}_s \mathbf{Q})^{-1}\right\|_2^2, \quad \text{for } j \in [d] \;,$$

$$y_{t,s}^{(3+d+j)} \stackrel{\text{def}}{=} \begin{cases} \left\|x_t^\top \mathbf{Z}\mathbf{Z}^\top (\boldsymbol{\Sigma}/\lambda_{k+1})^j \mathbf{P}_s \mathbf{Q} (\mathbf{V}^\top \mathbf{P}_s \mathbf{Q})^{-1}\right\|_2^2, & \text{for } j \in \{0,1,\ldots,t-s-1\}; \\ (1 - \eta_s \lambda_k) \cdot y_{t,s-1}^{(3+d+j)}, & \text{for } j \in \{t-s,\ldots,T-1\}. \end{cases}$$

We again consider upper bounds

$$\phi_{t,s}^{(1)} \stackrel{\text{def}}{=} 2\Xi_{\mathbf{Z}}, \quad \phi_{t,s}^{(2)} \stackrel{\text{def}}{=} \begin{cases} 2\Xi_{\mathbf{Z}} & s < T_0; \\ 2 & \text{otherwise.} \end{cases}, \quad \text{and } \phi_{t,s}^{(3)} = \cdots = \phi_{t,s}^{(2+d+T)} \stackrel{\text{def}}{=} 2\Xi_x^2 \;.$$

For each $t \in [T]$, define event $\mathcal{C}_t'$ and $\mathcal{C}_t''$ in the same way as before:

$$\mathcal{C}_t' \stackrel{\text{def}}{=} \left\{(x_1, \ldots, x_{t-1}) \text{ satisfies } \Pr_{x_t}\left[\exists i \in [3+d]: y_{t,t-1}^{(i)} > \phi_{t,t-1}^{(i)} \mid \mathcal{F}_{t-1}\right] \leq q\right\}$$

$$\mathcal{C}_t'' \stackrel{\text{def}}{=} \left\{(x_1, \ldots, x_t) \text{ satisfies } \forall i \in [3+d]: y_{t,t-1}^{(i)} \leq \phi_{t,t-1}^{(i)}\right\}$$

and denote by $\mathcal{C}_t \stackrel{\text{def}}{=} \mathcal{C}_t' \wedge \mathcal{C}_t''$ and $\mathcal{C}_{\leq t} \stackrel{\text{def}}{=} \bigwedge_{s=1}^{t} \mathcal{C}_s$.

**Verification of Assumption (A1) in Lemma i.D.1.**

Suppose $\mathbb{E}[x_s x_s^\top \mid \mathcal{C}_{\leq s}, \mathcal{F}_{\leq s-1}] = \boldsymbol{\Sigma} + \boldsymbol{\Delta}$, then we have $\|\boldsymbol{\Delta}\|_2 \leq \frac{q_1}{1-q_1} \leq \frac{q}{1-q}$ using the same proof as before.

This time, we use Lemma iii.K.1 (instead of Lemma i.B.1) with $\phi_t = 2\Xi_x$ to obtain the following tighter bounds for $i \in [T+d+2]$:

$$\mathbb{E}\left[y_{t,s+1}^{(i)} \mid \mathcal{F}_t, \mathcal{F}_{\leq s}, \mathcal{C}_{\leq s}\right] \leq f_s^{(i)}(y_{t,s}, q) \quad \text{and} \quad \mathbb{E}\left[|y_{t,s+1}^{(i)} - y_{t,s}^{(i)}|^2 \mid \mathcal{F}_t, \mathcal{F}_{\leq s}, \mathcal{C}_{\leq s}\right] \leq h_s^{(i)}(y_{t,s}, q)$$

where we define[25]

$$f_s^{(i)}(y,q) \stackrel{\text{def}}{=} \left(1 + O(\Lambda\eta_{s+1}^2\Xi_x^2)\right)y^{(i)} + O\left(\Lambda\eta_{s+1}^2\Xi_x + Err\right) \qquad \text{for } i = 1, 3, 4, \ldots 2+d$$

$$f_s^{(i)}(y,q) \stackrel{\text{def}}{=} \left(1 - 2\eta_{s+1}\rho + O(\Lambda\eta_{s+1}^2\Xi_x^2)\right)y^{(i)} + O\left(\Lambda\eta_{s+1}^2\Xi_x + Err\right) \qquad \text{for } i = 2$$

$$f_s^{(i)}(y,q) \stackrel{\text{def}}{=} \left(1 - \eta_{s+1}\lambda_k + O(\Lambda\eta_{s+1}^2\Xi_x^2)\right)y^{(i)} + \eta_{s+1}\lambda_k y^{(i+1)} + O\left(\Lambda\eta_{s+1}^2\Xi_x^2 + Err\right)$$
$$\text{for } i = 3+d, \ldots, 2+d+T$$

$$h_s^{(i)}(y,q) \stackrel{\text{def}}{=} O\left(\Lambda\eta_{s+1}^2\Xi_x^2(y^{(i)})^2 + \Lambda\eta_{s+1}^2\Xi_x^2 y^{(i)} + \Lambda\eta_{s+1}^4\Xi_x^2 + Err\right) \qquad \text{for } i = 1,2,3,\ldots 2+d$$

$$h_s^{(i)}(y,q) \stackrel{\text{def}}{=} O\left(\Lambda\eta_{s+1}^2\Xi_x^2(y^{(i)})^2 + \Lambda\eta_{s+1}^2\Xi_x^2 y^{(i)} + \Lambda\eta_{s+1}^4\Xi_x^4\right) \qquad \text{for } i = 3+d, \ldots, 2+d+T.$$

Above, we denote by $Err \stackrel{\text{def}}{=} \eta_{s+1}\Xi_x(\Xi_{\mathbf{Z}}^{3/2} + d\Xi_x^2) \cdot \frac{q}{1-q}$ the error term similar to the proof of Lemma Main 1. Obviously if $\frac{2q(\Xi_{\mathbf{Z}}^{3/2}+d\Xi_x^2)}{\Lambda} \le \eta_s$ is satisfied then the $Err$ term can be absorbed into the big-$O$ notation.

For every $i \in [2+d+T]$, we consider the same $g_s$ as defined in the proof of Lemma Main 1:
$$g_s^{(i)}(y) = 20\eta_{s+1}\Xi_x \cdot y^{(i)} + 42\eta_{s+1}^2\Xi_x^2$$
and it satisfies whenever $\mathcal{C}_{\le s+1}$ holds then $|y_{t,s+1}^{(i)} - y_{t,s}^{(i)}| \le g_s^{(i)}(y_{t,s})$.

Putting the above bounds together, we finish verifying assumption (A1) of Lemma i.D.1.

**Verification of Assumption (A2) of Lemma i.D.1.**

This step is exactly the same as the proof of Lemma Main 1 so ignored here.

**Verification of Assumption (A3) of Lemma i.D.1.**

For every $t \in [T]$, at a high level assumption (A3) is satisfied once we plug in the following three sets of parameter choices to Corollary i.C.4 and Corollary i.C.5: define $\kappa \stackrel{\text{def}}{=} 1/\sqrt{\Lambda} > 1$ and for every $s \in [T-1]$,

$$\beta_{s,1} = 0, \qquad \delta_{s,1} = 0, \qquad \tau_{s,1} = O(\eta_{s+1}\Xi_x \cdot \sqrt{\Lambda})$$
$$\beta_{s,2} = 2\eta_{s+1}\rho, \qquad \delta_{s,2} = 0, \qquad \tau_{s,2} = O(\eta_{s+1}\Xi_x \cdot \sqrt{\Lambda})$$
$$\beta_{s,3} = \eta_{s+1}\rho, \qquad \delta_{s,3} = \eta_{s+1}\lambda_k \qquad \tau_{s,3} = O(\eta_{s+1}\Xi_x \cdot \sqrt{\Lambda})$$

More specifically, for every $t \in [T]$, let $\{z_s\}_{s=0}^{t-1}$ be the *arbitrary* random vector satisfying (i.D.1) of Lemma i.D.1. Define $q_2 = q^2/(8+2d)$.

- For coordinate $i = 1, 3, 4, \ldots, 2+d$ of $\{z_s^{(i)}\}_{s=0}^{t-1}$, apply Corollary i.C.4 with $\{\beta_{s,1}, \delta_{s,1}, \tau_{s,1}\}_{s=0}^{t-2}$, $q = q_2$, $D = 1$, and $\kappa$;

- For coordinate $i = 2$ of $\{z_s\}_{s=0}^{t-1}$,
  - if $t < T_0$, apply Corollary i.C.4 with $\{\beta_{s,2}, \delta_{s,2}, \tau_{s,2}\}_{s=0}^{t-2}$, $q = q_2$, $D = 1$, and $\kappa$;
  - if $t \ge T_0$, apply Corollary i.C.5 with $\{\beta_{s,2}, \delta_{s,2}, \tau_{s,2}\}_{s=0}^{t-2}$, $q = q_2$, $D = 1$, $\gamma = 1$, and $\kappa$;

- For coordinates $i = 2+d+1, \ldots, 2+d+T$ of $\{z_s\}_{s=0}^{t-1}$,
  - apply Corollary i.C.4 with $\{\beta_{s,3}, \delta_{s,3}, \tau_{s,3}\}_{s=0}^{t-2}$, $q = q_2$, $D = T$, and $\kappa$.

---

[25] We refer readers to Footnote 18 in the proof of Lemma Main 1 on page 30 for a detailed discussion as well as a careful treatment for out-of-bound indices. We remark here that to derive such bounds, one needs to use the fact that when $w = x\mathbf{Z}\mathbf{Z}^\top$ for some unit norm vector $x$ (such as $x = x_t$ for some $t \in [T]$ or $x = \nu_j$ for some $j \in [d]$), the quantity $\frac{\eta_t}{\lambda_k}\|w^\top \mathbf{\Sigma L}_{t-1}\|_2^2$ that appeared in Lemma iii.K.1-(a) can be upper bounded by $\frac{\eta_t \lambda_{k+1}^2}{\lambda_k}\|w^\top \mathbf{L}_{t-1}\|_2^2 \le \eta_t \lambda_k \|w^\top \mathbf{L}_{t-1}\|_2^2$.

Note that we can apply Corollary i.C.4 because for every $i = 1, 2, 3$, our assumptions on $\eta_s$ imply $\sum_{s=0}^{T-1} \tau_{s,i}^2 \leq \frac{1}{100} \ln^{-1} \frac{4T}{q_2}$ and $\tau_{s,i} \leq \frac{\sqrt{\Lambda}}{24 \log(4T/q_2)}$.

We can apply Corollary i.C.5 with $\gamma = 1$ because our assumption $\eta_s \leq O\left(\frac{\rho}{\Lambda \cdot \Xi_x^2 \ln \frac{dT}{q}}\right)$ implies $\beta_{s,2} \geq 10 \ln \frac{3(T+d)}{q_2} \cdot \tau_{s,2}^2$ for every $s$, and our assumption $\sum_{s=0}^{T_0-1} \beta_{s,2} \geq 1 + \ln \Xi_\mathbf{Z}$ implies $\sum_{s=0}^{t-1} \beta_s - 10 \ln \frac{3t}{q_2} \tau_s^2 \geq \ln \Xi_\mathbf{Z} + 1 - 1 = \ln \Xi_\mathbf{Z}$ whenever $t > T_0$.

Therefore, the conclusion of Corollary i.C.4 and Corollary i.C.5 imply that
$$\mathbf{Pr}[\exists i \in [3+d] : z_{t-1}^{(i)} > \phi_{t,t-1}^{(i)}] \leq (3+d)q_2 < q^2/2$$
so assumption (A3) of Lemma i.D.1 holds.

**Application of Lemma i.D.1.** Applying Lemma i.D.1, we have $\mathbf{Pr}[\overline{\mathcal{C}_T}] \leq 2qT$ which implies our desired bounds and this finishes the proof of Lemma Main 4. □

### iii.L.2 After Warm Start

**Lemma Main 5** (after warm start). *In the same setting as Lemma Main 4, suppose in addition there exists $\delta \leq 1/\sqrt{8}$ such that*
$$\frac{T_0}{\ln^2 T_0} \geq \frac{9 \ln((8+2d)/q^2)}{\delta^2} \ , \quad \forall s \in \{T_0+1, \ldots, T\}: \quad 2\eta_s \rho - \eta_s^2 \Xi_x^2 \geq \frac{\Omega(1)}{s-1} \quad \text{and} \quad \eta_s \leq \frac{O(1)}{\sqrt{\Lambda(s-1)}\delta \Xi_x} \ .$$
*Then, with probability at least $1 - 2qT$ (over the randomness of $x_1, \ldots, x_T$):*

- $\|\mathbf{W}^\top \mathbf{P}_t \mathbf{Q}(\mathbf{V}^\top \mathbf{P}_t \mathbf{Q})^{-1}\|_F^2 \leq \frac{5T_0/\ln^2(T_0)}{t/\ln^2 t}$ *for every $t \in \{T_0, \ldots, T\}$.*

*Proof of Lemma Main 5.* For every $t \in [T]$ and $s \in \{0, 1, \ldots, t-1\}$, consider *the same* random vectors $y_{t,s} \in \mathbb{R}^{2+d+T}$ defined in the proof of Lemma Main 4. This time, we define upper bounds:
$$\phi_{t,s}^{(1)} \stackrel{\text{def}}{=} 2\Xi_\mathbf{Z}, \quad \phi_{t,s}^{(2)} \stackrel{\text{def}}{=} \begin{cases} 2\Xi_\mathbf{Z} & \text{if } s < T_0; \\ 2 & \text{if } s = T_0; \\ \frac{5T_0/\ln^2(T_0)}{s/\ln^2 s} & \text{if } s > T_0. \end{cases}, \quad \text{and } \phi_{t,s}^{(3)} = \cdots = \phi_{t,s}^{(2+d+T)} \stackrel{\text{def}}{=} 2\Xi_x^2 \ .$$

Also consider the same events $\mathcal{C}_t'$, $\mathcal{C}_t''$, $\mathcal{C}_t \stackrel{\text{def}}{=} \mathcal{C}_t' \wedge \mathcal{C}_t''$ and $\mathcal{C}_{\leq t} \stackrel{\text{def}}{=} \bigwedge_{s=1}^t \mathcal{C}_s$ defined in the proof of Lemma Main 4. We again want to apply the decoupling Lemma i.D.1.

**Verification of Assumption (A1) in Lemma i.D.1.**

The same functions $f_s^{(i)}$, $g_s^{(i)}$, and $h_s^{(i)}$ used in the proof of Lemma Main 4 still apply here. We make minor changes on the second coordinate (and this similar modification was also done in Lemma Main 2): whenever $s \geq T_0$, define
$$g_s^{(2)}(y) \stackrel{\text{def}}{=} 45\eta_{s+1}\Xi_x\sqrt{y^{(2)}} + 40\eta_{s+1}^2\Xi_x^2 \quad \text{and} \quad h_s^{(2)}(y,q) \stackrel{\text{def}}{=} O\left(\Lambda \eta_{s+1}^2 \Xi_x^2 y^{(2)} + \Lambda \eta_{s+1}^4 \Xi_x^2 + Err\right) \ .$$

Note that we can make this change for $g_s^{(2)}$ owing to exactly the same reason as the proof of Lemma Main 2. We can do so for $h_s^{(2)}$ because whenever $\mathcal{C}_{\leq s+1}$ holds for some $s \geq T_0$ (which implies $y_{t,s}^{(2)} \leq 5$), we have $(y^{(2)})^2 = O(y^{(2)})$ so the formulation of $h_s^{(2)}$ can be simplified as above.

These choices of $f_s^{(i)}$, $g_s^{(i)}$, and $h_s^{(i)}$ satisfy assumption (A1) of Lemma i.D.1.

**Verification of Assumption (A2) of Lemma i.D.1.**

Same as before.

**Verification of Assumption (A3) in Lemma i.D.1.**

Same as the proof of Lemma Main 2, for every $t \in [T]$, let $\{z_s\}_{s=0}^{t-1}$ be the *arbitrary* random vector satisfying (i.D.1) of Lemma i.D.1. Choosing $q_2 = q^2/(8 + 2d)$ again, the same argument before indicates that it suffices to focus on $t \geq T_0 + 2$ and prove

$$\mathbf{Pr}[z_{t-1}^{(2)} > \phi_{t,t-1}^{(2)} \mid z_{T_0}^{(2)} \leq 2] \leq q_2 \ . \tag{iii.L.2}$$

We want to apply Corollary i.C.3. Recall that for every $t \in \{T_0+2, \ldots, T\}$, the random sequence $\{z_s^{(2)}\}_{s=T_0}^{t-1}$ satisfies (i.D.1) with

$$f_s^{(2)}(y, q) \stackrel{\text{def}}{=} (1 - 2\eta_{s+1}\rho + O(\Lambda\eta_{s+1}^2\Xi_x^2))y^{(2)} + O\left(\Lambda\eta_{s+1}^2\Xi_x + Err\right) \ ,$$

$$h_s^{(2)}(y, q) \stackrel{\text{def}}{=} O\left(\Lambda\eta_{s+1}^2\Xi_x^2 y^{(2)} + \Lambda\eta_{s+1}^4\Xi_x^2 + Err\right) \ ,$$

$$g_s^{(2)}(y) \stackrel{\text{def}}{=} 45\eta_{s+1}\Xi_x\sqrt{y^{(2)}} + 40\eta_{s+1}^2\Xi_x^2$$

Therefore, $\{z_s^{(2)}\}_{s=T_0}^{t-1}$ satisfies (i.C.1) with $\kappa = 2/\sqrt{\Lambda}$ and $\tau_s = \frac{1}{\delta s}$ because the following holds from our assumptions:

$$q\Xi_{\mathbf{Z}}^{3/2} \leq \eta_{s+1} \qquad \delta\tau_s = \frac{1}{s} \leq 2\eta_{s+1}\rho - \Omega(\eta_{s+1}^2\Xi_x^2)$$

$$\tau_s^2 = \frac{1}{\delta^2 s^2} \geq \Omega(\Lambda\eta_{s+1}^2\Xi_x^2) \qquad \kappa^2\tau_s^4 = \frac{1}{\Lambda\delta^4 s^4} \geq \Omega(\Lambda\eta_{s+1}^4\Xi_x^2) \qquad \kappa\tau_s = \frac{2}{\sqrt{\Lambda}\delta s} \geq \Omega(\eta_{s+1}\Xi_x)$$

Finally, we are ready to apply Corollary i.C.3 with $q = q_2$, $t_0 = T_0$, and $\kappa = 2/\sqrt{\Lambda}$. Because $q_2 \leq e^{-2}$, $z_{T_0}^{(2)} \leq 2$, $\delta \leq 1/\sqrt{8}$ and $\frac{T_0}{\ln^2 T_0} \geq \frac{9\ln(1/q_2)}{\delta^2}$, the conclusion of Corollary i.C.3 tells us $\mathbf{Pr}[z_{t-1}^{(2)} > \phi_{t,t-1}^{(2)} \mid z_{T_0}^{(2)} \leq 2] \leq q_2$, which is exactly (iii.L.2) so this finishes the verification of assumption (A3).

**Application of Lemma i.D.1.** Applying Lemma i.D.1, we have $\mathbf{Pr}[\overline{\mathcal{C}_T}] \leq 2qT$ which implies our desired bounds and this finishes the proof of Lemma Main 5. $\square$

**Parameter iii.L.1.** There exists constants $C_1, C_2, C_3 > 0$ such that for every $q > 0$ that is sufficiently small (meaning $q < 1/\mathsf{poly}(T, \Xi_{\mathbf{Z}}, \Xi_x, 1/\mathsf{gap})$), the following parameters both satisfy Lemma Main 4 and Lemma Main 5:

$$\frac{T_0}{\ln^2(T_0)} = C_1 \cdot \frac{\Lambda\Xi_x^2 \ln \frac{dT}{q} \ln^2 \Xi_{\mathbf{Z}}}{\rho^2} \ , \quad \eta_t = C_2 \cdot \begin{cases} \frac{\ln \Xi_{\mathbf{Z}}}{T_0 \cdot \rho} & t \leq T_0; \\ \frac{1}{t \cdot \rho} & t > T_0. \end{cases} \ , \quad \text{and} \quad \delta = C_3 \cdot \frac{\rho}{\sqrt{\Lambda}\Xi_x} \ .$$

### iii.L.3 Under-Sampling Lemma

The previous two sections together (namely, Lemma Main 4 and Lemma Main 5 together), analyze the behavior of a rank-$k$ Oja's algorithm with an eigen-partition $(\lambda_k, \lambda_k - \rho)$, meaning $\mathbf{V}$ consists of all eigenvectors with eigenvalues $\geq \lambda_k$, $\mathbf{Z}$ consists of all the eigenvectors with eigenvalues $< \lambda_k$, and $\mathbf{W}$ consists of all the eigenvectors with eigenvalues $\leq \lambda_k - \rho$.

In this section, we consider a more general scenario

**Definition iii.L.2.** *Given $\lambda \in [0, \lambda_k]$ and $\rho \in (0, \lambda)$, we say that $\mathbf{V}, \mathbf{Z}, \mathbf{W}$ form an $(\lambda, \lambda - \rho)$ eigen-partition if $\mathbf{V}$ consists of all eigenvectors with eigenvalues $\geq \lambda$, $\mathbf{Z}$ consist of all the eigenvectors with eigenvalues $< \lambda$, and $\mathbf{W}$ consist of all the eigenvectors with eigenvalues $\leq \lambda - \rho$.*

The following lemma studies the behavior of a rank-$k'$ Oja's algorithm for an arbitrary $k' \leq r$.

**Lemma Main 6** (under sampling). *Let $\mathbf{V}, \mathbf{Z}, \mathbf{W}$ be an $(\lambda, \lambda - \rho)$ eigen-partition where the $\mathbf{V} \in \mathbb{R}^{d \times r}, \mathbf{Z} \in \mathbb{R}^{d \times (d-r)}, \mathbf{W} \in \mathbb{R}^{d \times (d-m-r)}$. We study rank $k'$ Oja's algorithm for some $k' \leq r$.*

*For every $q \in \left(0, \frac{1}{2}\right]$, $\Xi_{\mathbf{Z}} \geq 2$, $\Xi_x \geq 2$, and fixed matrix $\mathbf{Q} \in \mathbb{R}^{d \times k'}$, suppose it satisfies*

- $\|\mathbf{Z}^\top \mathbf{Q}(\mathbf{Q}^\top \mathbf{V}\mathbf{V}^\top \mathbf{Q})^{-1/2}\|_F^2 \leq \Xi_{\mathbf{Z}}$,
- $\mathbf{Pr}_{x_t}\left[\forall j \in [T], \left\|x_t^\top \mathbf{Z}\mathbf{Z}^\top (\mathbf{\Sigma}/\lambda_{r+1})^{j-1} \mathbf{Q}(\mathbf{Q}^\top \mathbf{V}\mathbf{V}^\top \mathbf{Q})^{-1/2}\right\|_2 \leq \Xi_x\right] \geq 1 - q^2/2$ for every $t \in [T]$,
- $\left\|\nu_j^\top \mathbf{Z}\mathbf{Z}^\top \mathbf{Q}(\mathbf{Q}^\top \mathbf{V}\mathbf{V}^\top \mathbf{Q})^{-1/2}\right\|_2 \leq \Xi_x$ for every $j \in [d]$.

*Suppose also the learning rates $\{\eta_s\}_{s \in [T]}$ satisfy the all conditions in Lemma Main 4 and Lemma Main 5 with $\Lambda$ replaced by $\Lambda' \stackrel{\text{def}}{=} \sum_{i=1}^{r} \lambda_i + \Xi_x^2 \sum_{j=r+1}^{r+m} \lambda_j$.*

*Then, letting $\mathbf{Q}_t \in \mathbb{R}^{d \times k'}$ be the output of the rank-$k'$ Oja's algorithm with respect to input $\mathbf{Q}$, with probability at least $1 - 2qT$ (over the randomness of $x_1, \ldots, x_T$), it satisfies*

- $\|\mathbf{W}^\top \mathbf{Q}_t\|_F^2 \leq \frac{5T_0/\ln^2(T_0)}{t/\ln^2 t}$ for every $t \in \{T_0, \ldots, T\}$.

*Proof of Lemma Main 6.* Since $\mathbf{V}^\top \mathbf{Q} \in \mathbb{R}^{r \times k'}$ for $r \geq k'$, we can always find a (column) orthonormal matrix $\mathbf{S} \in \mathbb{R}^{r \times (r-k')}$ such that $\mathbf{S}^\top \mathbf{V}^\top \mathbf{Q} = 0$. Letting $\widetilde{\mathbf{Q}} = \mathbf{V}\mathbf{S} \in \mathbb{R}^{d \times (r-k')}$, we have
$$\widetilde{\mathbf{Q}}^\top \mathbf{V}\mathbf{V}^\top \widetilde{\mathbf{Q}} = \mathbf{I}, \quad \widetilde{\mathbf{Q}}^\top \mathbf{V}\mathbf{V}^\top \mathbf{Q} = 0, \quad \mathbf{Z}^\top \widetilde{\mathbf{Q}} = 0, \quad \mathbf{Z}^\top \mathbf{\Sigma}^{j-1} \widetilde{\mathbf{Q}} = 0 \ .$$
Therefore, for every $\mathbf{X}^\top = \mathbf{Z}^\top$, $\mathbf{X}^\top = x_t^\top \mathbf{Z}\mathbf{Z}^\top(\mathbf{\Sigma}/\lambda_{r+1})^{j-1}$, or $\mathbf{X}^\top = \nu_j \mathbf{Z}\mathbf{Z}^\top$, we always have
$$\begin{aligned}
\|\mathbf{X}^\top[\mathbf{Q}, \widetilde{\mathbf{Q}}](\mathbf{V}^\top[\mathbf{Q}, \widetilde{\mathbf{Q}}])^{-1}\|_F^2 &= \mathbf{Tr}\left(\mathbf{X}^\top[\mathbf{Q}, \widetilde{\mathbf{Q}}]([\mathbf{Q}, \widetilde{\mathbf{Q}}]^\top \mathbf{V}\mathbf{V}^\top[\mathbf{Q}, \widetilde{\mathbf{Q}}])^{-1}[\mathbf{Q}, \widetilde{\mathbf{Q}}]^\top \mathbf{X}\right) \\
&= \mathbf{Tr}\left(\mathbf{X}^\top \mathbf{Q}(\mathbf{Q}^\top \mathbf{V}\mathbf{V}^\top \mathbf{Q})^{-1}\mathbf{Q}^\top \mathbf{X}\right) \\
&= \|\mathbf{X}^\top \mathbf{Q}(\mathbf{Q}^\top \mathbf{V}\mathbf{V}^\top \mathbf{Q})^{-1/2}\|_F^2 \leq \Xi_{\mathbf{X}}^2 \ .
\end{aligned}$$
This implies this matrix $[\mathbf{Q}, \widetilde{\mathbf{Q}}] \in \mathbb{R}^{d \times r}$ now satisfies

- $\|\mathbf{Z}^\top[\mathbf{Q}, \widetilde{\mathbf{Q}}](\mathbf{V}^\top[\mathbf{Q}, \widetilde{\mathbf{Q}}])^{-1}\|_F^2 \leq \Xi_{\mathbf{Z}}$,
- $\mathbf{Pr}_{x_t}\left[\forall j \in [T], \left\|x_t^\top \mathbf{Z}\mathbf{Z}^\top (\mathbf{\Sigma}/\lambda_{r+1})^{j-1} [\mathbf{Q}, \widetilde{\mathbf{Q}}](\mathbf{V}^\top[\mathbf{Q}, \widetilde{\mathbf{Q}}])^{-1}\right\|_2 \leq \Xi_x\right] \geq 1 - q^2/2$ for all $t \in [T]$,
- $\left\|\nu_j^\top \mathbf{Z}\mathbf{Z}^\top[\mathbf{Q}, \widetilde{\mathbf{Q}}](\mathbf{V}^\top[\mathbf{Q}, \widetilde{\mathbf{Q}}])^{-1}\right\|_2 \leq \Xi_x$ for every $j \in [d]$.

Therefore, we can apply Lemma Main 5 with initial matrix $[\mathbf{Q}, \widetilde{\mathbf{Q}}]$ and $k = r$ and conclude

$$\text{with probability } 1 - 2qT: \quad \|\mathbf{W}^\top \mathbf{P}_t[\mathbf{Q}, \widetilde{\mathbf{Q}}](\mathbf{V}^\top \mathbf{P}_t[\mathbf{Q}, \widetilde{\mathbf{Q}}])^{-1}\|_F^2 \leq \frac{5T_0/\ln^2(T_0)}{t/\ln^2 t} \ .$$

Now, let $\overline{\mathbf{Q}}_t = \mathsf{QR}(\mathbf{P}_t[\mathbf{Q}, \widetilde{\mathbf{Q}}]) \in \mathbb{R}^{d \times r}$ be the output of the rank-$r$ Oja's algorithm on input $[\mathbf{Q}, \widetilde{\mathbf{Q}}]$ after $t$ steps. Owing to Lemma 2.2, we have
$$\|\mathbf{W}^\top \overline{\mathbf{Q}}_t\|_F^2 \leq \|\mathbf{W}^\top \mathbf{P}_t[\mathbf{Q}, \widetilde{\mathbf{Q}}](\mathbf{V}^\top \mathbf{P}_t[\mathbf{Q}, \widetilde{\mathbf{Q}}])^{-1}\|_F^2 \leq \frac{5T_0/\ln^2(T_0)}{t/\ln^2 t} \ .$$

However, recall that $\mathsf{QR}$ decomposition orthonormalizes a column vectors only with respect to its previous columns. This implies, $\mathbf{Q}_t$, which is the output of the rank-$k'$ Oja's algorithm on input $\mathbf{Q}$, is *exactly* identical to the first $k'$ columns of $\overline{\mathbf{Q}}_t$. Therefore, we have
$$\|\mathbf{W}^\top \mathbf{Q}_t\|_F^2 \leq \|\mathbf{W}^\top \overline{\mathbf{Q}}_t\|_F^2 \leq \frac{5T_0/\ln^2(T_0)}{t/\ln^2 t} \ . \qquad \square$$

## iii.M  Proof of Theorems 1 and 2 (for Oja)

We now prove the final Theorem of Oja's algorithm in gap-free case:

> **Theorem 2** (restated). *For every $\rho, \varepsilon, p \in (0,1)$, let $k+m$ be the number of eigenvalues of $\boldsymbol{\Sigma}$ in $(\lambda_k - \rho, 1]$. Let*
> $$\Lambda_1 = \min\left\{\sum_{i=1}^k \lambda_i + \frac{k}{p^2}\sum_{j=k+1}^{k+m} \lambda_j\,,\; 1\right\}, \quad \Lambda_2 = \sum_{i=1}^{k+m} \lambda_i \ .$$
> *Define learning rates*
> $$T_0 = \widetilde{\Theta}\left(\frac{k\Lambda_1}{\rho^2 p^2}\right), \quad T_1 = \widetilde{\Theta}\left(\frac{\Lambda_2}{\rho^2}\right), \quad \eta_t = \begin{cases} \widetilde{\Theta}\left(\frac{1}{\rho \cdot T_0}\right) & t \leq T_0; \\ \widetilde{\Theta}\left(\frac{1}{\rho \cdot T_1}\right) & t \in (T_0, T_0+T_1]; \\ \widetilde{\Theta}\left(\frac{1}{\rho \cdot (t - T_0)}\right) & t > T_0 + T_1. \end{cases}$$
> *Let $\mathbf{W}$ be the column orthonormal matrix consisting of all eigenvectors of $\boldsymbol{\Sigma}$ with values no more than $\lambda_k - \rho$. Then, the output $\mathbf{Q}_T \in \mathbb{R}^{d \times k}$ of Oja's algorithm satisfies with prob. at least $1 - p$:*
> $$\text{for every} \quad T = T_0 + T_1 + \widetilde{\Theta}\left(\frac{T_1}{\varepsilon}\right) \quad \text{it satisfies} \quad \|\mathbf{W}^\top \mathbf{Q}_T\|_F^2 \leq \varepsilon \ .$$
> *Above, $\widetilde{\Theta}$ hides poly-log factors in $\frac{1}{p}, \frac{1}{\rho}$ and $d$.*

Note that Theorem 1 is a direct corollary of Theorem 2 by setting $\rho \leftarrow \mathsf{gap}$ and noticing that $m = 0$ so $\Lambda_1 = \lambda_1 + \cdots + \lambda_k$ and $\Lambda_2 = 0$.

*Proof of Theorem 2.* First for a sufficiently large constant $C$, we can apply Lemma 5.1 with $p' = \frac{p}{8}$ and some sufficiently small $q = \mathsf{poly}(1/T, 1/d, p)$, with probability at least $1 - p' - q^2 \geq 1 - p/4$ over the random choice of $\mathbf{Q}$, the following holds:

$$\begin{cases} \left\|(\mathbf{Z}^\top \mathbf{Q})(\mathbf{V}^\top \mathbf{Q})^{-1}\right\|_F^2 \leq O\!\left(\frac{dk}{p^2} \ln \frac{d}{p}\right) \ , \text{ and} \\ \mathbf{Pr}_{x_1,\ldots,x_T}\!\left[\exists i \in [T], \exists t \in [T], \left\|x_t^\top \mathbf{Z}\mathbf{Z}^\top \left(\boldsymbol{\Sigma}/\lambda_{k+1}\right)^{i-1} \mathbf{Q}(\mathbf{V}^\top \mathbf{Q})^{-1}\right\|_2 \geq \Omega\!\left(\frac{\sqrt{k \ln \frac{T}{q}}}{p}\right)\right] \leq \frac{q^2}{2} \\ \left\|\nu_j^\top \mathbf{Z}\mathbf{Z}^\top \mathbf{Q}(\mathbf{V}^\top \mathbf{Q})^{-1}\right\|_2 \leq O\!\left(\frac{\sqrt{k \ln \frac{T}{q}}}{p}\right) \text{ for every } j \in [d] \ . \end{cases}$$

Denote by $\mathcal{C}_1$ the union of the above three events, and we have $\mathbf{Pr}_\mathbf{Q}[\mathcal{C}_1] \geq 1 - p/4$.

Let us introduce now some notations for this proof. As illustrated in Figure 2, besides the definitions of $\mathbf{V}, \mathbf{W}, \mathbf{Z}$, let $\mathbf{W}_1$ be the (column) orthogonal matrix consisting of all eigenvectors of $\boldsymbol{\Sigma}$ with eigenvalue $\leq \lambda_k - \frac{\rho}{2}$, and $\mathbf{V}_1$ be the (column) orthogonal matrix consisting of all eigenvectors of $\boldsymbol{\Sigma}$ with eigenvalue $> \lambda_k - \frac{\rho}{2}$.

We now wish to apply our main lemma twice. Once with $(\mathbf{V}, \mathbf{Z}, \mathbf{W}) = (\mathbf{V}, \mathbf{Z}, \mathbf{W}_1)$, and once with $(\mathbf{V}, \mathbf{Z}, \mathbf{W}) = (\mathbf{V}_1, \mathbf{W}_1, \mathbf{W})$.

**Application One.** For every fixed $\mathbf{Q}$, whenever $\mathcal{C}_1$ holds, we can let
$$\Xi_\mathbf{Z} = \Theta\!\left(\frac{dk}{p^2} \ln \frac{d}{p}\right), \quad \Xi_x = \Theta\!\left(\frac{\sqrt{k \ln \frac{T}{p}}}{p}\right) \ ,$$
so the initial conditions in Lemma Main 4 is satisfied. Also, according to Parameter iii.L.1, our parameter choices $\eta_t$ for $t \leq T_0$ satisfy the assumptions in Lemma Main 4 with $\Lambda$ replaced by $\Lambda_1$. We can therefore apply Lemma Main 4 with

$$\boxed{\mathbf{Q} = \mathbf{Q}, \quad (\mathbf{V}, \mathbf{Z}, \mathbf{W}) = (\mathbf{V}, \mathbf{Z}, \mathbf{W}_1), \quad \lambda_k = \lambda_k, \quad \rho = \rho/2, \quad \Xi_\mathbf{Z}, \Xi_x}$$

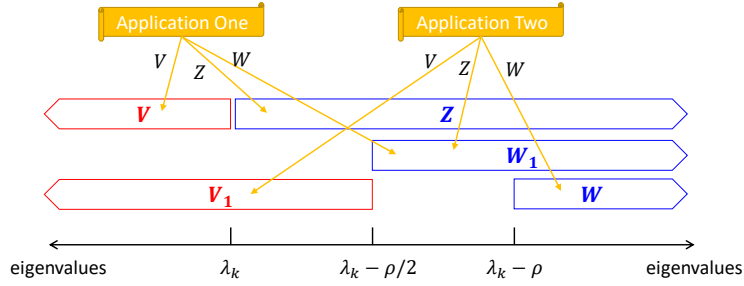

Figure 2: Notations $\mathbf{V}_1, \mathbf{W}_1$ for the proof of Theorem 2

and derive that
$$\Pr_{x_1,\ldots,x_{T_0}}\left[\|\mathbf{W}_1^\top \mathbf{P}_{T_0}\mathbf{Q}(\mathbf{V}^\top \mathbf{P}_{T_0}\mathbf{Q})^{-1}\|_F^2 \geq 2 \,\big|\, \mathcal{C}_1\right] \leq 2qT_0 \leq p/4 \ .$$

Now, let $\mathcal{C}_2$ be the event where $\|\mathbf{W}_1^\top \mathbf{P}_{T_0}\mathbf{Q}(\mathbf{V}^\top \mathbf{P}_{T_0}\mathbf{Q})^{-1}\|_F^2 \leq 2$ holds, and by union bound $\Pr_{\mathbf{Q},x_1,\ldots,x_{T_0}}[\mathcal{C}_2] \geq 1 - p/4 - p/4 = 1 - p/2$.

**Application Two.** If $\mathcal{C}_2$ is true then
$$\|\mathbf{W}_1^\top \mathbf{P}_{T_0}\mathbf{Q}(\mathbf{Q}^\top \mathbf{P}_{T_0}^\top \mathbf{V}_1\mathbf{V}_1^\top \mathbf{P}_{T_0}\mathbf{Q})^{-1/2}\|_F^2 \overset{\text{①}}{\leq} \|\mathbf{W}_1^\top \mathbf{P}_{T_0}\mathbf{Q}(\mathbf{Q}^\top \mathbf{P}_{T_0}^\top \mathbf{V}\mathbf{V}^\top \mathbf{P}_{T_0}\mathbf{Q})^{-1/2}\|_F^2$$
$$\overset{\text{②}}{=} \|\mathbf{W}_1^\top \mathbf{P}_{T_0}\mathbf{Q}(\mathbf{V}^\top \mathbf{P}_{T_0}\mathbf{Q})^{-1}\|_F^2 \leq 2 \ . \qquad \text{(iii.M.1)}$$

Above, inequality ① is because $\mathbf{V}_1\mathbf{V}_1^\top \succeq \mathbf{V}\mathbf{V}^\top$ and this gives us $(\mathbf{Q}^\top \mathbf{P}_{T_0}^\top \mathbf{V}_1\mathbf{V}_1^\top \mathbf{P}_{T_0}\mathbf{Q})^{-1} \preceq (\mathbf{Q}^\top \mathbf{P}_{T_0}^\top \mathbf{V}\mathbf{V}^\top \mathbf{P}_{T_0}\mathbf{Q})^{-1}$; equality ② is because $\mathbf{V}^\top \mathbf{P}_{T_0}\mathbf{Q}$ is a square matrix and we therefore have $(\mathbf{Q}^\top \mathbf{P}_{T_0}^\top \mathbf{V}\mathbf{V}^\top \mathbf{P}_{T_0}\mathbf{Q})^{-1} = (\mathbf{V}^\top \mathbf{P}_{T_0}\mathbf{Q})^{-1}\big((\mathbf{V}^\top \mathbf{P}_{T_0}\mathbf{Q})^{-1}\big)^\top$.

Inequality (iii.M.1) also implies that, for every $x \in \mathbb{R}^d$ with $\|x\|_2 \leq 1$, it satisfies
$$\|x^\top \mathbf{W}_1\mathbf{W}_1^\top \mathbf{P}_{T_0}\mathbf{Q}_0(\mathbf{Q}_0^\top \mathbf{P}_{T_0}^\top \mathbf{V}_1\mathbf{V}_1^\top \mathbf{P}_{T_0}\mathbf{Q}_0)^{-1/2}\|_2 \leq \sqrt{2} \ .$$

In sum, whenever $\mathcal{C}_2$ holds, we have that the initial conditions in Lemma Main 6 is satisfied. Also, according to Parameter iii.L.1, our parameter choices $\eta_t$ for $t = T_0+1,\ldots,T$ —once shifted left by $T_0$— satisfy the assumptions in Lemma Main 6 with $\Lambda$ replaced by $\Lambda_2$. We can therefore apply Lemma Main 6 with

$$\boxed{\mathbf{Q} = \mathbf{P}_{T_0}\mathbf{Q}, \quad (\mathbf{V},\mathbf{Z},\mathbf{W}) = (\mathbf{V}_1, \mathbf{W}_1, \mathbf{W}), \quad \lambda = \lambda_k - \rho/2, \quad \rho = \rho/2, \quad \Xi_x = \Xi_{\mathbf{Z}} = 2}$$

and conclude that
$$\Pr_{x_{T_0+1},\ldots,x_T}\left[\forall T_0 + T_1 \leq t \leq T \colon \|\mathbf{W}^\top \mathbf{Q}_t\|_F^2 \leq \frac{5T_1/\ln^2(T_1)}{(t-T_0)/\ln^2(t-T_0)} \,\bigg|\, \mathcal{C}_2\right] \geq 1 - 2qT \geq 1 - p/4 \ .$$

In sum, we have with probability at least $\Pr[\mathcal{C}_2](1-p/4) \geq 1-p$ over the random choices of $\mathbf{Q}$, and $x_1,\ldots,x_T$, for every $t = T_0 + T_1 + \widetilde{\Theta}(T_1/\varepsilon)$, it satisfies $\|\mathbf{W}^\top \mathbf{Q}_t\|_F^2 \leq \varepsilon$. □

## iii.N  Proof of Theorem 4 and 5 (for Oja$^{++}$)

---

**Algorithm 1** Oja$\big(\{x_t\}_{t=1}^T, \{\eta_t\}_{t=1}^T, \mathbf{Q}\big)$

---

**Input:** vectors $\{x_t \in \mathbb{R}^d\}_{t=1}^T$, learning rates $\{\eta_t \in \mathbb{R}_{>0}\}_{t=1}^T$, an initial matrix $\mathbf{Q} \in \mathbb{R}^{d \times k}$.
1: $\mathbf{Q}_0 \leftarrow \mathbf{Q}$.
2: **for** $t \leftarrow 1$ **to** $T$ **do**
3: $\quad \mathbf{Q}_t \leftarrow \mathsf{QR}\big((\mathbf{I} + \eta_t x_t x_t^\top)\mathbf{Q}_{t-1}\big)$
4: **end for**
5: **return** $\mathbf{Q}_T$

---

**Algorithm 2** Oja$^{++}\big(\{x_t\}_{t=1}^T, \{\eta_t\}_{t=1}^T, \{(T^{(i)}, \mathbf{Q}^{(i)})\}_{i=1}^s\big)$

---

**Input:** vectors $\{x_t \in \mathbb{R}^d\}_{t=1}^T$; learning rates $\{\eta_t \in \mathbb{R}_{>0}\}_{t=1}^T$; initial matrices $\{\mathbf{Q}^{(i)} \in \mathbb{R}^{d \times r_i}\}_{i=1}^s$; lengths $T^{(1)}, \ldots, T^{(s)}$ satisfying $T^{(1)} + \cdots + T^{(s)} = T$.
**Output:** $\mathbf{Q}_T \in \mathbb{R}^{d \times (r_1 + \cdots + r_s)}$.
1: $\mathbf{Q}_0 \leftarrow []; t \leftarrow 0$.
2: **for** $i \leftarrow 1$ **to** $s$ **do**
3: $\quad \mathbf{Q}_t \leftarrow [\mathbf{Q}_t, \mathbf{Q}^{(i)}]$
4: $\quad$ **for** $t' \leftarrow 1$ **to** $T^{(i)}$ **do**
5: $\quad\quad t \leftarrow t + 1$;
6: $\quad\quad \mathbf{Q}_t \leftarrow \mathsf{QR}\big((\mathbf{I} + \eta_t x_t x_t^\top)\mathbf{Q}_{t-1}\big)$
7: $\quad$ **end for**
8: **end for**
9: **return** $\mathbf{Q}_T$

---

We formally write the pseudocode of Oja$^{++}$ in Algorithm 2. We also write down the pseudocode of Oja because we shall use it in the analysis of this section. We emphasize here that Oja$^{++}$ spends the same per-iteration running time and space complexity as Oja.

We prove the following main theorem:

**Theorem 5** (Oja$^{++}$, restated). *For every $\rho, \varepsilon, p \in (0, 1)$, let $k + m$ be the number of eigenvalues of $\Sigma$ in $(\lambda_k - \rho, 1]$. Define*

$$\Lambda_1 = \sum_{i=1}^{k} \lambda_i + \frac{1}{p^2} \sum_{j=k+1}^{m} \lambda_j, \quad \Lambda_2 = \sum_{i=1}^{k+m} \lambda_i, \quad s = \lceil \log(k+1) \rceil, \quad T_0 = \widetilde{\Theta}\left(\frac{\Lambda_1}{\rho^2 p^2}\right), \quad T_1 = \widetilde{\Theta}\left(\frac{\Lambda_2}{\rho^2}\right)$$

*and learning rates (where $C$ is some fixed value that is only $\widetilde{\Theta}(1)$):*

$$\eta_t = \begin{cases} C \cdot \left(\frac{1}{\rho T_0}\right) & \text{if } t \in [11iT_0, 11iT_0 + T_0) \text{ for some } i \in \{0, 1, ..., s-1\}; \\ C \cdot \left(\frac{1}{\rho(t - 11iT_0)}\right) & \text{if } t \in [11iT_0 + T_0, 11iT_0 + 11T_0) \text{ for some } i \in \{0, 1, ..., s-1\}^{26}; \\ C \cdot \left(\frac{1}{\rho T_1}\right) & \text{if } t \in (11sT_0, 11sT_0 + T_1); \\ C \cdot \left(\frac{1}{\rho(t - 11sT_0)}\right) & \text{if } t > 11sT_0 + T_1. \end{cases}$$

*For each $i \in [s]$, let $\mathbf{Q}^{(i)} \in \mathbb{R}^{d \times (\lfloor k/2^{i-1} \rfloor - \lfloor k/2^i \rfloor)}$ be a random matrix with entries i.i.d. $\mathcal{N}(0, 1)$. Then, the output*

$$\mathbf{Q}_T \leftarrow \mathsf{Oja}^{++}\left(\{x_t\}_{t=1}^T, \{\eta_t\}_{t=1}^T, \{(T_0, \mathbf{Q}^{(i)})\}_{i=1}^{s-1} \cup \{(T - 11(s-1)T_0, \mathbf{Q}^{(s)})\}\right)$$

*satisfies with probability at least $1 - p$*

$$\text{for every} \quad T = 11sT_0 + T_1 + \widetilde{\Theta}\left(\frac{T_1}{\varepsilon}\right) \quad \text{it satisfies} \quad \|\mathbf{W}^\top \mathbf{Q}_T\|_F^2 \leq \varepsilon \ .$$

*Above, $\widetilde{\Theta}$ hides poly-log factors in $\frac{1}{p}, \frac{1}{\rho}$ and $d$.*

Note that Theorem 4 is a simple corollary of Theorem 5 by setting $\rho \leftarrow \mathsf{gap}$ and noticing that $m = 0$ so $\Lambda_1 = \lambda_1 + \cdots + \lambda_k$ and $\Lambda_2 = 0$.

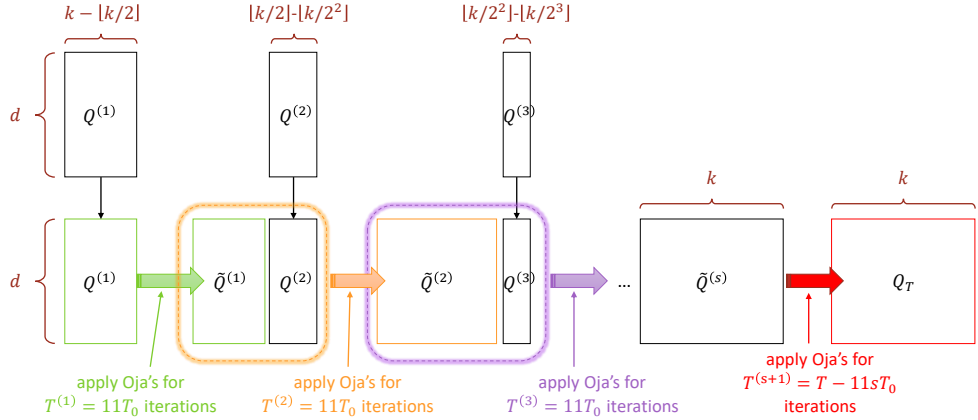

Figure 3: Illustration of our Oja$^{++}$ algorithm.

*Proof of Theorem 5.* Denote by $\widetilde{\mathbf{Q}}^{(0)} \stackrel{\text{def}}{=} []$ being empty matrix and for every $j \in [s]$ by

$$\widetilde{\mathbf{Q}}^{(j)} \stackrel{\text{def}}{=} \mathsf{Oja}^{++}\left(\{x_1, \ldots, x_{11jT_0}\}, \{\eta_1, \ldots, \eta_{11jT_0}\}, \{(11T_0, \mathbf{Q}^{(j)})\}_{i=1}^j\right)$$

---

[26]In fact, the intermediate learning rates for $t \in \{1, 2, \ldots, 11sT_0\}$ can all be set to the same value $\frac{C}{\rho T_0}$. We make them slightly decrease between epochs just for the sake of having a cleaner proof.

the output of Oja$^{++}$ if we run the outer loop only for $j$ iterations. By definition, it satisfies
$$\forall j \in [s]: \quad \widetilde{\mathbf{Q}}^{(j)} = \mathsf{Oja}\left(\{x_{11(j-1)T_0+b}\}_{b=1}^{T_0}, \{\eta_{11(j-1)T_0+b}\}_{b=1}^{T_0}, [\widetilde{\mathbf{Q}}^{(j-1)}, \mathbf{Q}^{(j)}]\right) ,$$
so we can view each $\widetilde{\mathbf{Q}}^{(j)}$ as the output of the old Oja's algorithm when initialized on the previous output $\widetilde{\mathbf{Q}}^{(j-1)}$ appended with a new random matrix $\mathbf{Q}^{(j)}$. We illustrate this pictorially in Figure 3.

We also introduce some notations (see Figure 4 for an illustration). For each $i \in [s]$, we define $\mathbf{V}_i$ to be the (column) orthonormal matrix consisting of all eigenvectors of $\boldsymbol{\Sigma}$ with eigenvalue $> \lambda_k - \frac{i}{s+1}\rho$, and $\mathbf{W}_i$ to be the (column) orthonormal matrix consisting of all eigenvectors of $\boldsymbol{\Sigma}$ with eigenvalue $\leq \lambda_k - \frac{i}{s+1}\rho$. Define $\mathbf{V}_0 = \mathbf{V}$ and $\mathbf{W}_0 = \mathbf{Z}$. Let $k_i \geq k$ be the column dimension of $\mathbf{V}_i$.

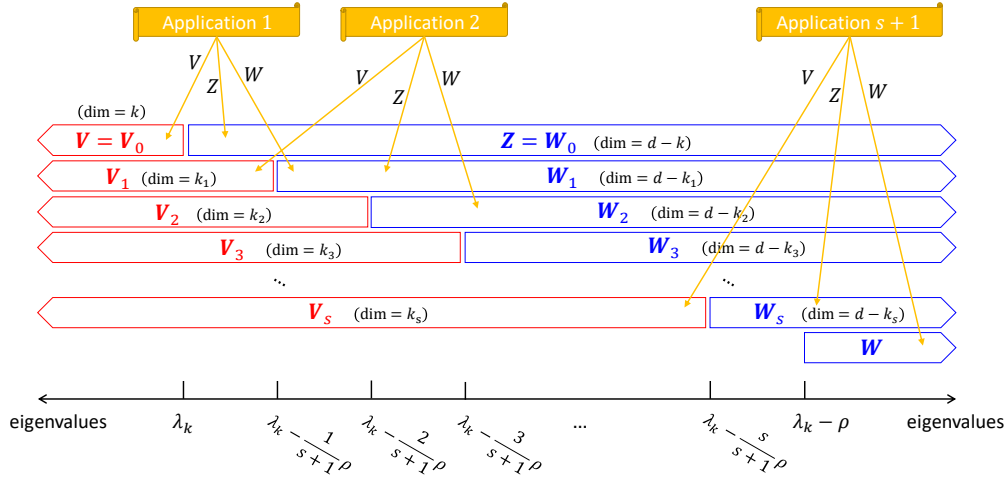

Figure 4: Notations $\mathbf{V}_0, \ldots, \mathbf{V}_s, \mathbf{W}_0, \ldots, \mathbf{W}_s$ for the proof of Theorem 5.

Below, we wish to apply our Lemma Main 6 a total of $s+1$ times each corresponding to one outer loop of Oja$^{++}$. Each time we shall use our new initialization lemma (i.e., Lemma iii.J.2) in order to satisfy the preassumption of Lemma Main 6. The first $s$ applications serve as a graduated warm-start phase and the last application provides the final convergence rate.

**Application 1 Through $s$.** Define event $\mathcal{C}_i \stackrel{\text{def}}{=} \{\forall j \in \mathbb{N}, j \leq i : \|\mathbf{W}_j^\top \widetilde{\mathbf{Q}}^{(j)}\|_F^2 \leq \frac{1}{2}\}$. In this first step we prove $\mathbf{Pr}[\mathcal{C}_i] \geq 1 - 2ip$ for all $i = 0, 1, \ldots, s$ by induction. The base case $\mathbf{Pr}[\mathcal{C}_0] = 1$ is obvious because $\widetilde{\mathbf{Q}}^{(0)}$ is an empty matrix.

Suppose $\mathbf{Pr}[\mathcal{C}_{i-1}] \geq 1 - 2(i-1)p$ holds true for some $i \in [s]$ and we wish to bound $\mathbf{Pr}[\mathcal{C}_i]$. We first note that event $\mathcal{C}_{i-1}$ implies (or if $i=1$ then $\widetilde{\mathbf{Q}}^{(i-1)}$ is an empty matrix)
$$\left(\widetilde{\mathbf{Q}}^{(i-1)}\right)^\top \mathbf{V}_{i-1} \mathbf{V}_{i-1}^\top \widetilde{\mathbf{Q}}^{(i-1)} = \mathbf{I} - \left(\widetilde{\mathbf{Q}}^{(i-1)}\right)^\top \mathbf{W}_{i-1} \mathbf{W}_{i-1}^\top \widetilde{\mathbf{Q}}^{(i-1)} \succeq \frac{1}{2}\mathbf{I} .$$
Next, we have $\widetilde{\mathbf{Q}}^{(i-1)} \in \mathbb{R}^{d\times\alpha}$ and $\mathbf{Q}^{(i)} \in \mathbb{R}^{d\times\beta}$ where $\alpha = (k - \lfloor k/2^{i-1} \rfloor)$ and $\beta = \lfloor k/2^{i-1} \rfloor - \lfloor k/2^i \rfloor$. Since $k_{i-1} - \alpha \geq k - \alpha \geq 2\beta - 1$, we have
$$\frac{\beta}{(\sqrt{k_{i-1} - \alpha} - \sqrt{\beta - 1})^2} \leq \frac{\beta}{(\sqrt{k - \alpha} - \sqrt{\beta - 1})^2} \leq \frac{\beta}{(\sqrt{2\beta - 1} - \sqrt{\beta - 1})^2} < 6 .$$
Therefore, we can apply Lemma iii.J.2 on
$$\boxed{\widetilde{\mathbf{Q}} = \widetilde{\mathbf{Q}}^{(i-1)}, \quad \mathbf{Q} = \mathbf{Q}^{(i)}, \quad (\mathbf{V}, \mathbf{Z}) = (\mathbf{V}_{i-1}, \mathbf{W}_{i-1}).}$$
and derive that, denoting by $\overline{\mathbf{Q}} = [\mathbf{Q}^{(i-1)}, \mathbf{Q}^{(i)}]$, then with probability $1 - p$ over the random choice of $\mathbf{Q}^{(i)}$, the following event $\mathcal{B}_i$ holds (for some polynomially small $q$):

$$\mathcal{B}_i \stackrel{\text{def}}{=} \begin{cases} \left\|(\mathbf{W}_{i-1}^\top \overline{\mathbf{Q}})\left(\overline{\mathbf{Q}}^\top \mathbf{V}_{i-1}\mathbf{V}_{i-1}^\top \overline{\mathbf{Q}}\right)^{-1/2}\right\|_F^2 \leq \widetilde{O}\left(\frac{d}{p^2}\right) & \text{and} \\ \Pr_{x_1,\ldots,x_T}\left[\substack{\exists \ell \in [T] \\ \exists t \in [T]} \left\|x_t^\top \mathbf{W}_{i-1}\mathbf{W}_{i-1}^\top (\mathbf{\Sigma}/\lambda_{k+1})^{\ell-1} \overline{\mathbf{Q}}\left(\overline{\mathbf{Q}}^\top \mathbf{V}_{i-1}\mathbf{V}_{i-1}^\top \overline{\mathbf{Q}}\right)^{-1/2}\right\|_2^2 \geq \widetilde{\Omega}\left(\frac{1}{p^2}\right)\right] \leq q \\ \left\|\nu_j^\top \mathbf{W}_{i-1}\mathbf{W}_{i-1}^\top \overline{\mathbf{Q}}\left(\overline{\mathbf{Q}}^\top \mathbf{V}_{i-1}\mathbf{V}_{i-1}^\top \overline{\mathbf{Q}}\right)^{-1/2}\right\|_2^2 \leq \widetilde{O}\left(\frac{1}{p^2}\right) & \text{for every } j \in [d]. \end{cases}$$

Whenever $\mathcal{B}_i$ holds, because $\widetilde{\mathbf{Q}}^{(j)}$ can be viewed as the output of the original Oja's algorithm on input $\overline{\mathbf{Q}}$ and applied with learning rates $\{\eta_{11(i-1)T_0+b}\}_{b=1}^{T_0}$, we can apply Lemma Main 6 with

$$\boxed{\begin{array}{c} \mathbf{Q} = \overline{\mathbf{Q}}, \quad (\mathbf{V}, \mathbf{Z}, \mathbf{W}) = (\mathbf{V}_{i-1}, \mathbf{W}_{i-1}, \mathbf{W}_i), \quad (\lambda, \rho) = \left(\lambda_k - \frac{i-1}{s+1}\rho, \frac{1}{s+1}\rho\right), \\ (\Xi_{\mathbf{Z}}, \Xi_x) = (\widetilde{O}(d/p^2), \widetilde{O}(1/p^2)) \end{array}}$$

We emphasize that these learning rates, $\{\eta_{11(i-1)T_0+b}\}_{b=1}^{T_0}$ —once shifted left by $11(i-1)T_0$— are exactly

$$\boxed{\eta_t = \begin{cases} \widetilde{\Theta}\left(\frac{1}{\rho T_0}\right) & t \leq T_0 \\ \widetilde{\Theta}\left(\frac{1}{\rho t}\right) & t > T_0 \end{cases}, \quad T = 11T_0, \quad \Lambda = \Lambda_1}$$

so they satisfy the assumption of Lemma Main 6 according to Parameter iii.L.1. Finally, the conclusion of Lemma Main 6 tells us

$$\Pr\left[\|\mathbf{W}_i^\top \widetilde{\mathbf{Q}}^{(i)}\|_F^2 \geq \frac{1}{2} \Big| \mathcal{C}_{i-1} \cap \mathcal{B}_i\right] \leq p \enspace,$$

and by union bound, we have (recalling $\Pr[\mathcal{B}_i] \geq 1 - p$):

$$\begin{aligned} \Pr[\mathcal{C}_i] &= \Pr\left[\mathcal{C}_{i-1} \wedge \|\mathbf{W}_i^\top \widetilde{\mathbf{Q}}^{(i)}\|_F^2 \leq \frac{1}{2}\right] \geq \Pr[\mathcal{C}_{i-1} \cap \mathcal{B}_i] \times (1-p) \\ &= (1-p)\Pr[\mathcal{C}_{i-1}]\Pr[\mathcal{B}_i \mid \mathcal{C}_{i-1}] \geq (1-p)(1-2(i-1)p)(1-p) \geq 1 - 2ip \enspace. \end{aligned}$$

This finishes the proof that $\Pr[\mathcal{C}_s] \geq 1 - 2sp$.

**Application $s+1$.** We now focus on the last outer loop of Oja$^{++}$, which satisfies

$$\mathbf{Q}_T = \mathsf{Oja}\left(\{x_t\}_{t=11sT_0+1}^T, \{\eta_t\}_{t=11sT_0+1}^T, \widetilde{\mathbf{Q}}^{(s)}\right) \enspace.$$

Similar to the first $s$ outer loops, under event $\mathcal{C}_s$ we have $\|\mathbf{W}_s^\top \widetilde{\mathbf{Q}}^{(s)}\|_F^2 \leq \frac{1}{2}$ and thus $(\widetilde{\mathbf{Q}}^{(s)})^\top \mathbf{V}_s \mathbf{V}_s^\top \widetilde{\mathbf{Q}}^{(s)} \succeq \frac{1}{2}\mathbf{I}$. This implies $\|((\widetilde{\mathbf{Q}}^{(s)})^\top \mathbf{V}_s \mathbf{V}_s^\top \widetilde{\mathbf{Q}}^{(s)})^{-1/2}\|_2 \leq \sqrt{2}$ and therefore

$$\|\mathbf{W}_s^\top \widetilde{\mathbf{Q}}^{(s)}((\widetilde{\mathbf{Q}}^{(s)})^\top \mathbf{V}_s \mathbf{V}_s^\top \widetilde{\mathbf{Q}}^{(s)})^{-1/2}\|_F^2 \leq \|\mathbf{W}_s^\top \widetilde{\mathbf{Q}}^{(s)}\|_F^2 \|((\widetilde{\mathbf{Q}}^{(s)})^\top \mathbf{V}_s \mathbf{V}_s^\top \widetilde{\mathbf{Q}}^{(s)})^{-1/2}\|_2^2 \leq 2.$$

Thus, for every $x \in \mathbb{R}^d$ with $\|x\|_2 \leq 1$:

$$\|x^\top \mathbf{W}_s \mathbf{W}_s^\top \widetilde{\mathbf{Q}}^{(s)}((\widetilde{\mathbf{Q}}^{(s)})^\top \mathbf{V}_s \mathbf{V}_s^\top \widetilde{\mathbf{Q}}^{(s)})^{-1/2}\|_2 \leq \|\mathbf{W}_s^\top \widetilde{\mathbf{Q}}^{(s)}\|_2 \|((\widetilde{\mathbf{Q}}^{(s)})^\top \mathbf{V}_s \mathbf{V}_s^\top \widetilde{\mathbf{Q}}^{(s)})^{-1/2}\|_2 \leq 2 \enspace.$$

Now we can apply Lemma Main 6 again with parameters

$$\boxed{\mathbf{Q} = \widetilde{\mathbf{Q}}^{(s)}, \quad (\mathbf{V}, \mathbf{Z}, \mathbf{W}) = (\mathbf{V}_s, \mathbf{W}_s, \mathbf{W}), \quad (\lambda, \rho) = \left(\lambda_k - \frac{s}{s+1}\rho, \frac{1}{s+1}\rho\right), \quad \Xi_{\mathbf{Z}} = \Xi_x = 2}$$

We emphasize that these learning rates, $\{\eta_t\}_{t=11sT_0+1}^T$ —once shifted left by $11sT_0$— are exactly

$$\boxed{\eta_t = \begin{cases} \widetilde{\Theta}\left(\frac{1}{\rho T_1}\right) & t \leq T_1 \\ \widetilde{\Theta}\left(\frac{1}{\rho t}\right) & t > T_1 \end{cases}, \quad T = T - 11sT_0, \quad \Lambda = \Lambda_2}$$

so they satisfy the assumption of Lemma Main 6 according to Parameter iii.L.1. We can conclude from Lemma Main 6 that
$$\mathbf{Pr}\left[\|\mathbf{W}^\top \mathbf{Q}_T\|_F^2 = \widetilde{\Omega}\left(\frac{T_1}{T - 11sT_0}\right) \;\Big|\; \|\mathbf{W}_s^\top \widetilde{\mathbf{Q}}^{(s)}\|_F^2 \leq \frac{1}{2}\right] \leq p \;.$$

Taking into account $\mathbf{Pr}\left[\|\mathbf{W}_s^\top \widetilde{\mathbf{Q}}^{(s)}\|_F^2 \geq \frac{1}{2}\right] \leq 2sp$ we complete the proof of Theorem 5 if we replace $p$ with $p/(2s+1)$. □